\documentclass[12pt]{article}    
\pdfoutput=1
\usepackage{amscd,amssymb,stmaryrd}
\usepackage{amsmath}
\usepackage{amscd,amssymb,stmaryrd}
\usepackage{graphicx}
\usepackage{subfigure,fancybox}
\usepackage{amscd,amstext}
\usepackage{amsmath}
\usepackage{mathtools}
\usepackage{setspace}
\usepackage{amssymb}
\usepackage{slashed}
\usepackage{pstricks,pst-node,pst-plot,pst-coil}
\usepackage{graphicx}
\usepackage{ntheorem}
\usepackage{upgreek}
\usepackage{enumerate}
\usepackage{graphicx}
\usepackage{cite}
\usepackage{tikz}
\usepackage{authblk}
\usepackage{txfonts}
\usepackage[colorlinks,linkcolor=black,anchorcolor=black,citecolor=black,CJKbookmarks=True]{hyperref}
\theorembodyfont{\normalfont}
\newtheorem{thm}{Theorem}[section]
\newtheorem{lem}{Lemma}[section]
\newtheorem*{pf}{PROOF:}
\newtheorem{prop}{Proposition}[section]

\newtheorem{remark}{Remark}[section]

\numberwithin{equation}{section}
\makeatletter

\newcommand{\Rmnum}[1]{\expandafter\@slowromancap\romannumeral #1@}
\makeatother

\title{\large Shifted shock formation for the 3D compressible Euler equations with damping and variation of the vorticity}
\author[1]{\normalsize Chen Zhendong\footnote{Acknowledgement:
This is part of the Ph. D thesis of the author written under the supervision of Professor Zhouping Xin at the Institute of Mathematical Sciences of The Chinese University of Hong Kong. }}

\affil[1]{Institute of Mathematical Science, The Chinese University of Hong Kong, Shatin, NT, Hong Kong.}
\date{}
\begin{document}           
\def\ra{\rightarrow}\def\Ra{\Rightarrow}\def\lra{\Longleftrightarrow}\def\fai{\varphi}
\def\s{\slashed}\def\Ga{\Gamma}\def\til{\tilde}
\def\de{\Delta}\def\fe{\phi}\def\sg{\slashed{g}}\def\m{\mathbb}
\def\a{\alpha}\def\p{\rho}\def\ga{\gamma}\def\lam{\lambda}\def\ta{\theta}\def\sna{\slashed{\nabla}}
\def\pa{\partial}\def\be{\beta}\def\da{\delta}\def\ep{\varepsilon}\def\dc{\underline{\chi}}
\def\si{\sigma}\def\Si{\Sigma}\def\wi{\widetilde}\def\wih{\widehat}\def\beeq{\begin{eqnarray*}}
\def\eeq{\end{eqnarray*}}\def\na{\nabla}\def\lie{{\mathcal{L}\mkern-9mu/}}\def\rie{\mathcal{R}}
\def\ud{\underline}\def\les{\lesssim}\def\ka{\kappa}
\def\dl{\underline{L}}\def\du{\underline{u}}
\def\hs{\hspace*{0.5cm}}\def\bee{\begin{equation}}\def\ee{\end{equation}}\def\been{\begin{enumerate}}
\def\een{\end{enumerate}}\def\bes{\begin{split}}\def\zpi{\prescript{(Z)}{}{\pi}}\def\szpi{\prescript{(Z)}{}{\slashed{\pi}}}
\def\ees{\end{split}}\def\hra{\hookrightarrow}\def\udx{\ud{\xi}_{A}}\def\ude{\ud{\eta}_{A}}
\def\tpi{\prescript{(T)}{}{\pi}}\def\lpi{\prescript{(L)}{}{\pi}}\def\qpi{\prescript{(Q)}{}{\pi}}
\def\stpi{\prescript{(T)}{}{\slashed{\pi}}}\def\sqpi{\prescript{(Q)}{}{\slashed{\pi}}}
\def\zgpi{\prescript{(Z)}{}{\wi{\pi}}}\def\szgpi{\prescript{(Z)}{}{\wi{\slashed{\pi}}}}
\def\tgpi{\prescript{(T)}{}{\wi{\pi}}}\def\lgpi{\prescript{(L)}{}{\wi{\pi}}}\def\qgpi{\prescript{(Q)}{}{\wi{\pi}}}
\def\stgpi{\prescript{(T)}{}{\wi{\slashed{\pi}}}}\def\sqgpi{\prescript{(Q)}{}{\wi{\slashed{\pi}}}}
\def\pre #1 #2 #3{\prescript{#1}{#2}{#3}}\def\spi{\slashed{\pi}}\def\rgpi{\prescript{(R_{i})}{}{\wi{\pi}}}
\def\rpi{\prescript{(R_{i})}{}{\pi}}\def\srpi{\prescript{(R_{i})}{}{\s{\pi}}}
\def\pp{\uprho}\def\vw{\upvarpi}\def\srgpi{\pre {(R_{i})} {} {\wi{\s{\pi}}}}
\def\supnormda{L^{\infty}(\Si_{t}^{\da})}\def\supnormu{L^{\infty}(\Si_{t}^{u})}
\def\normda2{L^{2}(\Si_{t}^{\da})}\def\normu2{L^{2}(\Si_{t}^{u})}
\def\les{\lesssim}\def\sdfai{|\s{d}\fai|}\def\sd{\s{d}}
\def\oli{\overline}\def\vta{\vartheta}\def\ol{\overline}
\def\ea+2{E_{\leq|\a|+2}}\def\dea+2{\ud{E}_{\leq|\a|+2}}\def\fa+2{F_{\leq|\a|+2}}\def\dfa+2{\ud{F}_{\leq|\a|+2}}
\def\tea+2{\wi{E}_{\leq|\a|+2}}\def\tdea+2{\ud{\wi{E}}_{\leq|\a|+2}}\def\tfa+2{\wi{F}_{\leq|\a|+2}}\def\tdfa+2{\ud{\wi{F}}_{\leq|\a|+2}}
\def\ntea{\wi{E}_{\leq|\a|+1}}\def\ntdea{\ud{\wi{E}}_{\leq|\a|+1}}\def\ntfa{\wi{F}_{\leq|\a|+1}}\def\ntdfa{\ud{\wi{F}}_{\leq|\a|+1}}
\def\ba2{b_{|\a|+2}}\def\Ka{K_{\leq|\a|+2}}\def\mumba{\mu_m^{-2\ba2}}\def\nba{b_{|\a|+1}}
\def\endpf{\hfill\raisebox{-0.3cm}{\rule{0.1mm}{3mm}\rule{0.5mm}{0.1mm}\rule[-0.6mm]{2.3mm}{0.7mm}\raisebox{3mm}[0pt][0pt]{\makebox[0pt][r]{\rule{2.9mm}{0.1mm}}}\rule[-0.6mm]{0.1mm}{3.5mm}\rule[-0.6mm]{0.5mm}{3mm}} \vspace{0.5cm}\\}\def\mbr{\mathbb{R}}\def\mr{\mathring}\def\mcl{\mathcal{L}}

\pagestyle{myheadings} \thispagestyle{empty} \markright{}

\maketitle
\begin{abstract}
\quad In this paper, we consider the shock formation problem for the 3-dimensional(3D) compressible Euler equations with damping inspired by the work \cite{BSV3Dfulleuler}. It will be shown that for a class of large data the damping can not prevent the formation of point shock and the damping effect shifts the shock time and the wave amplitude while the shock location and the blow up direction remain same with the information of this point shock being computed explicitly. Moreover, the vorticity is concentrated in the non-blow up direction which varies exponentially due to the damping effect. Our proof is based on the estimates for the modulated self-similar variables and lower bounds for the Lagrangian trajectories.
\end{abstract}

\tableofcontents
\section{\textbf{Introduction}}
We consider the following 3D Euler system with damping
\begin{equation}\label{euler}
\bes
\pa_t\p+u\cdot\nabla_x\p&=-\p div_x u\\
\pa_tu+u\cdot\nabla_x u&=-\dfrac{1}{\p}\nabla_x p-a u\\
\pa_tk+u\cdot\nabla_x k&=0,
\end{split}
\end{equation}
with the space-time variables $(t,x)=(t,x_1,x_2,x_3)\in \mbr\times \mbr^3$, the velocity $u: \mbr\times \mbr^3\ra \mbr^3$, the positive density $\p: \mbr^3\times\mbr\ra \mbr_+$, the entropy $k :\mbr\times\mbr^3\ra \mbr$ and the pressure $p=p(\p,k): \mbr\times\mbr^3\ra \mbr_+$ is given by
\[
p(\p,k)=\dfrac{1}{\ga}\p^{\ga}e^k.
\]
with the adiabatic constant $\ga>1$. These equations describe the motion of a perfect fluid which follow from the conservation of mass, momentum and energy respectively. The term $-au$ on the right hand side of momentum equation represents the damping effect to the fluid with the damping constant $a\neq0$.\\
\hs Instead of using the density, it's convenient to rewrite the equation by the sound speed, which is defined as
\[
c=\sqrt{\dfrac{\pa p}{\pa \p}}=\p^{\a}e^{\frac{k}{2}},
\]
where $\a=\dfrac{\ga-1}{2}.$ For simplicity, we denote $\sigma$ to be the re-scaled sound speed
\[
\sigma=\dfrac{1}{\a}c=\dfrac{1}{\a}\p^{\a}e^{\frac{k}{2}}.
\]
Then, \eqref{euler} can be rewritten as
\begin{equation}\label{euler2}
\bes
\pa_t\si+u\cdot\nabla_x\si+\a\si div_xu&=0,\\
\pa_tu+u\cdot\nabla_xu+\a\si div_x\si&=\dfrac{\a}{2\ga}\si^2\nabla_xk-au,\\
\pa_t k+u\nabla_xk&=0.
\end{split}
\end{equation}
 Define the specific vorticity $\zeta$ as $\zeta=\dfrac{\omega}{\p}=\dfrac{curl\ u}{\p}$. Then, $\zeta$ satisfies the following equation:
\bee\label{equationzeta}
\pa_t\zeta+u\cdot\nabla_x\zeta=\zeta\cdot\nabla_x u+\underbrace{\dfrac{\a}{\ga}\dfrac{\si}{\p}\nabla_x\si\times\nabla_xk}_{=\frac
{1}{\p^3}\nabla_x\p\times\nabla_x p}-a\zeta.
\ee
The first and second terms on the right-hand side of the equation\eqref{equationzeta} represent the consequences of vortex stretching and the baroclinic torque. The third term arises from the intersection of acoustic and entropy waves and results from the non-aligned pressure and density gradients. In the scenario where the system is isentropic, such that $\nabla_x\p\parallel\nabla_x p$, this term vanishes. However, when the system is non-isentropic, this term plays a crucial role in the generation of vorticity. Irrotational data can produce vorticity instantaneously, as demonstrated in \cite{BSV3Dfulleuler}. For example, if the density gradient is perpendicular to the pressure gradient, the lighter gas will be accelerated faster than the denser gas, resulting in the creation of vorticity.\\
\hs It is well-known that damping can prevent shock formation when the energy of the initial  data is small. However, for large data sets, such as short pulse data or data with large energy, we will demonstrate that damping may not sufficiently prevent the formation of shock. Furthermore, we will show that the damping effect shift the time at which a shock occurs, while the blow-up location remains unchanged compared with the undamped case, with the degree of shift reliant on the value of $a$.
\subsection{\textbf{Review of prior results for the Euler system}}
\hs For the system\eqref{euler}, there are many studies both in $1D$ and multi-dimensions. For the one-dimensional Euler equations with damping, the global existence of smooth solutions with small data was proved by Nishida\cite{1978NonlinearHE} and Slemrod\cite{slemrod} showed that for small data ($L^{\infty}$ sense), the $1D$ Euler equations with damping admit a global smooth solution while for large data, the equations can develop a shock in finite time. Later, these results were generalized by many authors, see \cite{liu1992,NISHIHARA2000191,MARCATI1990,MARCATI2000,MARCATI2005,NISHIHARA2000} and the references therein.\\
\hs In multi-dimensional case, the global existence and $L^2$ estimates for the solutions to the $3D$ isentropic Euler system with damping was obtained by Kawashima\cite{Kawashima1984}. Later, Sideris-Thomases-Wang in\cite{sideris2003long} showed that the size of the smooth initial data plays a key role for the lifespan to the $3D$ isentropic Euler equations with damping $a\p u$. If the initial data are small, then damping can prevent the development of singularities; while if the initial data are large, the damping is not strong enough to prevent the formation of singularities in finite time. However, the authors in\cite{sideris2003long} only obtained the finite lifespan of the solutions and did not show the shock formation (including the shock mechanism) in finite time. The long time behavior of the solutions was obtained and generalized by many authors, see \cite{WANG2001410,Kawashima2004,PAN2009581,TAN20121546,sideris2003long} and the references therein. \\
\hs In \cite{christodoulou2014compressible}, Christodoulou achieved a significant breakthrough in understanding the shock mechanism for hyperbolic systems in three dimensions. Specifically, he considered the classical, non-relativistic, compressible Euler's equations in three spatial dimensions, taking the data to be irrotational and isentropic. By imposing certain conditions on the initial data (i.e., the short pulse assumption), he obtained a complete geometric description of the maximal classical development. Notably, he provided a detailed analysis of the behavior of the solution at the boundary of the domain of the maximal classical solution, including a comprehensive description of its geometry. This work represents a significant advancement in understanding of shock formation in hyperbolic systems in three dimensions.\\
\hs In \cite{Ontheformationofshocks}, Yu-Miao applied Christodoulou's framework to the quasilinear wave equation in three spatial dimensions. By constructing a family of short pulse initial data that was first introduced by Christodoulou in \cite{christodoulou2007}, they demonstrated how the solution breaks down near the singularity and gave a sufficient condition on the initial data which leads to the shock formation in finite time. Additionally, J.Speck and J.Luk applied Christodoulou's framework to the 2D Euler system without the assumption of irrotation. In \cite{luk2018shock}, they studied plane-symmetric initial data with short pulse perturbation. For such initial data, they showed the shock mechanism such that the first derivatives of $u$ and $\rho$ blow up while $u$ and $\rho$ remain bounded near the shock. Then, in \cite{Specklukshockformation}, they generalized these results to the 3D case. Specifically, they considered the 1D Euler equation of a simple small-amplitude solution as a plane-symmetric solution in 3D. They perturbed this solution in the $(x_2,x_3)$ directions as a nearly plane-symmetric initial data for 3D isentropic Euler equations. They proved that the shock formation mechanism is stable under small and compactly supported perturbations with non-trivial vorticity and provided a precise description of the first singularity.\\
\hs Recently, Buckmaster, Shkoller, and Vicol published several results regarding shock formation in the multi-dimensional Euler system. Specifically, in their paper \cite{BSV2Disentropiceuler}, they studied the 2D isentropic Euler equations under azimuthal symmetry (which differs from the 1D problem) with smooth initial data of finite energy and nontrivial vorticity. By utilizing modulated self-similar variables, they were able to obtain point shock forms in finite time, with explicit computation of the blow-up time and location. Furthermore, the solutions near the shock exhibit a cusp type.\\
\hs In a subsequent paper, \cite{BSV3Disentropiceuler}, Buckmaster, Shkoller, and Vicol generalized their results to the 3D isentropic case for the ideal gas without any symmetry assumptions. In addition to the findings reported in \cite{BSV2Disentropiceuler}, they were able to demonstrate the precise direction of blow-up and the geometric structure of the tangent surface of the shock profile. They also provided homogeneous Sobolev bounds for the fluid variables with an initial datum of large energy. Later, in\cite{BSV3Dfulleuler}, the authors extended their results to the full Euler equations. In this work, they primarily investigated the evolution and creation of vorticity, demonstrating that the vorticity remains bounded up to the shock formation. Notably, one of the primary differences between the isentropic and non-isentropic cases is the baroclinic torque in the equation of vorticity, which arises from the interaction between sound waves and entropy waves. Additionally, they constructed a set of irrotational data that results in the instantaneous creation of vorticity, which remains non-zero up to the shock.\\
\hs It is important to note that the point shock in the three aforementioned results is stable. That is, for any small, smooth, and generic perturbation of the given initial data, the corresponding Euler system results in a smooth solution that blows up in a small neighborhood of the original shock time and location. However, in \cite{formationofunstableshock}, the authors demonstrated the existence of an open set of initial data that leads to the formation of an unstable shock. The primary difference between stable and unstable shocks is the set of background solutions for the self-similar variable $W$, specifically the solutions of the various self-similar Burgers equations \eqref{burgers3}. For further details, one may refer to \cite{Eggersselfsimilar}.\\
\hs The major tool used in the works mentioned above is the method of self-similar coordinates. This method was first introduced by Y.Giga and R.Kohn in \cite{GigaKohnselfsimilar} to study the asymptotic behavior of the solution to the nonlinear heat equation near the point of singularity. Y.Giga and R.Kohn considered the following nonlinear heat equation
 \bee\label{heat}
 u_t-\de u-|u|^{p-1}u=0,
 \ee
 where $p>1$ and $(x,t)\in\mathbb{R}^n\times(-1,0)$. They aimed to show the behavior of the solution near the singularity. To this end, they proposed the following self-similar transformation\footnote{Note that the transformation degenerates at $t=0$.}:
 \bee
 y=e^{\frac{1}{2}s}x,\hs s=-\ln(-t),\hs w(y,s)=e^{-\frac{1}{2(p-1)}s}u(x,t),
 \ee
 which transforms \eqref{heat} into
 \bee\label{heat2}
 w_s-\de w+\frac{1}{2}y\cdot\nabla_y w+\frac{1}{p-1}w-|w|^{p-1}w=0.
 \ee
 This transformation is motivated by the scaling property of \eqref{heat}. That is, for any $u$ solves \eqref{heat}, $u_{\lam}:=\lam^{\frac{2}{p-1}}u(\lam x,\lam^2t)$ solves \eqref{heat} as well. Based on the analysis of \eqref{heat2}, they were able to demonstrate the asymptotic behavior of $u$ near the blow-up point $(0,0)$. Later, such method has been applied to various other dynamical systems, including the Sch$\ddot{o}$dinger equation \cite{Merleschrodinger}, the Prandtl equations \cite{Prandtlselfsimilar}, the transverse viscous Burgers equation \cite{Masmoudiselfsimilarburgers}, and the semilinear wave equation \cite{Merleselfsimilar}. The method of self-similar variables can provide precise information about the singularity of a given system by adding modulation variables to enforce certain orthogonality conditions and track the position of the singularity. It can also be used to characterize the type of singularities present, such as stable fixed points or center manifold and etc..., based on the behavior of the solutions near the singularity. For example, the stable shock in \cite{BSV2Disentropiceuler}, \cite{BSV3Disentropiceuler}, and \cite{BSV3Dfulleuler} is based on the fact that the solutions approach the background solution near the shock exponentially with respect to the self-similar variables. For further details, one may refer to \cite{Eggersselfsimilar}.
\subsection{\textbf{Outline}}
\hs In section \ref{section2}, the fundamental ideas of \cite{BSV3Disentropiceuler} and \cite{BSV3Dfulleuler} will be illustrated by applying them to the standard Burgers equation and the Burgers equation with damping. Our results will demonstrate that the modulated variables can accurately track the location and time of the shock, and that the damping effect only shifts the blow-up time while leaving the blow-up location unchanged compared with the undamped case.\\
\hs Section \ref{section3} is devoted to reformulating of the Euler system \eqref{euler} into its modulated self-similar version \eqref{WZAK} by adapting the framework of \cite{BSV3Dfulleuler}. Additionally, we extend the solution of the self-similar Burgers equation, which is introduced in Section \ref{section2}, to three dimensions and derive higher-order derivatives of the evolution equations for the self-similar variables $(W,Z,A,K)$.\\
\hs In Section \ref{section4}, we present an explicit construction of the initial datum in \cite{BSV3Dfulleuler}. Note that the initial data for variable $W$ is dependent on the background solution $\bar{W}$, which is a refined solution of the self-similar Burgers equation. Furthermore, we outline the self-similar Bootstrap assumptions for various variables and their derivatives, including the modulation variables and the self-similar variables $(W, Z, A, K)$. These assumptions are more restrictive than the initial data constructed earlier and will be recovered in the subsequent sections. Finally, we state the main theorem (Theorem \ref{mainthm}) of the paper at the end of this section.\\
\hs In Section \ref{section5}, we present multiple estimates based on bootstrap assumptions, including the estimation of damping and forcing terms. Additionally, we state the $\dot{H}^m$ energy bounds for $(W,Z,A,K)$, whose proof solely relies on the bootstrap assumptions and the Friedrich's energy estimates for the hyperbolic system. These estimates lead to higher-order derivatives estimates for $(W,Z,A,K)$, which in turn help to further refine the estimates for the forcing terms.\\
\hs Section \ref{section6} is a key aspect of our work. In this section, we recover the bootstrap assumptions for the modulation variables under 10 constraints regarding $\pa^{\ga}W(0,s)$, where $|\ga|\leq 2$. These constraints are initially satisfied and remain valid over time as long as the estimates for the modulation variables hold. In addition, it will be clearly shown that how the damping effect can influence the information of shock within this section.\\
\hs In Section \ref{section7}, we derive lower bounds for the Lagrangian trajectories and introduced the crucial lemma\ref{keylemma}, which provides the key estimate required for recovering the remaining bootstrap assumptions. Moreover, we investigate the evolution of the vorticity, which is concentrated on the non-blow up direction, and demonstrated how it is affected by the damping effect. Inspired the work in\cite{BSV3Dfulleuler}, it will be shown that the initial region which genuinely affects the shock formation.\\
\hs By following the estimates presented in sections \ref{section5} to \ref{section7}, we are able to recover the bootstrap assumptions for $(W, Z, A, K)$ as outlined in section \ref{section8}. This recovery is based on the significant lemma \ref{keylemma} and the utilization of the weighted framework \eqref{framework}. Note that to recover the bootstrap assumption for $\pa_1A_{\nu}$, one has to establish the relation between $\pa_1A_{\nu}$ and the specific vorticity $\Omega_{\nu}$.\\
\hs Section \ref{section9} serves to establish the $\dot{H}^m$ energy bounds for the self-similar variables, thereby concluding the proof for the main theorem \ref{mainthm}. Instead of studying the system for $(W, Z, A, K)$, we utilize various Sobolev inequalities to derive the energy bounds for the equations of the velocity $U$, the pressure $P$, and the entropy $H$ (as seen in \eqref{UPHsystem}). This approach proves to be more convenient to apply the Friedrich's energy estimates.
\subsection{\textbf{Notations}}
Through the whole paper, the following notations will be used unless stated otherwise.
\begin{itemize}
\item Latin indices $\{i,j,k,l,\cdots\}$ take the values $1, 2, 3, $ and Greek indices $\{\a,\be,\ga,\cdots\}$
take the values $2,3$. Repeated indices are meant to be summed.
 \item For a three-component vector $v$, denote the last two components of $v$ simply as $\check{v}$. For example, one can rewrite the gradient operator as $\nabla=(\pa_1,\check{\nabla}).$
 \item The convention $f\les h$ means that there exists a universal positive constant $C$ such that $f\leq Ch$. The convention $A\sim B$ means that there exists a universal positive constant such that $A=O(B)$.
 \item For any function $A(y,s)$, denote $A^0$ to be $A^0=A^0(s)=A(y,s)|_{y=0}$.
\end{itemize}
\section{\textbf{Shock formation for the Burger's equation in self-similar fancy}}\label{section2}
We consider the following Cauchy problem for 1D Burger's equation
\bee\label{burgers}
\left\{\bes
&\pa_tu+u\pa_xu=0,\\
&u(x,t=-1)=f(x).
\end{split}\right.
\ee
For simplicity, we assume additionally that $f(0)=0$ and $\min \pa_xf(x)=\pa_xf(0)=-1$\footnote{Indeed, once can assume generally that $\pa_xf(0)=-c<0$. Then the framework here also applies by slightly modified, i.e. $T_{\ast}=-1+\frac{1}{c}$ and $x_{\ast}=0$, and the corresponding self-similar transformation becomes
\bee
s=-\ln(-1+\frac{1}{c}-t),\quad y=xe^{\frac{3}{2}s}.
\ee}. Then by standard characteristic method, the solution of \eqref{burgers} will form a shock at time $T_{\ast}=0$ and the location $x_{\ast}=0$ with $\pa_xu(0,t)\ra-\infty$ as $t\to 0$. However, this method has some drawbacks in the following sense:
\begin{itemize}
\item Is there any singularity before $t=0$ and whether $\pa_x u$ is the first quantity which blows up or not?
\item How does the shock profile look like, and how does the solution behave near the singularity?
\end{itemize}
 The above questions can be solved by using self-similar coordinates which is defined as follows:
\bee
s=-\ln(-t),\quad y=xe^{\frac{3}{2}s},
\ee
and the corresponding unknown is given as
\bee
u(x,t)=e^{-\frac{s}{2}}U(y,s).
\ee
Then, \eqref{burgers} is transformed as
\bee\label{selfsimilarburgers}
(\pa_s-\dfrac{1}{2})U+(\dfrac{3}{2}y+U)\pa_yU=0.
\ee
\begin{remark}
In general, the self-similar transformation should be
\bee\label{selfsimilargeneral}
s=-\ln(\tau(t)-t),y=(x-\xi(t))e^{\a s}, u(x,t)=e^{-\be s}U(y,s),
\ee
 where the parameters $\tau(t)$ and $\xi(t)$ represent the shock time and location, respectively. Here, one already knows the blow-up point is $(0,0)$, which means one could take $\tau(t)=\xi(t)=0$. Then, substituting \eqref{selfsimilargeneral} into \eqref{selfsimilarburgers} yields
\bee\label{burgers2}
(\pa_s-\be)U+[\a y+U]e^{(\a-\be-1)s}\pa_yU=0.
\ee
In order to guarantee the global existence of \eqref{burgers2}, one has to choose
\bee
\a=\be+1.
\ee
The choice of $\be=\dfrac{1}{2}$ is used to guarantee the stability of shock.\footnote{This means the solution of \eqref{burgers} is approaching exponentially to the solution of the self-similar Burger's equation (see\eqref{ssburgers}) in self-similar variable $s$, i.e.
\bee
\lim_{t\to0}|u(x,t)-(-t)^{\frac{1}{2}}\bar{U}(\frac{x}{(-t)^{\frac{3}{2}}})|=0,\quad\forall\ x\in \mbr,
\ee
where $\bar{U}$ is given by \eqref{solssburgers}. See\cite{Eggersselfsimilar} for more details.}
\end{remark}
Note that the Jacobian of the transformation is given by
\bee
\dfrac{\pa(y,s)}{\pa(x,t)}=\left|\begin{array}{cc}
e^{\frac{3}{2}s}&\frac{3}{2}ye^s\\
0&e^s
\end{array}\right|=e^{\frac{5}{2}s}=\dfrac{1}{(-t)^{\frac{5}{2}}}.
\ee
Hence, the self-similar transformation degenerates as $t\to 0$. Later we will prove that \eqref{burgers2} admits a global solution on $[0,\infty)$, and therefore, the only possibility of singularity formation is the transformation between the Cartesian coordinates and the self-similar coordinates becomes degenerated.\\
\hs Moreover, one can show that $U(y,s)$ converges to $\bar{U}$ pointwisely as $s\to\infty$. That is
\bee
\lim_{s\to\infty}|U(y,s)-\bar{U}(y)|=0,\quad \forall y\in R,
\ee
where $\bar{U}$ is the solution of the following self-similar Burger's equation
\bee\label{ssburgers}
-\dfrac{1}{2}\bar{U}+\left(\dfrac{3}{2}y+\bar{U}\right)\pa_y\bar{U}=0.
\ee
Indeed, $\bar{U}$ can be solved as a implicit function as
\bee\label{solssburgers}
y=-\bar{U}-\bar{U}^3,\hs y\in(-\infty,\infty),
\ee
which is globally defined.\\
\hs Hence, $u_x$ blows up only at the shock point and location, i.e.
\bee
\lim_{t\to0}\pa_xu(0,t)=\lim_{t\to 0}e^s\pa_yU(0,s)=\lim_{t\to 0}-\dfrac{1}{t}\to-\infty\footnotemark,
\ee
\footnotetext{The last equality comes from $\pa_yU(0,s)=-1$, which can be derived by differentiating \eqref{selfsimilarburgers} with respect to (w.r.t) $y$ and then evaluating at $y=0$. Indeed, one can show that $U(0,s)=0,\pa_yU(0,s)=-1$ hold for all $s$.}
and all the other quantities are bounded\footnote{Obviously, $|u(x,t)|\les1$ and it can be shown that for any $z\neq 0$, $\lim_{t\to0}|\pa_xu(z,t)|\les\frac{1}{z^{\frac{2}{3}}}.$}. The behavior of the solution near shock is given as
\bee
u(x,t)=(-t)^{\frac{1}{2}}\bar{U}\left(\frac{x}{(-t)^{\frac{3}{2}}}\right).
\ee
\begin{remark}
To obtain the global existence of \eqref{selfsimilarburgers}, one can check that if one defines the characteristics as
\[
\dfrac{dY}{ds}=U(Y,s)+\dfrac{3}{2}Y,\hs Y(0)=y_0,
\]
along which $U\circ Y=e^{\frac{s}{2}}U_0(y_0)$, then $e^{-\frac{3s}{2}}Y=y_0+(1-e^{-s})U_0(y_0)$, which means $U$ can be solved implicitly as $U(y,s)=e^{\frac{s}{2}}U_0\left(
e^{-\frac{3s}{2}}y-e^{-\frac{s}{2}}(1-e^{-s})U\right)$. In order to require $U$ to be defined globally, it suffices to show
\bee
1+(1-e^{-s})U'_0\neq 0.
\ee
This is guaranteed by assumption on the initial data.
\end{remark}
\subsection{\textbf{The geometric structure of shock front}}
\hs Consider the Surface $\Ga: (x,t,u(x,t))$ in $R^3$. Then, the normal vector $N$ of $\Ga$ and the Gauss curvature $K$ of $\Ga$ at each point are given as follows
\bee
N=J^{-1}(-u_x,-u_t,1),\hs K=\dfrac{-u_x^4}{(1+u_t^2+u_x^2)^2},
\ee
where $J=\sqrt{1+u_x^2+u_t^2}$. Initially, $N_0=\dfrac{1}{\sqrt{2}}(1,0,1)$, $K_0=-\dfrac{1}{4}$.
In the self-similar coordinate, consider the evolutions of $N,K$ at $y=0$,
\begin{align}
N&=\dfrac{1}{\sqrt{1+e^{2s}U_y^2+e^sU^2U_y^2}}\left(-e^{s}U_y,-e^{\frac{s}{2}}UU_y,1\right)\left|_{y=0}\right.,\\
K&=-\dfrac{1}{(1+e^{-2s}U_y^{-2}+e^{-s}U^2)^2}\left|_{y=0}\right..
\end{align}
Then, as $s\to\infty$ ($t\to 0$), since $U_y(0,s)\to-1$, $U(0,s)\to\bar{U}(0)=0$, it holds that
\bee
N\to (-1,0,0),\hs K\to -1.
\ee
Therefore, shock formation to the Burgers' equation \eqref{burgers} is equivalent to that
\begin{itemize}
   \item the transformation between the physical variables (the Cartesian coordinates) and the self-similar coordinates becomes degenerate;
   \item the normal $N(t)$ of the shock front become horizontal at shock point. 
   \end{itemize}
\subsection{\textbf{Shifted singularity for the Burgers' equation with damping}}
We consider the following Cauchy problem for 1D Burgers' equation with damping:
\bee\label{burgersdamp}
\left\{\bes
&\pa_tu+u\pa_xu=-au,\\
&u(x,t=-1)=f(x),
\end{split}\right.
\ee
where $a\neq 0$ is the damping constant and $f(x)$ is the same as in \eqref{burgers}. It will be shown that when $a<1$, the damping effect is small and \eqref{burgersdamp} behaves like the standard Burgers' equation such that smooth data can form a shock in finite time. Moreover, the damping will shift the blow-up time according to the sign of $a$ as follows.
\begin{itemize}
\item If $0<a<1=-\min\frac{1}{\pa_xf}$, the damping will delay the formation of shock and the closer $a$ is to $1$, the larger the blow-up time is;
\item if $a<0$, then the blow-up time will be in advanced.
\end{itemize}
Precisely, the blow-up time can be clarified as $T_{\ast}=-\frac{1}{a}\ln(1-a)-1$\footnote{If one assumes generally that $\pa_xf(0)=-c$, then the first case becomes that $0<a<\frac{1}{c}$ while $T_{\ast}=-\frac{1}{a}\ln(1-\frac{a}{c})-1$.} and in both case, the blow-up location remains the same compared with the undamped case (i.e. $x_{\ast}=0$). When $a\geq1$, the damping effect is strong enough and \eqref{burgersdamp} admits a global solution on $[-1,\infty)$. To this end, define the following self-similar coordinates\footnote{The choice of modulation variables $\tau(t)$ and $\xi(t)$ is mainly dependent on the invariance of \eqref{burgersdamp} under space translation and time translation. That is, given any $x_0,t_0\in\mathbb{R}$, then for any solution $u$ to \eqref{burgersdamp}, $\tilde{u}:=u(x-x_0,t-t_0)$ solves \eqref{burgersdamp} as well.}:
\begin{equation}
\bes
s&=-\ln(\tau(t)-t),\\
y&=e^{\frac{3}{2}s}(x-\xi(t)),\\
u(x,t)&=e^{-\frac{s}{2}}U(y,s),
\end{split}
\end{equation}
where $\tau(t)$ and $\xi(t)$ represent the blow-up time and location respectively, with the initial condition $\tau(-1)=\xi(-1)=0$. Note that the Jacobian of the transformation is given by
\bee
\dfrac{\pa (y,s)}{\pa (x,t)}=\left|\begin{array}{cc}
e^{\frac{3}{2}s}& \ast\\
0 &(1-\dot{\tau})e^s
\end{array}\right|=(1-\dot{\tau}(t))\dfrac{1}{(\tau(t)-t)^{\frac{5}{2}}}.
\ee
Hence, the transformation will degenerate at time $T_{\ast}=\tau(T_{\ast})$ (it will be shown that $1-\dot{\tau}(t)>0$) and is a diffeomorphism for $t<\tau(t)$. Then, \eqref{burgersdamp} will be transformed into
\bee\label{ssburgersdamp}
\left(\dfrac{\pa}{\pa s}-\dfrac{1}{2}\right)U+\left(\dfrac{3}{2}y+\dfrac{U}{1-\dot{\tau}}-e^{\frac{s}{2}}\dfrac{\dot{\xi}(t)}{1-\dot{\tau}}\right)\pa_y U=
-\dfrac{ae^{-s}U}{1-\dot{\tau}},
\ee
with the initial condition
\bee\label{burgersinitial}
U(y,0)=f(x), U(0,0)=0, \pa_yU(0,0)=-1.
\ee
To obtain a global solution for \eqref{ssburgersdamp} is similar to the former case\footnote{In this case, define the characteristics as
\bee
\dfrac{dY}{ds}=\frac{3}{2}Y+\frac{1}{1-\dot{\tau}}U(Y,s),\quad Y(0)=y_0,
\ee
along which $U\circ Y=e^{\frac{s}{2}}e^{-a(\tau-e^{-s}+1)}U_0(y_0):=g(s)U_0(y_0)$. Then, $U$ can be solved implicitly as $U(y,s)=g(s)U_0(e^{-\frac{3}{2}s}y-\int_0^s\frac{e^{\frac{3}{2}s'}}{1-\dot{\tau}}g(s')ds'\frac{1}{g(s)} U(y,s))$. To show the global existence of $U$, it suffices to show
\bee
1+\frac{1}{a}(1-e^{-a(t+1)})U_0'\neq0,\quad \forall s\in[0,\infty),
\ee
which is guaranteed by the initial data and the definition of $T_{\ast}$ (see\eqref{xittaut}).}. Here, we focus on the evolutions of $\tau(t)$ and $\xi(t)$. Differentiating \eqref{ssburgersdamp} w.r.t $y$ yields
\bee\label{sspayu}
\left(\dfrac{\pa}{\pa s}-\dfrac{1}{2}\right)\pa_yU+\left(\dfrac{3}{2}+\dfrac{1}{1-\dot{\tau}}\pa_yU\right)\pa_yU+\left(\dfrac{3}{2}y+\dfrac{U}{1-\dot{\tau}}-
e^{\frac{s}{2}}\dfrac{\dot{\xi(t)}}{1-\dot{\tau}}\right)\pa_y^2U=-\dfrac{ae^{-s}\pa_yU}{1-\dot{\tau}}.
\ee
Now we postulate that $U(0,s)=0$ and $\pa_yU(0,s)=-1$ for all $s\in[0,\infty)$. This can be achieved by choosing $\tau(t)$ and $\xi(t)$ suitably since initially, it holds that
$U(0,0)=0$ and $\pa_yU(0,0)=-1$. Evaluating \eqref{ssburgersdamp} and \eqref{sspayu} at $y=0$ yields
\bee
\left\{\bes
&\dot{\xi}(t)=0,\\
&\dot{\tau}(t)=ae^{-s}=a(\tau(t)-t),
\end{split}\right.
\ee
which implies
\bee\label{xittaut}
\left\{\bes
&\xi(t)=\xi(-1)=0,\\
&\tau(t)=t+e^{a(1+t)}+\dfrac{1}{a}(1-e^{a(1+t)}).
\end{split}\right.
\ee
To conclude, it holds that
\begin{itemize}
\item if $a\geq1$, then $\tau(t)\geq t$ for all $t\geq -1$, and one obtains a global classical solution to \eqref{burgersdamp};
\item if $a<1$\footnote{If one assumes generally that $\pa_xf(0)=-c$, then this case becomes "if $a<\frac{1}{c}$, then a shock forms at $T_{\ast}=-\dfrac{1}{a}\ln(1-\frac{a}{c})-1$." Therefore, for fixed $a$, small data(i.e.$c\leq\frac{1}{a}$) leads to the global solution while large data leads to the shock formation in finite time.}, then the transformation between Cartesian coordinates and self-similar coordinates will degenerate at time $T_{\ast}=\tau(T_{\ast})$ which can be computed as
$T_{\ast}=-\dfrac{1}{a}\ln(1-a)-1$, and
      \bee
      \lim_{t\to T_{\ast}}\pa_xu(\xi(t),t)=\lim_{t\to T_{\ast}}e^s\pa_yU(0,s)=\lim_{t\to T_{\ast}}\dfrac{-1}{\tau(t)-t}\to -\infty,
      \ee
      i.e. the solution will form a shock at $(x,t)=(0,T_{\ast})$. In this case, one sees that the blow-up location doesn't shift.
\end{itemize}
\subsection{\textbf{Global existence to the Burgers' equation}}
\hs One can also use the self-similar coordinates to show the global existence to \eqref{burgers} with the initial data being everywhere increasing. For simplicity, we assume that $\pa_xu_0=\pa_xu(x,t=1)$(start from $t=1$) attains its min at $x=0$ with $\pa_xu_0(0)\geq0$. Then, the classical results show that \eqref{burgers} admits a global "rarefaction wave" solution. To this end, define the following self-similar transformation:
\begin{equation}\label{globalss}
\bes
s&=\ln(t),\\
y&=e^{-\frac{3s}{2}}x,\\
u(x,t)&=e^{\frac{s}{2}}U(y,s).
\end{split}
\end{equation}
Then, \eqref{burgers} is transformed into
\bee\label{burgers2'}
\left(\pa_s+\dfrac{1}{2}\right)U+\left(U-\dfrac{3}{2}y\right)\pa_yU=0,
\ee
with the initial condition
\bee
\pa_yU(y,0)=\pa_yU_0(y)=\pa_xu_0\geq0.
\ee
Note that in this case, the Jacobian of the transform is given as
\bee
\dfrac{\pa(y,s)}{\pa(x,t)}=\left|\begin{array}{cc}
e^{-\frac{3s}{2}}&-\dfrac{3}{2}ye^{-s}\\
0&e^{-s}
\end{array}\right|=e^{-\frac{5s}{2}}=\dfrac{1}{t^{\frac{5}{2}}}>0,
\ee
which implies the transform \eqref{globalss} is a global diffeomorphism. Hence, if one obtains a global solution in the self-similar coordinates, then the global existence of original Burgers' equation in Cartesian coordinates is automatically obtained. It can be shown in the similar way that \eqref{burgers2'} admits a global solution on $s\in[0,\infty)$ which converges to $\hat{U}$ pointwisely,
\[
\lim_{s\to\infty}|U(y,s)-\hat{U}(y)|=0,\quad \forall y\in \mbr,
\]
where $\hat{U}=\hat{U}(y)$ is the solution of
\bee
\dfrac{1}{2}\hat{U}+(\hat{U}-\dfrac{3}{2}y)\pa_y\hat{U}=0.
\ee
Indeed, $\hat{U}$ can be solved implicitly as
\bee
y=\hat{U}+\hat{U}^3,\hs y\in(-\infty,+\infty).
\ee
Then, for all $x_0\in \mbr$, it holds that
\bee
|\pa_xu(x_0,t)|=e^{-s}|\pa_yU(y_0,s)|\leq C\dfrac{1}{t},\quad \forall t\geq1.
\ee
\begin{remark}
It seems that the condition $\pa_yU_0(0)=\pa_xu_0(0)\geq0$ which is the global minimum doesn't use for the global existence of \eqref{burgers2'}. However, one can check that if one defines the characteristics of \eqref{burgers2'} as
\[
\dfrac{dY}{ds}=U(Y,s)-\dfrac{3}{2}Y,\hs Y(0)=y_0,
\]
along which $U\circ Y=e^{-\frac{s}{2}}U_0(y_0)$, then $e^{\frac{3s}{2}}Y=y_0+U_0(y_0)(e^{s}-1)$, which allows one to solve $U$ implicitly as $U(y,s)=e^{-\frac{s}{2}}U_0\left(
e^{\frac{3s}{2}}y-e^{\frac{s}{2}}(e^s-1)U\right)$. Then, this formula indeed defines a global solution $U$ provided
\bee
1+(e^s-1)U_0'\neq 0,
\ee
which is guaranteed by the assumption on the initial data.
\end{remark}
\section{\textbf{Coordinates transformations and the self-similar Euler system}}\label{section3}
\hs To study the structure of shock front and introduce the self-similar coordinates, we will construct 10 modulation variables to control the following entities\footnote{The choice of these modulation variables are derived from the invariance of the equations\eqref{euler} under time re-scaling and Galilean transformations which including the space transformation, the time translation, the union motion (the shear transformation) and the space rotation.}:
\bee
\left\{
\bes
&\text{the speed shock formation} \gets\kappa(t)\in\mbr,\\
&\text{the blow up time} \gets\tau(t)\in\mbr,\\
&\text{the blow up location} \gets\xi(t)\in\mbr^3,\\
&\text{the blow up direction} \gets n(t)\in\mbr^3(\text{ only 2 freedoms}),\\
&\text{the tangent surface of shock front} \gets \phi(t)\in\mbr^3.
\end{split}\right.
\ee
Then, the following coordinates transformations will be used:
\begin{align*}
&\text{original physical variables }\xrightarrow[\text{rotation,translation}]{\text{time-rescaling}} \text{the Galilean coordinates}
\xrightarrow{\text{shear-transformation}}\\
 &\text{flatted coordinates w.r.t the shock front}\xrightarrow[\text{transformation}]{\text{self-similar}} \text{self-similar coordinates}
\end{align*}
\subsection{\textbf{The coordinates under Galilean transform}}
   First, do time-rescaling as
   \[
   \mathrm{t}\ra t=\dfrac{1+\a}{2}\mathrm{t},
   \]
   so that all the modulation variables are defined w.r.t the rescaled time (i.e. $t$).\\
   \hs Given an unit vector $n(t)=(\sqrt{1-n_2^2(t)-n_3^2(t)},n_2(t),n_3(t))=(n_1,\breve{n}(t))$, one can generate a rotation matrix which transforms $e_1$ to $n_1(t)$. Precisely, let
   \bee
   R(t)=\left(\begin{array}{ccc}
   n_1&-n_2&-n_3\\
   n_2&1-\frac{n_2^2}{1+n_1}&-\frac{n_2n_3}{1+n_1}\\
   n_3&-\frac{n_2n_3}{1+n_1}&1-\frac{n_3^2}{1+n_1}
   \end{array}\right).
   \ee
   Then, the rotation matrix $R(t)$ rotates $e_1$ to $n(t)$. Define the new coordinates as
   \bee\label{Galilean}
   (\tilde{x},t)=(R^T(\mathit{x}-\xi(t)),t)
   \ee
   and the corresponding fluid variables as
   \bee
   \tilde{\si}(\tilde{x},t)=\si(\mathit{x},t),\quad \tilde{u}(\tilde{x},t)=
   R^T(t)u(\mathit{x},t),\quad \tilde{k}(\tilde{x},t)=k(x,t).
   \ee
   Note that
   \begin{align}
   \dfrac{\pa}{\pa t}&=\dfrac{1+\a}{2}\dfrac{\pa}{\pa\mathit{t}}+\tilde{v}\dfrac{\pa}
   {\pa\tilde{x}_i},\\
   \dfrac{\pa}{\pa \mathit{x}_i}&=R_{ik}\dfrac{\pa}{\pa \tilde{x}_k},
   \end{align}
   where
   \bee\label{definitionQ}
   \tilde{v}=\dot{Q}\tilde{x}-R^T\dot{\xi},\quad  \dot{Q}=\dot{R}^TR,\quad \dot{Q}_{ij}=-\dot{Q}_{ji}.
   \ee
   Hence, the Euler system \eqref{euler2} is transformed in $(\tilde{x},t)$ as
   \begin{equation}\label{euler3}
   \bes
   \dfrac{1+\a}{2}\pa_t\tilde{\si}+(\tilde{u}+\tilde{v})\cdot\nabla_{\tilde{x}}\tilde{\si}+\a\si div_{\tilde{x}}\tilde{u}&=0,\\
   \dfrac{1+\a}{2}\pa_t\tilde{u}-\dot{Q}\tilde{u}+(\tilde{u}+\tilde{v})\cdot\nabla_{\tilde{x}}\tilde{u}+
   \a\tilde{\si}\nabla_{\tilde{x}}\tilde{\si}&=\dfrac{\a}{2\ga}\tilde{\si}^2\nabla_{\tilde{x}}\tilde{k}-a\tilde{u},\\
   \dfrac{1+\a}{2}\pa_t\tilde{k}+(\tilde{u}+\tilde{v})\cdot\nabla_{\tilde{x}}\tilde{k}&=0,
   \end{split}
   \end{equation}
   together with the vorticity equation
   \bee
   \dfrac{1+\a}{2}\pa_t\tilde{\zeta}-\dot{Q}\tilde{\zeta}+(\tilde{u}+\tilde{v})\cdot\nabla_{\tilde{x}}\tilde{\zeta}-
   \tilde{\zeta}\cdot\nabla_{\tilde{x}}\tilde{u}=\dfrac{\a}{\ga}\dfrac{\tilde{\si}}{\tilde{\p}}
   \nabla_{\tilde{x}}\tilde{\si}\times\nabla_{\tilde{x}}\tilde{k}-a\tilde{\zeta}.
   \ee
   \subsection{\textbf{The flatted coordinates with respect to shock front}}
   In order to see the geometry of the shock, define the "tangent" quadratic surface as
   \bee
   (f(\tilde{x}_2,\tilde{x}_3,t),\tilde{x}_2,\tilde{x}_3),
   \ee
   where $f$ is a quadratic form given as
   \bee
   f(\breve{\tilde{x}},t)=\dfrac{1}{2}\phi_{\a\be}\tilde{x}_{\a}\tilde{x}_{\be},
   \ee
   with $\phi$ being a symmetric $2-$tensor. Then, define the shear transformation as
   \bee
   (x_1,x_2,x_3)=(\tilde{x}_1-f(\tilde{x}_1,\tilde{x}_2,t),\tilde{x}_2,\tilde{x}_3),
   \ee
   which flattens the "tangent" surface. Note that together with the coordinates transform, the "basis" $(n(t),\tilde{e}_2=R(t)e_2,\tilde{e}_3=R(t)e_3)$ transforms automatically into
   \begin{align}\label{flatbasis}
   &N(t)=J^{-1}(1,-f_{,2},-f_{,3}),\\& T^{2}=R(N(t))e_2=\left(\dfrac{f_{,2}}{J},1-\dfrac{(f_{,2})^2}{J(J+1)},\dfrac{-f_{,2}f_{,3}}{J(J+1)}\right),\\
   &T^{3}=R(N(t))e_3=\left(\dfrac{f_{,3}}{J},\dfrac{-f_{,2}f_{,3}}{J(J+1)},1-\dfrac{(f_{,3})^2}{J(J+1)}\right),
   \end{align}
   where $J>0$ and $J^2=1+|f_{,2}|^2+|f_{,3}|^2$ and $R(N(t))$ is the rotation matrix generated by $N(t)$ (with the role of $n(t)$ replaced by $N(t)$).
   \begin{remark}
   Note that the $2^{nd}$ fundamental form of the "tangent" surface is given by $\text{II}=J^{-1}\phi_{\a\be}d\tilde{x}^{\a}d\tilde{x}^{\be}$. So, the functions $\fe$ indeed reveal the geometry structure of the shock front.
   \end{remark}
   Denote the corresponding fluid variables as
   \begin{align*}
   \mathring{u}(x,t)&=\tilde{u}(\tilde{x},t),\quad \mathring{\si}(x,t)=\tilde{\si}(\tilde{x},t),\quad \mathring{k}(x,t)=\tilde{k}(\tilde{x},t),\\
   \mathring{v}(x,t)&=\tilde{v}(\tilde{x},t),\quad \mathring{\zeta}(x,t)=\tilde{\zeta}(\tilde{x},t).
   \end{align*}
   Note that
   \begin{align*}
   \dfrac{\pa}{\pa t}&=\dfrac{\pa }{\pa t}-\dot{f}\dfrac{\pa}{\pa x_1},\quad \dfrac{\pa}{\pa \tilde{x}_1}=\dfrac{\pa}{\pa x_1},\quad \dfrac{\pa}{\pa \tilde{x}_{\a}}=\dfrac{\pa}{\pa x_{\a}}-f_{,\a}\dfrac{\pa }{\pa x_1}.
   \end{align*}
   Then, by defining the following constants
   \bee
   \be_1=\dfrac{1}{1+\a},\be_2=\dfrac{1-\a}{1+\a},\be_3=\dfrac{\a}{1+\a},\be_4=\dfrac{\be_3}{1+2\a},
   \ee
    the Euler system \eqref{euler3} is transformed into
   \begin{equation}\label{euler4}
   \bes
   &\pa_t\mathring{\si}+2\be_1\left[JN\cdot(\mathring{u}+\mathring{v})-\dfrac{\dot{f}}{2\be_1}\right]\pa_1\mathring{\si}+2\be_3(\mathring{u}_{\mu}+\mathring{v}_{\mu})\pa_{\mu}\mathring{\si}
   +2\be_3\mathring{\si}[\pa_1\mathring{u}\cdot JN+\pa_{\mu}\mathring{u}_{\mu}]=0,\\
   &\pa_t\mathring{u}-2\be_1\mathring{u}\dot{Q}+2\be_1\left[JN\cdot(\mathring{u}+\mathring{v})-\dfrac{\dot{f}}{2\be_1}\right]\pa_1\mathring{u}
   +2\be_3(\mathring{u}_{\mu}+\mathring{v}_{\mu})\pa_{\mu}\mathring{u}
   +2\be_3\mathring{\si}[JN\pa_1+\da^{\cdot\mu}\pa_{\mu}]\mathring{\si}\\
   &=\be_4\mathring{\si}^2[JN\pa_1+\da^{\cdot\mu}\pa_{\mu}]\mathring{k}-2\be_1a\mathring{u},\\
   &\pa_t\mathring{k}+2\be_1\left[JN\cdot(\mathring{u}+\mathring{v})-\dfrac{\dot{f}}{2\be_1}\right]\pa_1\mathring{k}+2\be_3(\mathring{u}_{\mu}+\mathring{v}_{\mu})\pa_{\mu}\mathring{k}=0,
   \end{split}
   \end{equation}
   together with the vorticity equation
   \bee
   \bes
   &\pa_t\mathring{\zeta}-2\be_1\mathring{\zeta}\dot{Q}+2\be_1\left[JN\cdot(\mathring{u}+\mathring{v})-\dfrac{\dot{f}}{2\be_1}\right]\pa_1\mathring{\zeta}
   +2\be_3(\mathring{u}_{\mu}+\mathring{v}_{\mu})\pa_{\mu}\mathring{\zeta}
   -2\be_1[JN\cdot\mathring{\zeta}\pa_1\mathring{u}+\mathring{\zeta}_{\mu}\pa_{\mu}\mathring{u}]\\
   &=2\be_4\dfrac{\mathring{\si}}{\mathring{\p}}\nabla_{\tilde{x}}\mathring{\si}\times\nabla_{\tilde{x}}\mathring{k}-2\be_1a\mathring{\zeta},
   \end{split}
   \ee
   where in terms of $(N,T^{\a})$ basis,
   \bee
   \bes
   \nabla_{\tilde{x}}\mathring{\si}\times\nabla_{\tilde{x}}\mathring{k}&=(\pa_{T^2}\mathring{\si}\pa_{T^3}\mathring{k}-\pa_{T^3}\mathring{\si}\pa_{T^2}\mathring{k})N+
   (\pa_{T^3}\mathring{\si}\pa_{N}\mathring{k}-\pa_{N}\mathring{\si}\pa_{T^3}\mathring{k})T^2\\
   &+(\pa_{N}\mathring{\si}\pa_{T^2}\mathring{k}-\pa_{T^2}\mathring{\si}\pa_{N}\mathring{k})T^3
   \end{split}
   \ee
   \begin{remark}
   The term $\pa_{N}\mathring{\si}$ in $T^{\a}$ direction is the potentially dangerous term since this term blows up at the shock time which will be proved later.
   \end{remark}
   In order to see the intersection of different wave families, it's convenient to introduce the Riemann variables
   \bee
   w=\mr{u}\cdot N+\mr{\si},\hs z=\mr{u}\cdot N-\mr{\si},\hs a_{\mu}=\mr{u}\cdot T^{\mu}.
   \ee
   \begin{remark}
   Later, one can see that actually only the quantity $\pa R_+=\pa_{N}w$ blows up at the shock time and location, the remaining quantities are bound up to the shock, which is coincide to the 1D case.
   \end{remark}
   In terms of Riemann variables, the Euler system \eqref{euler4} can be rewritten as the following $(w,z,a_{\mu},k)$ system:
   \bee\label{wzak}
   \bes
   &\pa_tw+v_{w}\cdot\nabla_xw=F_w,\\
   &\pa_tz+v_{z}\cdot\nabla_xz=F_z,\\
   &\pa_ta_{\mu}+v_{a}\cdot\nabla_xa_{\mu}=F_{a_{\mu}},\\
   &\pa_tk+v_k\cdot\nabla_xk=0,
   \end{split}
   \ee
   where
   \begin{align*}
   v_w&=(Jw+J\be_2z+2\be_1J\mr{v}\cdot N-\dot{f},wN_{\mu}-\be_2zN_{\mu}+2\be_1\mr{v}_{\mu}+2\be_1a_{\nu}T^{\nu}_{\mu}),\\
   v_z&=(Jz+J\be_2w+2\be_1J\mr{v}\cdot N-\dot{f},zN_{\mu}+\be_2wN_{\mu}+2\be_1\mr{v}_{\mu}+2\be_1a_{\nu}T^{\nu}_{\mu}),\\
   v_{a}&=(J\be_1z+J\be_1w+2\be_1J\mr{v}\cdot N-\dot{f},\be_1zN_{\mu}+\be_1wN_{\mu}+2\be_1\mr{v}_{\mu}+2\be_1a_{\nu}T^{\nu}_{\mu}),\\
   v_k&=(2\be_1J(\mr{u}+\mr{v})\cdot N-\dot{f}, 2\be_1(\mr{u}_{\mu}+\mr{v}_{\mu})),
   \end{align*}
   and
   \begin{align*}
   F_w&=\dot{N}\cdot a_{\nu}T^{\nu}+2\be_1\dot{Q}_{ij}a_{\nu}T^{\nu}_jN_i+2\be_1(\mr{v}_{\mu}+\mr{u}\cdot NN_{\mu}+a_{\nu}T^{\nu}_{\mu})N_{i,\mu}a_{\ga}T^{\ga}_i\\
   &-2\be_3\mr{\si}[\mr{u}\cdot N\pa_{\mu}N_{\mu}+\pa_{\mu}a_{\nu}\cdot T^{\nu}_{\mu}+a_{\nu}\pa_{\mu}T^{\nu}_{\mu}]\\
   &+\be_4\mr{\si}^2(J\pa_1\mr{k}+N_{\mu}\pa_{\mu}\mr{k})-2\be_1a\mr{u}\cdot N,\\
   F_z&=\dot{N}\cdot a_{\nu}T^{\nu}+2\be_1\dot{Q}_{ij}a_{\nu}T^{\nu}_jN_i+2\be_1(\mr{v}_{\mu}+\mr{u}\cdot NN_{\mu}+a_{\nu}T^{\nu}_{\mu})N_{i,\mu}a_{\ga}T^{\ga}_i\\
   &+2\be_3\mr{\si}[\mr{u}\cdot N\pa_{\mu}N_{\mu}+\pa_{\mu}a_{\nu}\cdot T^{\nu}_{\mu}+a_{\nu}\pa_{\mu}T^{\nu}_{\mu}]\\
   &+\be_4\mr{\si}^2(J\pa_1\mr{k}+N_{\mu}\pa_{\mu}\mr{k})-2\be_1a\mr{u}\cdot N,\\
   F_{a_{\nu}}&=-2\be_3\mr{\si}T^{\nu}_{\mu}\pa_{\mu}\mr{\si}+2\be_1\dot{Q}_{ij}[a_{\ga}T^{\ga}_j+\mr{u}\cdot NN_j]T^{\nu}_i+2\be_1\dot{T}^{\nu}_i[\mr{u}\cdot NN_i+a_{\ga}T^{\ga}_i]\\
   &+2\be_1(\mr{v}_{\mu}+\mr{u}\cdot NN_{\mu}+a_{\ga}T^{\ga}_{\mu})T^{\nu}_{i,\mu}[\mr{u}\cdot NN_i+a_{\ga}T^{\ga}_i]\\
   &+\be_4\mr{\si}^2T^{\nu}_{\mu}\pa_{\mu}\mr{k}-2\be_1aa_{\nu}.\\
   \end{align*}
   \subsection{\textbf{The self-similar coordinates}}
  \hs Define the following self-similar transformation
   \bee\label{selfsimilartranseuler}
   \bes
   s&=s(t)=-\ln(\tau(t)-t),\\
   y_1&=y_1(x_1,t)=\dfrac{x_1}{(\tau(t)-t)^{\frac{3}{2}}}=x_1e^{\frac{3}{2}s},\\
   y_{\mu}&=y_{\mu}(x_{\mu},t)=\dfrac{x_{\mu}}{(\tau(t)-t)^{\frac{1}{2}}}=x_1e^{\frac{1}{2}s},
   \end{split}
   \ee
   and the corresponding fluid variables
   \bee\label{selfsimilarfluid}
   \bes
   w(x,t)&=e^{-\frac{s}{2}}W(y,s)+\kappa(t),\\
   z(x,t)&=Z(y,s), \hs a_{\mu}(x,t)=A_{\mu}(y,s),\\
   k(x,t)&=K(y,s),\hs\mr{\zeta}(x,t)=\Omega(y,s),\\
   \mr{\si}(x,t)&=S(y,s),\\
   v(x,t)&=V(y,s)=\dot{Q}_{i1}\left(e^{-\frac{3}{2}s}y_1+\dfrac{1}{2}e^{-s}\fe_{\mu\nu}y_{\mu}y_{\nu}\right)+e^{-\frac{s}{2}}\dot{Q}_{i\mu}y_{\mu}-R_{ji}\cdot{\dot{\xi}}_j.
   \end{split}
   \ee
   \begin{remark}
   The reason of choosing the indexes in \eqref{selfsimilartranseuler} is similar to the Burger's case.
   \end{remark}
   Note that
   \bee
   \bes
   \dfrac{ds}{dt}&=(1-\dot{\tau(t)})e^s=e^s\cdot\dfrac{1}{\be_{\tau}},\\
   \dfrac{\pa y_1}{\pa t}&=\dfrac{3}{2}x_1e^{\frac{3}{2}s}\dfrac{ds}{dt}=\dfrac{3}{2}y_1\dfrac{e^s}{\be_{\tau}},\\
   \dfrac{\pa y_{\mu}}{\pa t}&=\dfrac{1}{2}y_{\mu}\dfrac{e^s}{\be_{\tau}},\\
   \dfrac{\pa y_1}{\pa x_1}&=e^{\frac{3}{2}s},\hs \dfrac{\pa y_{\mu}}{\pa x_{\mu}}=e^{\frac{s}{2}}.
   \end{split}
   \ee
   Then, the Jacobian of the transformation\eqref{selfsimilartranseuler} is given as
   \bee\label{jocabian}
   \dfrac{\pa(y,s)}{\pa(x,t)}=\dfrac{e^s}{\be_{\tau}}\left|\begin{array}{cccc}
   1&0&0&0\\
   \dfrac{3}{2}y_1&e^{\frac{3}{2}s}&0&0\\
   \dfrac{1}{2}y_2&0&e^{\frac{s}{2}}&0\\
   \dfrac{1}{2}y_3&0&0&e^{\frac{s}{2}}
   \end{array}\right|=\dfrac{e^{\frac{7}{2}s}}{\be_{\tau}}=\dfrac{1}{(\tau(t)-t)^{\frac{7}{2}}}\dfrac{1}{\be_{\tau}},
   \ee
   which means the transformation is a diffeomorphism as long as $\tau(t)>t$ (it will be shown that $\be_{\tau}>0$ and becomes degenerated as $\tau(t)\to t$ (this is how we define the shock time). Therefore, shock formation to \eqref{euler} is equivalent to that
   \begin{itemize}
   \item the gradient of Riemann variable $\mathcal{R}_+=w=\mr{u}\cdot N+\mr{\si}$ along $N$-direction become infinite;
   \item the transformation between the physical variables (the Cartesian coordinates) and the self-similar coordinates becomes degenerated;
   \item the normal $N(t)$ of the "tangent" surface $(f(\check{\tilde{x}},t),\tilde{x}_2,\tilde{x}_3)$ of the shock front becomes "horizontal".
   \end{itemize}
   \begin{remark}
   The normal vector $N(t)$ becomes "horizontal" above is w.r.t the $O(\ep)$ scale. Precisely, initially
   \bee
   \dfrac{|N(t)-e_1|}{\ep}\sim1,
   \ee
   and as time approaches to the shock time
   \bee
   \dfrac{|N(t)-e_1|}{\ep}\to 0.
   \ee
   This is because we set initial time at $t=-\ep$ instead of $-1$.
   \end{remark}
   In terms of the self-similar coordinates $(y,s)$, the $(w,z,a_{\mu},\mr{k})$ system \eqref{wzak} is transformed into the following $(W,Z,A_{\mu},K)$ system:
   \bee\label{WZAK}
   \bes
   \left(\pa_s-\dfrac{1}{2}\right)W+(V_{W}\cdot\nabla)W&=\mathcal{F}_{W},\\
   \pa_sZ+(V_Z\cdot\nabla)Z&=\mathcal{F}_Z,\\
   \pa_sA_{\nu}+(V_{U}\cdot\nabla)A_{\nu}&=\mathcal{F}_{A_{\nu}},\\
   \pa_sK+(V_U\cdot\nabla)K&=0,
   \end{split}
   \ee
   where
   \bee\label{defGWGZGU}
   \bes
   V_W&=(g_W+\dfrac{3}{2}y_1,h_W^{\mu}+\dfrac{1}{2}y_{\mu}),\\
   g_W&=\be_{\tau}JW+\be_{\tau}e^{\frac{s}{2}}\left(-\dot{f}+2\be_1JV\cdot N+J\be_2Z+J\kappa\right):=\be_{\tau}JW+G_W,\\
   h_W^{\mu}&=\be_{\tau}e^{-s}N_{\mu}W+\be_{\tau}e^{-\frac{s}{2}}\left(2\be_1V_{\mu}+2\be_1A_{\nu}T^{\nu}_{\mu}-\be_2ZN_{\mu}+N_{\mu}\kappa\right),\\
   V_Z&=(g_Z+\dfrac{3}{2}y_1,h_Z^{\mu}+\dfrac{1}{2}y_{\mu}),\\
   g_Z&=\be_2\be_{\tau}JW+\be_{\tau}e^{\frac{s}{2}}\left(-\dot{f}+2\be_1JV\cdot N+JZ+J\be_2\kappa\right):=\be_2\be_{\tau}JW+G_Z,\\
   h_Z^{\mu}&=\be_2\be_{\tau}e^{-s}N_{\mu}W+\be_{\tau}e^{-\frac{s}{2}}\left(2\be_1V_{\mu}+2\be_1A_{\nu}T^{\nu}_{\mu}+\be_2\kappa N_{\mu}+N_{\mu}Z\right),\\
   V_U&=(g_U+\dfrac{3}{2}y_1,h_U^{\mu}+\dfrac{1}{2}y_{\mu}),\\
   g_U&=\be_1\be_{\tau}JW+\be_{\tau}e^{\frac{s}{2}}\left(-\dot{f}+2\be_1JV\cdot N+J\be_1Z+J\be_1\kappa\right):=\be_1\be_{\tau}JW+G_U,\\
   h_U^{\mu}&=\be_1\be_{\tau}e^{-s}N_{\mu}W+\be_{\tau}e^{-\frac{s}{2}}\left(2\be_1V_{\mu}+2\be_1A_{\nu}T^{\nu}_{\mu}+\be_1\kappa N_{\mu}+\be_1N_{\mu}Z\right),
   \end{split}
   \ee
   and
   \begin{align}
   \mathcal{F}_W&=F_W-e^{-\frac{s}{2}}\be_{\tau}\dot{\kappa}-2\be_1\be_{\tau}e^{-\frac{s}{2}}aU\cdot N,\\
   \begin{split}\label{FW}
   F_W&=-2\be_3\be_{\tau}S\pa_{\mu}A_{\nu}T_{\mu}^{\nu}+\be_{\tau}e^{-\frac{s}{2}}A_{\nu}T^{\nu}\cdot \dot{N}+2\be_1\be_{\tau}e^{-\frac{s}{2}}\dot{Q}_{ij}
   A_{\nu}T^{\nu}_jN_i\\
   &+2\be_1\be_{\tau}(V_{\mu}+U\cdot NN_{\mu}+A_{\ga}T^{\ga}_{\mu})N_{i,\mu}A_{\ga}T^{\ga}_i-2\be_3\be_{\tau}e^{-\frac{s}{2}}S(\pa_{\mu}N_{\mu}U\cdot N+
   A_{\nu}\pa_{\mu}T^{\nu}_{\mu})\\
   &+\be_4\be_{\tau}S^2(Je^s\pa_1K+N_{\mu}\pa_{\mu}K),
   \end{split}\\
   \mathcal{F}_Z&=F_Z-e^{-s}\be_{\tau}\dot{\kappa}-2\be_1\be_{\tau}e^{-s}aU\cdot N,\\
   \begin{split}\label{FZ}
   F_Z&=2\be_3\be_{\tau}S\pa_{\mu}A_{\nu}T_{\mu}^{\nu}+\be_{\tau}e^{-s}A_{\nu}T^{\nu}\cdot \dot{N}+2\be_1\be_{\tau}e^{-s}\dot{Q}_{ij}
   A_{\nu}T^{\nu}_jN_i\\
   &+2\be_1\be_{\tau}e^{-s}(V_{\mu}+U\cdot NN_{\mu}+A_{\ga}T^{\ga}_{\mu})N_{i,\mu}A_{\ga}T^{\ga}_i+2\be_3\be_{\tau}e^{-\frac{s}{2}}S(\pa_{\mu}N_{\mu}U\cdot N+
   A_{\nu}\pa_{\mu}T^{\nu}_{\mu})\\
   &+\be_4\be_{\tau}S^2(Je^{\frac{s}{2}}\pa_1K+N_{\mu}e^{-\frac{s}{2}}\pa_{\mu}K),
   \end{split}\\
   \mathcal{F}_{A_{\nu}}&=F_{A_{\nu}}-2\be_1\be_{\tau}e^{-s}aA_{\nu},\\
   \begin{split}\label{FA}
   F_{A_{\nu}}&=-2\be_3\be_{\tau}S\pa_{\mu}ST_{\mu}^{\nu}+\be_{\tau}e^{-s}(U\cdot NN+A_{\ga}T^{\ga})\cdot\dot{T}^{\nu}+2\be_1\be_{\tau}e^{-s}\dot{Q}_{ij}
   [U\cdot N N_j+A_{\ga}T^{\ga}_j]T^{\nu}_i\\
   &+2\be_1\be_{\tau}e^{-s}(V_{\mu}+U\cdot NN_{\mu}+A_{\ga}T^{\ga}_{\mu})\pa_{\mu}T^{\nu}_{i}[U\cdot NN_i+A_{\ga}T^{\ga}_i]+\be_4\be_{\tau}e^{-\frac{s}{2}}S^2T^{\nu}_{\mu}\pa_{\mu}K.
   \end{split}
   \end{align}
   \hs Let $\ga$ be a multi-index $\ga=(\ga_1,\ga_2,\ga_3)=(\ga_1,\breve{\ga})$. Then, acting $\pa^{\ga}$ to the $(W,Z,A_{\nu},K)$ system yields
   \begin{align}
   \left(\pa_s+\dfrac{3\ga_1+\ga_2+\ga_3-1}{2}+\be_{\tau}J(1+\ga_1 1_{|\ga|\geq 2})\pa_1W\right)\pa^{\ga}W+V_W\cdot\nabla\pa^{\ga}W&=F_{W}^{(\ga)},\label{pagaw}\\
   \left(\pa_s+\dfrac{3\ga_1+\ga_2+\ga_3}{2}+\be_2\be_{\tau}J\ga_1 \pa_1W)\right)\pa^{\ga}Z+V_Z\cdot\nabla\pa^{\ga}Z&=F_{Z}^{(\ga)},\\
   \left(\pa_s+\dfrac{3\ga_1+\ga_2+\ga_3}{2}+\be_2\be_{\tau}J\ga_1 \pa_1W)\right)\pa^{\ga}A_{\nu}+V_U\cdot\nabla\pa^{\ga}A_{\nu}&=F_{A_{\nu}}^{(\ga)},\\
   \left(\pa_s+\dfrac{3\ga_1+\ga_2+\ga_3}{2}+\be_1\be_{\tau}J(1+\ga_1 \pa_1W)\right)\pa^{\ga}K+V_U\cdot\nabla\pa^{\ga}K&=F_{K}^{(\ga)},
   \end{align}
   where
   \bee
   \bes
   F_W^{(\ga)}&=\pa^{\ga}F_W-\sum_{\substack{0\leq|\be|\leq|\ga|-1,\\ \be\leq \ga}}\left(C_{\ga}^{\be}\pa^{\ga-\be}G_W\pa^{\be}\pa_1W+
   C_{\ga}^{\be}\pa^{\ga-\be}h_W^{\mu}\pa^{\be}\pa_{\mu}W\right)\\
   &-1_{|\ga|\geq3}\be_{\tau}\sum_{\substack{1\leq|\be|\leq|\ga|-2,\\ \be\leq\ga}}C_{\ga}^{\be}\pa^{\ga-\be}(JW)\pa^{\be}\pa_1W-1_{|\ga|\geq2}\be_{\tau}\sum_{\substack{
   |\be|=|\ga|-1,\\ \be_1=\ga_1}}C_{\ga}^{\be}\pa^{\ga-\be}(JW)\pa^{\be}\pa_1W\\
   &-2\be_1\be_{\tau}e^{-\frac{s}{2}}a\pa^{\ga}(U\cdot N)-
   \be_{\tau}\sum_{\substack{0\leq|\be|\leq|\ga|-1,\\ \be\leq\ga}}C_{\ga}^{\be}\pa^{\ga-\be}J\pa^{\be}W\pa_1W,\\
   F_Z^{(\ga)}&=\pa^{\ga}F_Z-\sum_{\substack{0\leq|\be|\leq|\ga|-1,\\ \be\leq \ga}}\left(C_{\ga}^{\be}\pa^{\ga-\be}G_Z\pa^{\be}\pa_1Z+
   C_{\ga}^{\be}\pa^{\ga-\be}h_Z^{\mu}\pa^{\be}\pa_{\mu}Z\right)\\
   &-1_{|\ga|\geq2}\be_2\be_{\tau}\sum_{\substack{0\leq|\be|\leq|\ga|-2,\\ \be\leq\ga}}C_{\ga}^{\be}\pa^{\ga-\be}(JW)\pa^{\be}\pa_1Z-1_{|\ga|\geq2}\be_2\be_{\tau}\sum_{\substack{
   |\be|=|\ga|-1,\\ \be_1<\ga_1}}C_{\ga}^{\be}\check{\nabla}(JW)\pa^{\be}\pa_1Z\\
   &-2\be_1\be_{\tau}e^{-s}a\pa^{\ga}(U\cdot N),\\
   F_{A_{\nu}}^{(\ga)}&=\pa^{\ga}F_{A_{\nu}}-\sum_{\substack{0\leq|\be|\leq|\ga|-1,\\ \be\leq \ga}}\left(C_{\ga}^{\be}\pa^{\ga-\be}G_U\pa^{\be}\pa_1A_{\nu}+
   C_{\ga}^{\be}\pa^{\ga-\be}h_U^{\mu}\pa^{\be}\pa_{\mu}A_{\nu}\right)\\
   &-1_{|\ga|\geq2}\be_1\be_{\tau}\sum_{\substack{0\leq|\be|\leq|\ga|-2,\\ \be\leq\ga}}C_{\ga}^{\be}\pa^{\ga-\be}(JW)\pa^{\be}\pa_1A_{\nu}-1_{|\ga|\geq2}\be_1\be_{\tau}\sum_{\substack{
   |\be|=|\ga|-1,\\ \be_1<\ga_1}}C_{\ga}^{\be}\pa^{\ga-\be}(JW)\pa^{\be}\pa_1A_{\nu}\\
   &-2\be_1\be_{\tau}e^{-s}a\pa^{\ga}A_{\nu},\\
   F_{K}^{(\ga)}&=-\sum_{\substack{0\leq|\be|\leq|\ga|-1,\\ \be\leq \ga}}\left(C_{\ga}^{\be}\pa^{\ga-\be}G_U\pa^{\be}\pa_1K+
   C_{\ga}^{\be}\pa^{\ga-\be}h_U^{\mu}\pa^{\be}\pa_{\mu}K\right)\\
   &-1_{|\ga|\geq2}\be_1\be_{\tau}\sum_{\substack{0\leq|\be|\leq|\ga|-2,\\ \be\leq\ga}}C_{\ga}^{\be}\pa^{\ga-\be}(JW)\pa^{\be}\pa_1K-1_{|\ga|\geq2}\be_1\be_{\tau}\sum_{\substack{
   |\be|=|\ga|-1,\\ \be_1<\ga_1}}C_{\ga}^{\be}\pa^{\ga-\be}(JW)\pa^{\be}\pa_1K.\\
   \end{split}
   \ee
   \subsection{\textbf{The $(U,P,H)$ system for energy estimates}}
   In order to derive the energy estimates for Euler system, we compute the equations for the velocity, the pressure and the entropy, instead. Let
   \bee
   \mathcal{P}(y,s)=\dfrac{1}{\a}(\ga P(y,s))^{\frac{\a}{\ga}}=\dfrac{1}{\a}(\ga p(x,t))^{\frac{\a}{\ga}},\hs H(y,s)=e^{K(y,s)}{2\ga}.
   \ee
   Then, $(U,\mathcal{P},H)$ satisfies the following system
   \bee\label{UPHsystem}
   \bes
   \pa_sU-2\be_1\be_{\tau}e^{-s}U\dot{Q}+V_U\cdot\nabla_yU+2\be_3\be_{\tau}H^2P\left(JNe^{\frac{s}{2}}\pa_1+e^{-\frac{s}{2}}\da^{\cdot\nu}\pa_{\nu}\right)P+
   2\be_1\be_{\tau}e^{-s}aU&=0,\\
   \pa_sP+V_U\cdot\nabla_yP+2\be_3\be_{\tau}P\left(JN\cdot e^{\frac{s}{2}}\pa_1U+e^{-\frac{s}{2}}\pa_{\mu}U_{\mu}\right)&=0,\\
   \pa_sH+V_U\cdot\nabla_yH&=0.
   \end{split}
   \ee
   \subsection{\textbf{The solution of the self-similar Burgers' equation}}
   Note that if the solution of the equation of $W$ is independent of $s$, then, this equation will "converge" to the following 3D self-similar Burgers' equation (later it will be shown that as $s\to\infty$, $G_W,h_W^{\mu}\to 0$)
    \bee\label{burgers3}
    -\dfrac{1}{2}\bar{W}+(\dfrac{3}{2}y_1+\bar{W})\pa_1\bar{W}+\dfrac{1}{2}y_{\mu}\pa_{\mu}\bar{W}=0.
    \ee
   Indeed, the solution of \eqref{burgers3} can be generated by the solution of 1D self-similar Burgers' equation. Let $\langle\breve{y}\rangle=1+|\breve{y}|^2=1+y_2^2+y_3^2$ and
   \bee\label{definitionbarWy}
   \bar{W}(y_1,\breve{y})=\langle\breve{y}\rangle^{\frac{1}{2}}\bar{U}\left(\langle\breve{y}\rangle^{-\frac{3}{2}}y_1\right),
   \ee
   where $\bar{U}(y)$ is given implicitly by $y=-\bar{U}-\bar{U}^3$ (see\eqref{burgers2}). Then, $\bar{W}(y_1,y_2,y_3)$ solves \eqref{burgers3} and direct computations lead to
   \bee
   |\bar{U}|\les(1+y^2)^{\frac{1}{6}},\quad |\bar{U}'|\les(1+y^2)^{-\frac{1}{3}}.
   \ee
   Moreover, if one defines $\eta(y)=1+y_1^2+|\breve{y}|^6$, then it holds that
   \bee\label{estpabarW}
   |\bar{W}|\les\eta^{\frac{1}{6}}(y),\hs |\pa_1\bar{W}|\les\eta^{-\frac{1}{3}}(y),\hs |\pa_{\mu}\bar{W}|\les 1,\hs |\pa_1\pa_{i}\bar{W}|\les\eta^{-\frac{1}{2}}(y),\hs
   |\breve{\nabla}^2\bar{W}|\les\eta^{-\frac{1}{6}}(y).
   \ee
   Later, it will be proven that $W$ will converge to $\bar{W}$ pointwisely as $s\to\infty$. Then, it's natural to consider the evolution equation of $\tilde{W}=W-\bar{W}$:
   \bee\label{asymeq}
   \left(\pa_s+\be_{\tau}J\pa_1W-\frac{1}{2}\right)\tilde{W}+V_W\cdot\nabla\tilde{W}=\tilde{F}_W,
   \ee
   where
   \bee
   \tilde{F}_W=\mathcal{F}_W-\left[(\be_{\tau}J-1)\bar{W}-G_W\right]\pa_1\bar{W}-h_W^{\mu}\pa_{\mu}\bar{W}.
   \ee
   Applying $\pa^{\ga}$ with $|\ga|\geq1$ to \eqref{asymeq} yields
   \bee\label{pagatildew}
   \left[\pa_s+\dfrac{3\ga_1+\ga_2+\ga_3-1}{2}+\be_{\tau}J(\pa_1\bar{W}+\ga_1\pa_1W)\right]\pa^{\ga}\tilde{W}+V_W\cdot\nabla\pa^{\ga}\tilde{W}=\tilde{F}^{(\ga)}_{W},
   \ee
   where
   \bee\label{tildeFWga}
   \bes
   \tilde{F}_W^{(\ga)}&=\pa^{\ga}\tilde{F}_W-\sum_{\substack{0\leq|\be|\leq|\ga|-1,\\ \be\leq \ga}}\left(C_{\ga}^{\be}\pa^{\ga-\be}G_W\pa^{\be}\pa_1\tilde{W}+
   C_{\ga}^{\be}\pa^{\ga-\be}h_W^{\mu}\pa^{\be}\pa_{\mu}\tilde{W}+\be_{\tau}C_{\ga}^{\be}\pa^{\ga-\be}(J\pa_1\bar{W})\pa^{\be}\tilde{W}\right)\\
   &-1_{|\ga|\geq2}\be_{\tau}\sum_{\substack{0\leq|\be|\leq|\ga|-2,\\ \be\leq\ga}}C_{\ga}^{\be}\pa^{\ga-\be}(JW)\pa^{\be}\pa_1\tilde{W}-1_{|\ga|\geq2}\be_{\tau}\sum_{\substack{
   |\be|=|\ga|-1,\\ \be_1=\ga_1}}C_{\ga}^{\be}\pa^{\ga-\be}(JW)\pa^{\be}\pa_1\tilde{W}.
   \end{split}
   \ee
   \section{\textbf{Initial data construction, Bootstrap assumptions and the main results}}\label{section4}
   \hs Set the initial time to be $t_0=-\ep$. The initial values of the modulation variables are set as follows:
   \begin{align}
   &\kappa(-\ep)=\kappa_0,\hs \tau(-\ep)=\tau_0=0,\hs \xi(-\ep)=\xi_0=0,\\
   &\check{n}(-\ep)=\check{n}_0=0,\hs \fe(-\ep)=\fe_0,
   \end{align}
   with $\kappa_0>1$ so that initially, the Galilean transformation \eqref{Galilean} is the identical map and the shear transformation is given as
   \bee
   x_1=x_1-f_0(\check{\tilde{x}}),\hs x_{\mu}=\tilde{x}_{\mu},
   \ee
   where
   \bee
   f_0(\check{\tilde{x}})=\fe_{0\mu\nu}\tilde{x}_{\a}\tilde{x}_{\be}.
   \ee
   Then, the initial basis $(N_0(\tilde{x}),T_0^{\nu}(\tilde{x}))$ in the flatted coordinates are set as in\eqref{flatbasis}. Finally, set the initial fluid variables as follows:
   \begin{align}
   &u(\tilde{x},\ep)=u_0(\tilde{x}),\hs \p(\tilde{x},\ep)=\p_0(\tilde{x}),\hs k(\tilde{x},-\ep)=\tilde{k}_0(\tilde{x}),\hs\sigma(x,-\ep)=\sigma_0(\tilde{x}),\\
   &\tilde{w}(\tilde{x},-\ep)=\tilde{w}_0(\tilde{x})=u_0(\tilde{x})\cdot N_0(\tilde{x})+\sigma_0(\tilde{x}),\\
   &\tilde{z}(\tilde{x},-\ep)=\tilde{z}_0(\tilde{x})=u_0(\tilde{x})\cdot N_0(\tilde{x})-\sigma_0(\tilde{x}),\\
   &\tilde{a}_{\nu}(\tilde{x},-\ep)=\tilde{a}_{0\nu}(\tilde{x})=u_0(\tilde{x})\cdot T_0^{\nu}(\tilde{x}).
   \end{align}
   We then construct the initial datum for $(\tilde{w}_0(\tilde{x}),\tilde{z}_0(\tilde{x}),\tilde{a}_{0\mu}(\tilde{x}),\tilde{k}_0(\tilde{x}))$ as follows.
   \begin{lem}\label{initialdata}
   Given any functions $\fai_i\in C^{\infty}_c(\mbr^3)$ (for $i=0,1,2,3$) such that $supp\fai_i
   \subset\{0\leq|\tilde{x}|\leq1\}$ for $i=0,1,2,3$. For simplicity, we also assume that $|\pa^{\ga}\fai_i|\leq 1$ for $|\ga|\leq 4$ for $\ga=(\ga_1,\check{\ga})$ being multi-index and $i=0,1,2,3$. Let $\tilde{\eta}(x)=\ep^3+\tilde{x}_1^2+|\check{\tilde{x}}|^6$. Then, define the initial data as follows:
   \begin{align}
   \til{z}_0(\tilde{x})&=\ep\fai_1\left(\ep^{3\a}\tilde{x}_1,\ep^{\a}\check{\tilde{x}}\right),\quad \text{with}\ -\frac{1}{3}\leq\a\leq-\frac{1}{6},\\
   \til{a}_0(\tilde{x})&=\ep\fai_2\left(\ep^{3\a}\tilde{x}_1,\ep^{\a}\check{\tilde{x}}\right),\quad \text{with}\ -\frac{1}{3}\leq\a\leq-\frac{1}{6},\\
   \til{k}_0(\tilde{x})&=\ep\fai_3\left(\ep^{-\frac{1}{2}}\tilde{x}_1,\ep^{-\frac{1}{6}}\check{\tilde{x}}\right),\\
   \til{w}_0(\tilde{x})&=\ep^{\frac{1}{2}}\tilde{\eta}^{\frac{1}{6}}(\tilde{x})\fai_0\left(\ep^{3\be}\tilde{x}_1,\ep^{\be}\check{\tilde{x}}\right)+\bar{W}_{\ep}(\tilde{x}), \quad
   \text{with}\ -\frac{2}{3}\leq\be\leq-\frac{1}{6},
   \end{align}
   where
   \[
   \bar{W}_{\ep}(\tilde{x})=\ep^{\frac{1}{2}}\bar{W}(\ep^{-\frac{3}{2}}\tilde{x}_1,\ep^{-\frac{1}{2}}\check{\tilde{x}}).
   \]
   \end{lem}
   \begin{remark}
   One can check that the initial data in Lemma\ref{initialdata} satisfies all the assumptions in \cite{BSV3Dfulleuler}.
   \end{remark}
   \begin{remark}
   We constrain the support of $(\til{w}_0,\til{z}_0,\til{a}_0,\til{k}_0)$ in $\{|\tilde{x}_1|\leq\ep^{\frac{1}{2}},\ |\check{\tilde{x}}|\leq\ep^{\frac{1}{6}}\}$. Indeed, the initial data in this Lemma can be taken generally such as $\til{z}_0(\tilde{x})=\ep\fai_1\left(\ep^{\a}\tilde{x}_1,\ep^{\be}\check{\tilde{x}}\right)$ with $ -1\leq\a\leq-\frac{1}{2}$ and $-\frac{1}{2}\leq\be\leq-\frac{1}{6}$ or more generally $\til{z}_0(\tilde{x})=\ep^{\a}\fai_1\left(\ep^{\be}\tilde{x}_1,\ep^{\ga}\check{\tilde{x}}\right)$ with $\be\geq\max\{-a,-\frac{3}{4}-\frac{\a}{2}\}$, $\ga\geq\frac{1}{2}-\a$, $\be+\ga\geq-\frac{1}{2}-\a$ and $\a\geq1$, $\be\leq-\frac{1}{2}$, $\ga\leq-\frac{1}{6}$.
   \end{remark}
   \begin{remark}
   Similar to the case for the Burgers equation \eqref{burgersinitial} in section2 to control the modulation variables, we further assume that
   \begin{align}\label{eulerinitial}
   \pa_{x_1}\tilde{w}(0)&=-\frac{1}{\ep}\footnotemark,\\
   \pa^{\ga}\tilde{w}(0)&=0,
   \end{align}
   \footnotetext{Similar to the case in section2, this condition leads to the shock formation in finite time. Later one can see that if $\ep$ is sufficiently small, then a shock forms at $T_{\ast}=O(\ep)$ while for large $\ep$, the solution exists globally. One can also assume generally that $\pa_{x_1}\tilde{w}(0)=-\frac{1}{\da}<0$. Then, for sufficiently small $\da$, a shock forms at $T_{\ast}=O(\da\ep)$.}
   for $|\ga|\leq2$ and $\ga\neq(1,0,0)$. As a consequence, one obtains $|\fe_0|\leq\ep$ due to
   \bee
   0=\pa_{x_{\a}x_{\be}}\tilde{w}(0)=\pa_{\tilde{x}_{\a}\tilde{x}_{\be}}\tilde{w}(0)-\frac{1}{\ep}\fe_{0\a\be}.
   \ee
   Then, initially,
   \bee
   \dfrac{|N_0-e_1|}{\ep}\leq 1,\hs\dfrac{|T_0^{\nu}-e_{\nu}|}{\ep}\leq 1.
   \ee
   \end{remark}
  In the self-similar coordinates, one has the following bounds due to $y_1=\ep^{-\frac{3}{2}}x_1$ and $y_{\mu}=\ep^{-\frac{1}{2}}x_{\mu}$ for $t=-\ep$ (note that the initial support of $(W,Z,A,K)$ in self-similar coordinates is $supp=\{
   |y_1|\leq\ep^{-1},|\check{y}|\leq\ep^{-\frac{1}{3}}\}$):
   \begin{itemize}
   \item For all $y$, it holds that
         \bee
         |\pa^{\ga}Z|\leq\left\{\begin{array}{cc}
         \ep^{\frac{3}{2}},& \text{if}\ |\ga_1|\geq1,|\check{\ga}|=0,1,\\
         \ep,& \text{if}\ |\ga_1|=0,|\check{\ga}|\leq2,
         \end{array}\right.
         \ee
         \bee
         |\pa^{\ga}A|\leq\left\{\begin{array}{cc}
         \ep^{\frac{3}{2}},& \text{if}\ |\ga_1|\geq1,|\check{\ga}|=0,\\
         \ep,& \text{if}\ |\ga_1|=0,|\check{\ga}|\leq2,
         \end{array}\right.
         \ee
         \bee
         |\pa^{\ga}K|\leq\left\{\begin{array}{cc}
         \ep^2,& \text{if}\ \ga_1=1,|\check{\ga}|=0,1,\\
         \ep^{\frac{9}{4}}\eta^{-\frac{1}{15}}(y),& \text{if}\ \ga_1=2,|\check{\ga}|=0,\\
         \ep,& \text{if}\ \ga_1=0,|\check{\ga}|=1,2,
         \end{array}\right.
         \ee
         where $\eta(y)=1+y_1^2+y_{\mu}^6$.
   \item For $|y|\leq l$, it holds that
               \bee
               |\pa^{\ga}\tilde{W}(y,-\log\ep)|\leq\ep^{\frac{1}{8}},\hs \text{for} |\ga|=4,
               \ee
               while at $y=0$, one has
               \bee
               |\pa^{\ga}\tilde{W}(0,-\log\ep)|\leq\ep^{\frac{1}{2}-\frac{4}{2m-7}},\hs \text{for}\ |\ga|=3.
               \ee
   \item For $|y|\leq\mcl$, it holds that
         \begin{align}
         |\tilde{W}(y,-\log\ep)|&\leq\ep^{\frac{1}{10}}\eta^{\frac{1}{6}}(y),\\
         |\pa_1\tilde{W}(y,-\log\ep)|&\leq\ep^{\frac{1}{11}}\eta^{-\frac{1}{3}}(y),\\
         |\check{\nabla}\tilde{W}(y,-\log\ep)|&\leq\ep^{\frac{1}{12}}.
         \end{align}
   \item For $|y|\geq\mcl$, it holds that
         \begin{align}
         |W(y,-\log\ep)|&\leq(1+\ep^{\frac{1}{11}})\eta^{\frac{1}{6}}(y),\\
         |\pa_1W(y,-\log\ep)|&\leq(1+\ep^{\frac{1}{12}})\eta^{-\frac{1}{3}}(y),\\
         |\check{\nabla}W(y,-\log\ep)|&\leq\frac{3}{4}.
         \end{align}
   \item For the second derivatives and all $y$, it holds that
         \begin{align}
         |\pa^{\ga}W(y,-\log\ep)|&\leq\eta^{-\frac{1}{3}}(y),\ \text{for}\ \ga_1=1,\ |\ga|=2,\\
         |\pa^{\ga}W(y,-\log\ep)|&\leq\eta^{-\frac{1}{3}}(y)(\eta^{-\frac{1}{2}}(y)+\ep^{\frac{4}{5}}\eta^{\frac{2}{5}})^{\frac{1}{2}},\ \text{for}\ \ga_1=2,\ |\ga|=2,\\
         |\check{\nabla}^2W(y,-\log\ep)|&\leq\eta^{-\frac{1}{6}}(y).
         \end{align}
         Furthermore, at $y=0$, it follows from \eqref{eulerinitial} that
         \bee\label{contraints0}
         W(0,-\log\ep)=0,\ \pa_1W(0,-\log\ep)=-1,\ \check{\nabla}W(0,-\log\ep)=0,\ \nabla^2W(0,-\log\ep)=0.
         \ee
         For $m\geq 30$, the following energy bounds hold
         \bee
         \ep||W(\cdot,-\log\ep)||_{\dot{H}^m}^2+||Z,A,K(\cdot,-\log\ep)||_{\dot{H}^m}^2\leq\ep.
         \ee

   \end{itemize}
   \begin{remark}
   Moreover, the following bounds for the specific vorticity and the sound speed hold due to the definition \eqref{selfsimilarfluid}:
   \begin{align}
   ||\Omega(\cdot,-\log\ep)\cdot N_0||_{L^{\infty}}&\leq\ep^{\frac{1}{4}},\\
   ||\Omega(\cdot,-\log\ep)\cdot T^{\nu}_0||_{L^{\infty}}&\leq1,\\
   ||S(\cdot,-\log\ep)-\frac{\kappa_0}{2}||_{L^{\infty}}&\leq\ep^{\frac{1}{7}}.
   \end{align}
   \end{remark}
   \subsection{\textbf{The Bootstrap assumptions}}
   \hs For convenience, we introduce the following notations.
   \begin{itemize}
   \item $A$ is a lower order term (l.o.t) compared with $B$ means $|A(y,s)|=\ep^{\a}e^{-\be s}|O(B(y,s))|$ where $\a,\be\geq0$ and $\a^2+\be^2\neq 0$.
   \item Define the function $\ta(\a,\be)(y,s)$ (simply written as $\ta(\a,\be)$) to be $\ta(\a,\be):=e^{-\a s}\eta^{-\be}(y)$ where $\eta(y)=1+|y_1|^2+|\check{y}|^6$;
   \end{itemize}
  Then, we assume the following Bootstrap assumptions.
   \begin{itemize}
   \item Bootstrap assumptions on modulation variables. Assume that
         \begin{align}
         &|\dot{\kappa}(t)+\be_1ae^{-\be_1a(t+\ep)}\kappa_0|\leq \ep^{\frac{1}{6}} ,\hs |\dot{\tau}(t)|\leq \be_1ae^{-s}+Me^{-s},\label{bsmodulation1}\\
          &|\dot{\xi}(t)|\leq M^{\frac{1}{4}},\hs |\dot{\breve{n}}(t)|\leq M^2\ep^{\frac{1}{2}},
         \hs |\dot{\fe}(t)|\leq M^2,\label{bsmodulation2}
         \end{align}
         \bee\label{bsmodulation3}
         \dfrac{1}{2}\kappa_0\leq\kappa(t)\leq2\kappa_0,\hs |\tau(t)|\leq M\ep,\hs |\xi(t)|\leq M^{\frac{1}{4}}\ep, \hs|\breve{n}(t)|\leq M^2\ep^{\frac{3}{2}},
         \hs |\fe(t)|\leq M^2\ep,
         \ee
         for all $-\ep\leq t<T_{\ast}$, where $T_{\ast}$ is the shock time which is defined as
         \bee
         T_{\ast}=\tau(T_{\ast})\quad \text{or equivalently}\quad T_{\ast}=\int_{-\ep}^{T_{\ast}}\dot{\tau}(t) dt.
         \ee
         \begin{remark}
         It follows from the definition of $T_{\ast}$ and \eqref{bsmodulation1} that
         \bee
         T_{\ast}\les \ep.
         \ee
         Hence, the $L^{\infty}$ estimates of the modulation variables is a direct consequence of corresponding estimates on their derivatives. It follows from \eqref{bsmodulation1}-\eqref{bsmodulation3} that
         \bee\label{boundsdotQbetau}
         \bes
         |\dot{Q}(t)|&\leq 2M^2\ep^{\frac{1}{2}},\\
          |1-\be_{\tau}J|&\leq|(\be_{\tau}-1)(J-1)|+|\be_{\tau}-1|+|J-1|\leq\ep^{\frac{1}{2}},
          \end{split}
         \ee
         which will be used without mention.
         \end{remark}
   \item Bootstrap assumption on support of $(W,Z,A_{\nu},K)$. Assume that $(W,Z,A_{\nu},K)$ are supported in
         \bee\label{bssupport}
         Y(s):=\{|y_1|\leq2\ep^{\frac{1}{2}}e^{\frac{3}{2}s},|\check{y}|\leq2\ep^{\frac{1}{6}}\ep^{\frac{s}{2}}\}.
         \ee
         \begin{remark}\label{taprop}
         It follows from \eqref{bssupport} that $\eta^{\frac{1}{6}}(y)\leq\ep^{\frac{1}{4}}e^{\frac{s}{2}} $, and the following properties of the function $\ta(\a,\be)$ hold.
         \been[(1)]
         \item If $\a\geq 0$, then
         \bee
         \ta(\a,\be)\leq\ep^{\frac{\ga}{2}}e^{-(\a-\ga)s}\eta^{-(\frac{\ga}{3}+\be)}(y)=\ep^{\frac{\ga}{2}}\ta(\a-\ga,\frac{\ga}{3}+\be),
         \ee
         for all $0\leq\ga\leq\a$;
         \item if $\be\leq 0$, then \bee
         \ta(\a,\be)\leq\ep^{-\frac{3}{2}\ga}e^{-(\a+3\ga)s}\eta^{-(\be-\ga)}(y)=\ep^{-\frac{3}{2}\ga}\ta(\a+3\ga,\be-\ga),
         \ee
         for all $\be\leq\ga\leq0$.
         \een
         \end{remark}
   \item Bootstrap assumptions on $(W,Z,A_{\nu},K)$. Assume that the following bounds hold for all $y$.
         \bee\label{bspaZ}
         |\pa^{\ga}Z|\leq\left\{\begin{array}{cc}
         M^{\frac{1+|\check{\ga}|}{2}}e^{-\frac{3}{2}s},& \text{if}\ |\ga_1|\geq1,|\check{\ga}|=0,1,\\
         M\ep^{\frac{2-|\check{\ga}|}{2}}e^{-\frac{|\check{\ga}|}{2}},& \text{if}\ |\ga_1|=0,|\check{\ga}|\leq2,
         \end{array}\right.
         \ee
         \bee\label{bspaAmu}
         |\pa^{\ga}A|\leq\left\{\begin{array}{cc}
         Me^{-\frac{3}{2}s},& \text{if}\ |\ga_1|\geq1,|\check{\ga}|=0,\\
         M\ep^{\frac{2-|\check{\ga}|}{2}}e^{-\frac{|\check{\ga}|}{2}},& \text{if}\ |\ga_1|=0,|\check{\ga}|\leq2,
         \end{array}\right.
         \ee
         \bee\label{bspaK}
         |\pa^{\ga}K|\leq\left\{\begin{array}{cc}
         \ep^{\frac{1}{4}}e^{-\frac{3}{2}s},& \text{if}\ \ga_1=1,|\check{\ga}|=0,\\
         \ep^{\frac{1}{8}}e^{-\frac{13}{8}s},& \text{if}\ \ga_1=1,|\check{\ga}|=1,\\
         \ep^{\frac{1}{8}}e^{-2s}\eta^{-\frac{1}{15}}(y),& \text{if}\ \ga_1=2,|\check{\ga}|=0,\\
         \ep^{\frac{1}{8}}e^{-\frac{|\check{\ga}|}{2}s},& \text{if}\ \ga_1=0,|\check{\ga}|=1,2.
         \end{array}\right.
         \ee
         For the variables $W$ and $\tilde{W}$, we divide the spatial region as follows.
         \begin{itemize}
         \item For $|y|\leq l$, assume that
               \begin{align}
               |\pa^{\ga}\tilde{W}(y,s)|&\leq (\log M)^4\ep^{\frac{1}{10}}|y|^{4-|\ga|}+M\ep^{\frac{1}{4}}|y|^{3-|\ga|}, \quad \text{for}\ |\ga|\leq 3\label{pa3tildew},\\
               |\pa^{\ga}\tilde{W}(y,s)|&\leq\ep^{\frac{1}{10}}(\log M)^{|\check{\ga}|},\hs \text{for}\ |\ga|=4,\label{bapa4tildew}
               \end{align}
               while at $y=0$, assume that
               \bee\label{pa30tildew}
               |\pa^{\ga}\tilde{W}(0,s)|\leq\ep^{\frac{1}{4}},\hs \text{for}\ |\ga|=3.
               \ee
         \item For $|y|\leq\mcl$, assume that
         \begin{align}
         |\tilde{W}(y,s)|&\leq\ep^{\frac{1}{11}}\eta^{\frac{1}{6}}(y),\\
         |\pa_1\tilde{W}(y,s)|&\leq\ep^{\frac{1}{12}}\eta^{-\frac{1}{3}}(y),\\
         |\check{\nabla}\tilde{W}(y,s)|&\leq\ep^{\frac{1}{13}}.\label{bschecknablatildew}
         \end{align}
         \item For all $y$, assume that
         \bee\label{bsallyW}
         |\pa^{\ga}W(y,s)|\leq\left\{\begin{array}{cc}
         (1+\ep^{\frac{1}{20}})\eta^{\frac{1}{6}}(y), &\ \text{if}\ |\ga|=0,\\
         \tilde{\eta}(\frac{y}{2})\mathbf{1}_{|y|\leq\mcl}+2\eta^{-\frac{1}{3}}(y)\mathbf{1}_{|y|\geq\mcl},& \ \text{if}\ \ga_1=1,|\check{\ga}|=0,\\
         1, & \ \text{if}\ \ga_1=0,|\check{\ga}|=1,\\
         M^{\frac{2}{3}}\eta^{-\frac{1}{3}}(y),& \ \text{if}\ \ga_1=1,|\check{\ga}|=1,\\
         M^{\frac{1}{3}}\eta^{-\frac{1}{3}}(y)\fe^{\frac{1}{2}},&\ \text{if}\ \ga_1=2,|\check{\ga}|=0,\\
         M\eta^{-\frac{1}{6}}(y),&\ \text{if}\ \ga_1=0,|\check{\ga}|=2,
         \end{array}\right.
         \ee
         where $\fe=\ta(0,\frac{1}{2})+\ta(\frac{4}{5},-\frac{2}{5})$\footnote{The function $\fe$ is to control the growth of the entropy $K$, which vanishes in the isentropic case. The second term in $\fe$ can be chosen as $\ta(\a,\be)$ with $\a+3\be<\frac{2}{5}$ and $\a+\be>\frac{7}{15},$ which can be derived in recovering the bootstrap assumptions for $\pa_1^2K$ and $\pa_1^2W$.}.
         \begin{remark}\label{wtildew}
         Note that in $|y|$-small region, the bootstrap assumptions for $\tilde{W}$ are stronger than those in $|y|$-large region. One can obtain the corresponding estimates for $\pa^{\ga}W$ due to \eqref{estpabarW} and \eqref{pa3tildew}-\eqref{bschecknablatildew}, which is stronger than those in \eqref{bsallyW}.
         \end{remark}
         \end{itemize}

   \end{itemize}
   \subsection{\textbf{The main result}}
   \begin{thm}\label{mainthm}
   Assume the initial data in the physical variables are set as in Lemma\ref{initialdata}. Then,
   \begin{itemize}
   \item in the self-similar coordinates, the bootstrap assumptions \eqref{bsmodulation1}-\eqref{bsallyW} hold for all $(y,s)\in Y(s)\times [-\log\ep,+\infty)$, and the sound speed $S$, the specific vorticity $\Omega$, the velocity $U$, the pressure $P$ are smooth together with their derivatives for all $(y,s)\in Y(s)\times [-\log\ep,+\infty)$. Furthermore, it holds that
   \bee
   (W,Z,A_{\nu},K)\in C([-\log\ep,+\infty),H^m)\cap C^1([-\log\ep,+\infty),H^{m-1}),\hs m\geq 30,
   \ee
   with the energy bounds
   \bee
   e^{-s}||W(\cdot,s)||_{\dot{H}^m}^2+||Z(\cdot,s),A(\cdot,s),K(\cdot,s)||_{\dot{H}^m}^2\leq16k_0^2\lam^{-m}\ep^{-1}e^{-2s}+e^{-s}(1-e^{-s}\ep^{-1})M^{4m},
   \ee
   for all $s\geq-\log\ep$ and some constant $\lam\in(0,1)$.
   \item In the physical variables, there are following two possibilities.
   \been[(1)]
   \item If $\frac{1}{\ep}\geq 2\be_1a$, then the damping effect is weak enough to allow for the formation of a shock. There exists a pair $(T_{\ast},\xi_{\ast})$, which can be computed explicitly with $T_{\ast}=O(\ep),\xi_{\ast}=O(\ep)$, such that a point shock forms at this point. Furthermore, the shock time is shifted while the shock location and blow up direction remain same compared with the undamped case (see section\ref{section6.1} for more details).
   \item If $\frac{1}{\ep}\leq\frac{\be_1a}{2}$, then the damping effect is strong enough to prevent the formation of shock and a smooth classical solution to the compressible Euler system\eqref{euler} on $[-\ep,+\infty)$ is obtained.
   \een
   In case $(1)$, the following additional results hold.
   \begin{itemize}
   \item[$\cdot$] The Jacobian of the transformation between the physical variables (Cartesian coordinates) and the self-similar coordinates become $0$ at $T_{\ast}$, never vanishing everywhere else.
   \item[$\cdot$] The first derivative of Riemann invariants $\tilde{w}$ along blow up direction (i.e. along $N$) blows up like $-\dfrac{1}{T_{\ast}-t}$ at shock point while being bound everywhere else. Precisely, it holds that
                \begin{align}
                &\lim_{t\to T_{\ast}}N\cdot\nabla_{\tilde{x}}\tilde{w}(\xi(t),t)=-\infty,\\
                &\lim_{t\to T_{\ast}}N\cdot\nabla_{\tilde{x}}(\tilde{u}\cdot N)(\xi(t),t)=\lim_{t\to T_{\ast}}N\cdot\nabla_{\tilde{x}}\tilde{\p}(\xi(t),t)=-\infty,\\
                &\lim_{t\to T_{\ast}}N\cdot\nabla_{\tilde{x}}\tilde{w}(x,t)\leq |x|^{-\frac{2}{3}},\hs x\neq 0.
                \end{align}
   \item[$\cdot$] The other quantities are bounded as follows.
                \begin{align}
                &\sup_{[-\ep,T_{\ast})}\left(||\tilde{u}\cdot N-\frac{1}{2}\kappa_0||_{L^{\infty}}+||\tilde{u}\cdot T^{\nu}||_{L^{\infty}}+\ep^{-\frac{1}{8}}||\tilde{\sigma}-\frac{1}{2}\kappa_0||_{L^{\infty}}+||\zeta||_{L^{\infty}}+||\tilde{k}||_{L^{\infty}}\right)\les 1,\\
                &\sup_{[-\ep,T_{\ast})}\left(||T^{\nu}\cdot\nabla_{\tilde{x}}\tilde{\p}||_{L^{\infty}}+||T^{\nu}\cdot\nabla_{\tilde{x}}\tilde{u}||_{L^{\infty}}
                +||N\cdot\nabla_{\tilde{x}}(\tilde{u}\cdot T^{\nu})||_{L^{\infty}}+\ep^{-\frac{1}{8}}||\nabla_{\tilde{x}}\tilde{k}||_{L^{\infty}}\right)\les 1.
                \end{align}
   \item[$\cdot$] If $a>0$, then the damping effect leads to the dissipation of the vorticity while for $a<0$, the anti-damping effect leads to an increase in the vorticity.
       \item[$\cdot$] The normal vector $N(t)$ become horizontal in $O(\ep)$ scale as $t\to T_{\ast}$, that is
                \bee
                \dfrac{|N(t)-e_1|}{\ep}\to 0\ \text{as} \ t\to T_{\ast}.
                \ee
   \end{itemize}
   In case $(2)$, it is shown that $\tau(t)>t$ and $\dot{\tau}(t)>1$ for all $t\geq-\ep$. Then, the Jacobian of the transformation $\dfrac{\pa(y,s)}{\pa(x,t)}=\dfrac{1}{(\tau(t)-t)^{\frac{7}{2}}}\dfrac{1}{\be_{\tau}}$ never vanishes for all $t\geq-\ep$. This implies the fluid variables $(u,\p,\sigma)$ are bounded together with their first derivatives. In particular, compared to the case $(1)$, it holds that
   \bee
   |N\cdot\nabla_{\tilde{x}}(\tilde{u}\cdot N)(\xi(t),t)|=|\frac{1}{\tau(t)-t}|<\infty,\text{for all} t\in[-\ep,+\infty).
   \ee
   \end{itemize}
   \end{thm}
   The proof of the main theorem can be illustrated in the following picture:

\begin{figure}[h]
  \centering
  \includegraphics[width=0.8\textwidth]{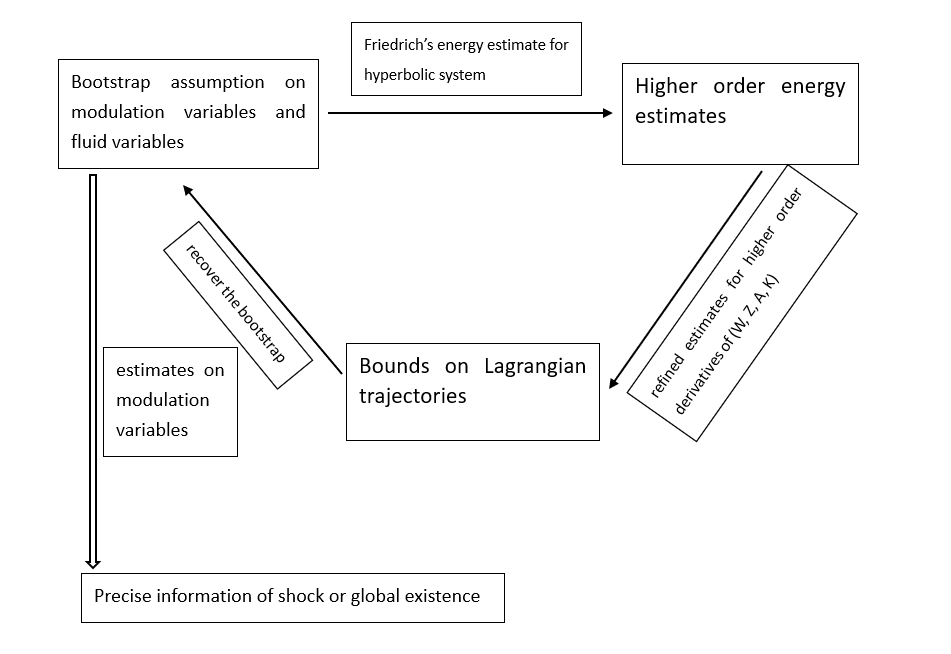}\\ 
   \caption{Main structure of the proof} 
   \label{mainstructure} 
\end{figure}

   \begin{remark}(The blow-up quantity)
   To see which quantities blow up as $t$ approaches $T_{\ast}$, one computes the following limits:
   \bee
   \lim_{t\to T_{\ast}}N\cdot\nabla_{\tilde{x}}\tilde{w}(\xi(t),t),\ \lim_{t\to T_{\ast}}T^2\cdot\nabla_{\tilde{x}}\tilde{w}(\xi(t),t),\
   \lim_{t\to T_{\ast}}N\cdot\nabla_{\tilde{x}}z(\xi(t),t),\ \lim_{t\to T_{\ast}}N\cdot\nabla_{\tilde{x}}\tilde{w}(x,t)
   \ee
   and observes that only the first quantity blows up, while the others remain bounded. Note that the above $z$ can be replaced by $a_v,k$ as well.
   \been[(1)]
   \item \begin{align*}
         N\cdot\nabla_{\tilde{x}}\tilde{w}(\xi(t),t)&=N\cdot\left(\dfrac{\pa}{\pa x_1},\dfrac{\pa}{\pa x_2}-f_{,2}\dfrac{\pa}{\pa x_1},\dfrac{\pa}{\pa x_3}-f_{,3}\dfrac{
         \pa}{\pa x_1}\right)w(\xi(t),t)\\
         &=J\dfrac{\pa}{\pa x_1}w|_{x=0}+N_{\mu}\dfrac{\pa}{\pa x_{\mu}}w|_{x=0}\\
         &=\dfrac{\pa}{\pa x_1}w|_{x=0}=e^{s}\dfrac{\pa}{\pa y_1}W|_{y=0}=e^s=\dfrac{-1}{\tau(t)-t}.
         \end{align*}
         Therefore,
         \bee\label{blowupquantity}
         \lim_{t\to T_{\ast}}N\cdot\nabla_{\tilde{x}}\tilde{w}(\xi(t),t)=-\infty.
         \ee
   \item \begin{align*}
         T^2\cdot\nabla_{\tilde{x}}\tilde{w}(\xi(t),t)&=T^2\cdot\left(\dfrac{\pa}{\pa x_1},\dfrac{\pa}{\pa x_2}-f_{,2}\dfrac{\pa}{\pa x_1},\dfrac{\pa}{\pa x_3}-f_{,3}\dfrac{
         \pa}{\pa x_1}\right)w(\xi(t),t)\\
         &=0\cdot\dfrac{\pa}{\pa x_1}w|_{x=0}+T^2_2\dfrac{\pa}{\pa x_{2}}w|_{x=0}+T^2_3\dfrac{\pa}{\pa x_3}w_{x=0}\\
         &=\dfrac{\pa}{\pa x_2}w|_{x=0}=\dfrac{\pa}{\pa y_2}W|_{y=0}\leq 1.
         \end{align*}
   \item Similarly,
         \bee
         N\cdot\nabla_{\tilde{x}}z(\xi(t),t)=\cdots=e^{\frac{3}{2}s}\dfrac{\pa}{\pa y_1}Z|_{y=0}\leq1.
         \ee
   \item \begin{align*}
         N\cdot \nabla_{\tilde{x}}\tilde{w}(x,t)&=J\dfrac{\pa}{\pa x_1}w(x,t)+N_{\mu}\dfrac{\pa}{\pa x_{\mu}}w(x,t)\\
         &=Je^{s}\dfrac{\pa}{\pa y_1}W(y,s)+N_{\mu}\dfrac{\pa}{\pa x_{\mu}}W(y,s)\\
         &\leq Je^sy_1^{-\frac{2}{3}}+1\leq 2x_1^{-\frac{2}{3}}.
         \end{align*}
   \een
   \end{remark}
   \hs We conclude this section by stating the following Sobolev-type inequalities, which can be verified directly.
   \begin{lem}\label{GNS}
   For $u: \mbr\to \mbr$, $1\leq q,r\leq\infty$, and $\dfrac{j}{m}\leq\a\leq1$ where $j,m\in\mathcal{N}$, if
   \bee
   \dfrac{1}{p}=\dfrac{j}{d}+\a\left(\dfrac{1}{r}-\dfrac{m}{d}\right)+\dfrac{1-\a}{q},
   \ee
   then
   \bee
   ||D^ju||_{L^p}\leq C||D^mu||_{L^r}^{\a}||u||_{L^q}^{1-\a}.
   \ee
   \end{lem}
    \hs The case $p=r=2,q=\infty,d=3$ will be used, which yields
    \bee
    ||u||_{\dot{H}^j}\les||u||_{\dot{H}^m}^{\a}||u||_{L^{\infty}}^{1-\a},
    \ee
    for $u\in H^m(\mbr^3)$ with compact support and $\a=\dfrac{2j-3}{2m-3}$.\\
    \hs It follows from Lemma\ref{GNS} that the following Lemma holds.
    \begin{lem}\label{sobolev2}
    Let $m\geq4$ and $0\leq l\leq m-3$. Then for $a+b=1-\dfrac{1}{2m-4}\in(0,1)$ and $q=\dfrac{6(2m-3)}{2m-1}$, it holds that
    \bee
    ||D^{2+l}\fe D^{m-1-l}\fai||_{L^2}\les ||D^m\fe||^{a}_{L^2}||D^m\fai||_{L^2}^b||D^2\fe||_{L^q}^{1-a}||D^2\fai||_{L^q}^{1-b}.
    \ee
    \end{lem}
    \begin{lem}\label{moser}
    For $\fe,\fai\in H^m(\mbr^3)$  with compact support, it holds that
    \bee
    ||\fe\fai||_{\dot{H}^m}\les||\fe||_{L^{\infty}}||\fai||_{\dot{H}^m}+||\fe||_{\dot{H}^m}||\fai||_{L^{\infty}}.
    \ee
    \end{lem}
   \section{\textbf{Preliminary estimates in self-similar coordinates under Bootstrap assumptions}}\label{section5}
   \hs First, we state the following homogeneous Sobolev energy estimates for $(W,Z,A_{\nu},K)$:
   \begin{prop}\label{energy}
   For some integers $m\geq 30$ and for some constant $\lam=\lam(m)\in(0,1)$, it holds that for all $s\geq-\log\ep$
   \bee
   e^{-s}||W(\cdot,s)||_{\dot{H}^m}^2+||Z(\cdot,s),A(\cdot,s),K(\cdot,s)||_{\dot{H}^m}^2\leq16k_0^2\lam^{-m}\ep^{-1}e^{-2s}+e^{-s}(1-e^{-s}\ep^{-1})M^{4m}.
   \ee
   \end{prop}
   The proof for this proposition will be given later, which is solely dependent on the bootstrap assumptions and standard Freidrich's energy estimates for the symmetric hyperbolic system. We mention the estimates here because the $(W,Z,A_{\nu},K)$ system loses one derivative. To derive the estimates for the higher order derivatives of forcing terms, it is necessary to obtain high order derivative estimates for $(W,Z,A_{\nu},K)$ using standard Sobolev inequalities.
   \begin{lem}\label{highorderwzak}
   For integer $m$ sufficiently large, the following refined estimates for higher order derivatives for $(W,Z,A,K)$ hold where the implicit constants are independent of $M$:
   \bee\label{highordera}
   |\pa^{\ga}A(y,s)|\les\left\{\begin{array}{cc}
   \ta(\frac{3}{2}-\frac{2|\ga|-1}{2m-5},0),&\quad  \text{if}\ \ga_1\geq 1, |\ga|=2,3,\\
   \ta(1-\frac{|\ga|-1}{2m-7},0),&\quad \text{if}\ |\ga|=3,4,5.
   \end{array}\right.
   \ee
   \bee\label{highorderz}
   |\pa^{\ga}Z(y,s)|\les\left\{\begin{array}{cc}
   \ta(\frac{3}{2}-\frac{3}{2m-7},0),&\quad  \text{if}\ \ga_1\geq 1, |\ga|=3,\\
   \ta(1-\frac{|\ga|-1}{2m-7},0),&\quad \text{if}\ |\ga|=3,4,5.
   \end{array}\right.
   \ee
    \bee\label{highorderw}
   |\pa^{\ga}W(y,s)|\les\left\{\begin{array}{cc}
   \ta(-\frac{2}{2m-7},\frac{1}{3}),&\quad \text{if}\ \ga_1=1, |\ga|=3,\\
   \ta(-\frac{1}{2m-7},\frac{1}{6}),&\quad \text{if}\ \ga_1=0,|\ga|=3,\\
   \ta(-\frac{3}{2m-7},\frac{1}{3})(\ta(0,\frac{1}{2})+\ta(\frac{2}{3},\frac{1}{2}))^{\frac{1}{2}},&\quad \text{if}\ \ga_1\geq 2,|\ga|=3.
   \end{array}\right.
   \ee
    \bee\label{highorderk}
   |\pa^{\ga}K(y,s)|\les\left\{\begin{array}{cc}
   \ta(\frac{13}{8}-\frac{9}{4(2m-7)},0),&\quad \text{if}\ \ga_1=1, |\ga|=3,\\
   \ta(2-\frac{4}{2m-7},\frac{1}{15}),&\quad \text{if}\ \ga_1\geq2,|\ga|=3,\\
   \ta(1-\frac{|\ga|-2}{2m-7},0),&\quad \text{if}\ |\ga|=3,4,5.
   \end{array}\right.
   \ee
   \end{lem}
   \begin{pf}
   The proof for the first and second cases in \eqref{highorderk} will be given and the proofs for \eqref{highordera}, \eqref{highorderz} and \eqref{highorderw} are similar. Taking $p=q=\infty$ in Lemma\ref{GNS} yields for any $u$
   \bee\label{GNS2}
   ||D^ju||_{L^{\infty}}\leq C||u||_{\dot{H}^{k}}^{\a}\cdot||u||_{L^{\infty}}^{1-\a},
   \ee
   where $\a=\dfrac{2j}{2k-3}$.
   \begin{itemize}
   \item For the first case in \eqref{highorderk}, take $u=\pa_1\nabla K$ in \eqref{GNS2} and then $\a=\frac{2}{2m-7}$. Then, it follows from the bootstrap assumptions and Proposition\ref{energy} that
         \begin{align}
         ||\pa^{\ga}K||_{L^{\infty}}&\leq C||\pa_1\nabla K||_{\dot{H}^{m-2}}^{\a}||\pa_1\nabla K||_{L^{\infty}}^{1-\a}\\
         &\leq\ep^{\frac{1}{8}(1-\a)}\ta(\frac{13}{8}-\frac{9}{8}\a,0)\leq\ta(\frac{13}{8}-\frac{9}{4(2m-7)},0).
         \end{align}
   \item For the second case in \eqref{highorderk}, take $u=\eta^{\frac{1}{15}}\pa_1^2K$ in \eqref{GNS2}. Note that
         \begin{align*}
         |\pa_1^2\nabla K\eta^{\frac{1}{15}}|&=|\nabla(\eta^{\frac{1}{15}}\pa_1^2K)|+|\pa_1^2K||\nabla\eta^{\frac{1}{15}}|
         \leq|\nabla(\eta^{\frac{1}{15}}\pa_1^2K)|+\ep^{\frac{1}{8}}\ta(2,\frac{1}{5}).
         \end{align*}
        Then, it suffices to show $|\nabla u|\les\ta(2-\frac{4}{2m-7},0)$. It follows from \eqref{GNS2} and the bootstrap assumptions that
         \bee
         ||\nabla(\eta^{\frac{1}{15}}\pa_1^2K)||_{L^{\infty}}\les||\eta^{\frac{1}{15}}\pa_1^2K||_{\dot{H}^{m-2}}^{\a}||\eta^{\frac{1}{15}}\pa_1^2K||_{L^{\infty}}^{1-\a}
         \leq||\eta^{\frac{1}{15}}\pa_1^2K||_{\dot{H}^{m-2}}^{\a}\ep^{\frac{1}{10}}\ta(2-\frac{4}{2m-7},0),
         \ee
         where $\a=\frac{2}{2m-7}$. By Moser inequality, it holds that
         \begin{align*}
         ||\eta^{\frac{1}{15}}\pa_1^2K||_{\dot{H}^{m-2}}&\leq||K||_{\dot{H}^{m}}||\eta^{\frac{1}{15}}||_{L^{\infty}}+||\eta^{\frac{1}{15}}||_{\dot{H}^{m-2}}||\pa_1^2K||_{L^{\infty}
         }\\
         &\leq\ta(\frac{1}{2},-\frac{1}{15})+\ta(2,0)\leq\ta(\frac{1}{10},0),
         \end{align*}
         due to Remark\ref{taprop}. Hence,
         \bee
         ||\nabla(\eta^{\frac{1}{15}}\pa_1^2K)||_{L^{\infty}}\leq\ep^{\frac{1}{10}}\ta(2-\frac{4}{2m-7},0).
         \ee
   \end{itemize}
   \end{pf}
   Due to the definition $U\cdot N=\dfrac{1}{2}(e^{-\frac{s}{2}}W+Z+\kappa)$ and $S=\dfrac{1}{2}(e^{-\frac{s}{2}}W-Z+\kappa)$, the following lemma holds by the bootstrap assumptions\eqref{bsmodulation1}-\eqref{bsallyW}.
   \begin{lem}\label{UNSestimates}
   For all $y\in Y(s)$ and $s\geq-\log\ep$, it holds that
   \bee
   |\pa^{\ga}U\cdot N|+|\pa^{\ga}S|\les\left\{\begin{array}{cc}
   M^{\frac{1}{4}},&|\ga|=0,\\
   M^{\frac{1}{3}}\ta(\frac{1}{2},\frac{1}{3}),&\ga=(1,0,0),\\
   \ta(\frac{1}{2},0),& |\ga_1|=0,|\breve{\ga}|=1,\\
   M^{\frac{2}{3}}\ta(\frac{1}{2},\frac{1}{3}),& |\ga_1|=1, |\check{\ga}|=1,\\
   M^{\frac{2}{3}}\ta(\frac{1}{2},\frac{1}{3}),& \ga=(2,0,0),\\
   M\ta(\frac{1}{2},\frac{1}{6}),& |\ga_1|=0,|\check{\ga}|=2,\\
   \ta(\frac{1}{2}-\frac{3}{2m-7},\frac{1}{3}),& |\ga_1|=1,|\check{\ga}|=2,\\
   \ta(\frac{1}{2}-\frac{1}{2m-7},\frac{1}{6}),& |\ga_1|=0,|\check{\ga}|=3,\\
   \ta(\frac{1}{2}-\frac{3}{2m-7},\frac{1}{3})\fe^{\frac{1}{2}},& |\ga_1|\geq 2, |\ga|=3.
   \end{array}\right.
   \ee
   \end{lem}
   \hs Similar argument leads to the following lemma.
   \begin{lem}\label{prelimi}
   For all $y\in Y(s)$ and $n\geq 1$, it holds that
   \bee
   |f|\les\ep^{\frac{1}{4}}e^{-s}|\check{y}|^2,\quad |\dot{f}|\les M^2e^{-s}|\check{y}|^2,
   \ee
   \begin{equation}
   \bes
   &|\breve{\nabla}^nf|+|\breve{\nabla}^n\dot{f}|+|\breve{\nabla}^n(N-N_0)|+|\breve{\nabla}^n\dot{N}|+|\breve{\nabla}^n(T^{\nu}-T^{\nu}_0)|+|\breve{\nabla}^n\dot{T}^{\nu}|\\
   &+|\breve{\nabla}^n(J-1)|+|\breve{\nabla}^n(J^{-1}-1)|\leq \ep^{\frac{1}{4}}\ta(\frac{n}{2},0),
   \end{split}
   \end{equation}
   \bee
   |\pa^{\ga}V|\les\left\{\begin{array}{cc}
   M^{\frac{1}{4}}, & |\ga|=0,\\
   M^2\ep^{\frac{1}{2}}\ta(\frac{3}{2},0), & \ga=(1,0,0),\\
   M^2\ep^{\frac{1}{2}}\ta(-\frac{1}{2},0), & |\ga_1|=0,|\breve{\ga}|=1,\\
   M^4\ep^{\frac{3}{2}}\ta(-1,0), & \ga=(2,0,0),\\
   0, & \text{else.}
   \end{array}\right.
   \ee
   \end{lem}
  The following bounds for the sound speed and the density can be obtained.
   \begin{lem}
   There exists a constant $C=C(\a,\kappa_0)>1$ such that
   \begin{align}
   \dfrac{1}{C}\leq||\p(\cdot,t)||_{L^{\infty}}&\leq C,\\
   ||\sigma-\frac{\kappa_0}{2}||_{L^{\infty}}=||S-\frac{\kappa_0}{2}||_{L^{\infty}}&\leq\ep^{\frac{1}{8}}.
   \end{align}
   \end{lem}
   \begin{pf}
   It follows from the definition\eqref{selfsimilarfluid} that
   \begin{align*}
   |S-\dfrac{\kappa_0}{2}|\leq|\dfrac{\kappa-\kappa_0}{2}|+|\dfrac{1}{2}(e^{-\frac{s}{2}}W-Z)|\leq\ep^{\frac{1}{8}},
   \end{align*}
   and hence,
   \bee
   |\p^{\a}-\a\dfrac{\kappa_0}{2}|=|\a\sigma e^{-\frac{k}{2}}-\a\dfrac{\kappa_0}{2}|\leq\ep^{\frac{1}{8}}.
   \ee
   \end{pf}
   \subsection{\textbf{Estimates on transport terms and forcing terms}}
   \begin{lem}\label{dampingtermsestimates}
   For $\ep>0$ sufficiently small and $y\in Y(s)$, it holds that
   \bee
   |\pa^{\ga}G_W|\les\left\{\begin{array}{cc}
   Me^{-(\frac{1}{2}-\frac{4}{2m-7})s}+M^{\frac{1}{2}}e^{-s}|y_1|+\ep^{\frac{1}{3}}|\breve{y}|, &|\ga|=0,\\
   M^2\ep^{\frac{1}{2}}, &\ga_1=0,|\breve{\ga}|=1,\\
   Me^{-\frac{s}{2}}, &|\ga_1|\leq 1, |\ga|=2,\\
   M^{\frac{1}{2}}e^{-s}, &\ga=(2,0,0).
   \end{array}\right.
   \ee
   \bee
   |\pa^{\ga}[G_Z+(1-\be_2)e^{\frac{s}{2}}\kappa_0]|+
   |\pa^{\ga}[G_U+(1-\be_1)e^{\frac{s}{2}}\kappa_0]|\les\left\{\begin{array}{cc}
   \ep^{\frac{1}{2}}\ta(-\frac{1}{2},0), &|\ga|=0,\\
   M^{2}\ep^{\frac{1}{2}}, &\ga_1=0,|\breve{\ga}|=1,\\
   M\ta(\frac{1}{2},0), &\ga_1\leq 1, |\ga|=2,\\
   M^{\frac{1}{2}}\ta(1,0), &\ga=(2,0,0).
   \end{array}\right.
   \ee
   \bee
   |\pa^{\ga}h_W|+|\pa^{\ga}h_Z|+|\pa^{\ga}h_U|\les\left\{\begin{array}{cc}
   \ta(\frac{1}{2},0), &|\ga|=0,\\
   \ta(1,0), &\ga_1=0, |\breve{\ga}|=1,\\
   \ta(1,\frac{1}{6}), &\ga_1\leq 1, |\ga|=2\\
   \ta(2-\frac{3}{2m-5},0), &\ga=(2,0,0).
   \end{array}\right.
   \ee
   \end{lem}
   \begin{pf}
   Note that
   \begin{align*}
   G_W&=\be_{\tau}e^{\frac{s}{2}}[-\dot{f}+J(\kappa+\be_2Z+2\be_1V\cdot N)]\\
   &=\be_{\tau}e^{\frac{s}{2}}[-\dot{f}+J(\underbrace{\kappa+\be_2Z^0-R_{ji}\dot{\xi}_j}_{G_W^0}+\be_2(Z-Z^0)+2\be_1(V\cdot N+R_{ji}\dot{\xi}_j))].
   \end{align*}
   Then\footnote{The estimate for $G_W^0$ can be found in \eqref{GW0hW0}, which only relies on the bootstrap assumptions.},
   \begin{align*}
   |G_W|&\les e^{\frac{s}{2}}\left(|\dot{f}|+|G_W^0|+|y_1||\pa_1Z|+|V\cdot N+R_{ji}\dot{\xi}_j|\right)\\
   &\les Me^{-(\frac{1}{2}-\frac{4}{2m-7})s}+M^{\frac{1}{2}}e^{-\frac{s}{2}}|y_1|+\ep^{\frac{1}{2}}|\breve{y}|.
   \end{align*}
   For $G_Z$, it follows from \eqref{defGWGZGU} that
   \bee
   G_W-G_Z-(1-\be_2)e^{\frac{s}{2}}\kappa_0=e^{\frac{s}{2}}(1-\be_2)[(\be_{\tau}J-1)\kappa+(\kappa-\kappa_0)-J\be_1Z],
   \ee
   and then the estimate for $G_Z$ follows from the bootstrap assumptions \eqref{bsmodulation1} and \eqref{bspaZ}. The estimate for $G_U$ is similar to $G_Z$. The estimates for the remaining term follows directly from the bootstrap assumptions, the definiton\eqref{defGWGZGU} and Lemma\ref{prelimi}.
   \end{pf}
   \begin{lem}\label{forcingtermsestimates}
   For all $s\in[-\log\ep,+\infty)$ and $y\in Y(S)$, the following bounds for the forcing terms $F_W^{\ga}$, $F_Z^{\ga}$ and $F_K^{\ga}$ hold.
   \bee\label{estimatesFWga}
   |F_W^{(\ga)}|\les\left\{\begin{array}{cc}
   M^{\frac{1}{4}}\ta(\frac{1}{2},0), &|\ga|=0,\\
   \ta(1,\frac{1}{15}), &\ga=(1,0,0),\\
   \ep^{\frac{1}{3}}\ta(0,\frac{5}{24}), &\ga_1=0, |\breve{\ga}|=1,\\
   M\ta(\frac{1}{2},\frac{1}{3})\fe^{\frac{1}{2}}+M^{\frac{1}{2}}\ta(1-\frac{4}{2m-7},\frac{1}{15}), & \ga=(2,0,0),\\
   M^{\frac{2}{3}}\ta(0,\frac{1}{3}), & |\ga_1|=1, |\breve{\ga}|=1,\\
   M^{\frac{2}{3}}\ta(0,\frac{1}{3})+M^{\frac{1}{2}}\ta(\frac{5}{8}-\frac{9}{4(2m-7)},0), &\ga_1=0, |\breve{\ga}|=2,
   \end{array}\right.
   \ee
   \bee
   |F_Z^{(\ga)}|\les\left\{\begin{array}{cc}
   M^{\frac{1}{4}}\ta(1,0), &|\ga|=0,\\
   \ep^{\frac{1}{8}}\ta(\frac{3}{2},\frac{1}{15}), &\ga=(1,0,0),\\
   \ep^{\frac{1}{8}}\ta(\frac{9}{8},0), &\ga_1=0, |\breve{\ga}|=1,\\
   M\ta(2-\frac{3}{2m-5},0), & \ga=(2,0,0),\\
   M\ta(\frac{3}{2}-\frac{4}{2m-7},\frac{1}{15})+M\ta(\frac{3}{2},0), & |\ga_1|=1, |\breve{\ga}|=1,\\
   M\ta(\frac{9}{8}-\frac{9}{4(2m-7)},0), &\ga_1=0, |\breve{\ga}|=2,
   \end{array}\right.
   \ee
   \bee
   |F_{A_{\nu}}^{(\ga)}|\les\left\{\begin{array}{cc}
   M^{\frac{1}{2}}\ta(1,0), &|\ga|=0,\\
   M\ta(1,\frac{1}{6}), &\ga_1=0, |\breve{\ga}|=1,\\
   M^{\frac{1}{4}}\ta(1-\frac{3}{2m-7},\frac{1}{6}), &\ga_1=0, |\breve{\ga}|=2,
   \end{array}\right.
   \ee
   \bee
   |F_K^{(\ga)}|\les\left\{\begin{array}{cc}
   \ep^{\frac{1}{8}}\ta(\frac{3}{2},\frac{1}{6}), &\ga=(1,0,0),\\
   \ep^{\frac{1}{8}}\ta(\frac{3}{2},0), &\ga_1=0, |\breve{\ga}|=1,\\
   M^{\frac{1}{3}}\ep^{\frac{1}{4}}\ta(\frac{3}{2},\frac{1}{3})\fe^{\frac{1}{2}}, & \ga=(2,0,0),\\
   \ep^{\frac{1}{8}}\ta(\frac{3}{2},\frac{1}{6}), & |\ga_1|=1, |\breve{\ga}|=1,\\
   \ep^{\frac{1}{8}}\ta(\frac{13}{8},0), &\ga_1=0, |\breve{\ga}|=2.
   \end{array}\right.
   \ee
   For $\tilde{F}_W^{(\ga)}$, the following bounds hold
   \bee
   |\tilde{F}_W^{(\ga)}|\les\left\{\begin{array}{cc}
   \ep^{\frac{1}{11}}\ta(0,\frac{2}{5}),&\ga=(1,0,0)\ \text{and}\ |y|\leq\mcl,\\
   \ep^{\frac{1}{12}}\ta(0,\frac{1}{3}),&\ga_1=0,|\check{\ga}|=1\ \text{and}\ |y|\leq\mcl,\\
   \ep^{\frac{1}{8}}+\ep^{\frac{1}{10}}(\log M)^{|\check{\ga}|-1},&|\ga|=4\ \text{and}\ |y|\leq l,
   \end{array}\right.
   \ee
   and for $|\ga|=3$, it holds that
   \bee\label{ftildew30}
   |\tilde{F}_W^{(\ga),0}|\les\ta(\frac{1}{2}-\frac{4}{2m-7},0).
   \ee
   \end{lem}
   \begin{pf}
   It follows from the Bootstrap assumptions and \eqref{FW} that
   \bee
   |\pa^{\ga}F_W|\les |\pa^{\ga}(ST^{\nu}_{\mu}\pa_{\mu}A_{\nu})|+e^{s}|\pa^{\ga}(JS^2\pa_1K)|+\text{l.o.ts},
   \ee
   and then
   \begin{align}
   |F_W^{(\ga)}|&\leq|\pa^{\ga}F_W|+\sum_{0\leq\be<\ga}(|\pa^{\ga-\be}G_W\pa^{\ga}\pa_1W|+|\pa^{\ga-\be}h_W^{\mu}\pa^{\be}\pa_{\mu}W|)+
   1_{|\ga|\geq 2}\be_{\tau}\sum_{|\be|=|\ga|-1,\be_1=\ga_1}|\check{\nabla}(JW)\pa^{\be}\pa_1W|\\
   &+1_{|\ga|\geq 3}\be_{\tau}\sum_{1\leq|\be|\leq|\ga|-2,\be\leq\ga}|\pa^{\ga-\be}(JW)\pa^{\be}\pa_1W|+\be_{\tau}
   \sum_{0\leq\be<\ga-1}|\pa^{\ga-\be}J\pa^{\be}W\pa_1W|+\be_1\be_{\tau}ae^{-\frac{s}{2}}\pa^{\ga}(U\cdot N)\\
   &\leq |\pa^{\ga}(ST^{\nu}_{\mu}\pa_{\mu}A_{\nu})|+e^{s}|\pa^{\ga}(JS^2\pa_1K)|+\sum_{0\leq\be<\ga}|\pa^{\ga-\be}G_W\pa^{\ga}\pa_1W|\\
   &+1_{|\ga|\geq 2}\be_{\tau}\sum_{|\be|=|\ga|-1,\be_1=\ga_1}|\check{\nabla}(JW)\pa^{\be}\pa_1W|+e^{-\frac{s}{2}}|\be_1\be_{\tau}a\pa^{\ga}(U\cdot N)|+\text{l.o.ts}.
   \end{align}
   Therefore, the following cases hold due to the bootstrap assumptions, Lemma\ref{highorderwzak},\ref{UNSestimates},\ref{prelimi} and Lemma\ref{dampingtermsestimates}.
   \begin{itemize}
   \item For $|\ga|=0,$
         \begin{align*}
         |F_W^{(\ga)}|&\leq|ST^{\nu}_{\mu}\pa_{\mu}A_{\nu}|+e^{s}|JS^2\pa_1K|+e^{-\frac{s}{2}}|\be_1\be_{\tau}a U\cdot N|\\
         &\leq \kappa_0e^{-\frac{s}{2}}+\ep^{\frac{1}{4}}e^{-s}\leq M^{\frac{1}{4}}\ta(\frac{1}{2},0);
         \end{align*}
   \item for $\ga_1=1,|\check{\ga}|=0$,
         \begin{align*}
         |F_W^{(\ga)}|&\leq|\pa_1(ST^{\nu}_{\mu}\pa_{\mu}A_{\nu})|+e^{s}|\pa_1(JS^2\pa_1K)|\\
         &+|\pa_1G_W\pa_1W|+e^{-\frac{s}{2}}|\be_1\be_{\tau}a\pa_1(U\cdot N)|\\
         &\leq M^{\frac{1}{4}}\ta(\frac{1}{3},\frac{1}{3})+M^{\frac{1}{2}}\ep^{\frac{1}{8}}
         \ta(1,\frac{1}{15})+M\ta(\frac{1}{2},\frac{1}{3})\\
         &\leq M^{\frac{1}{4}}\ta(\frac{1}{15},\frac{1}{3});
         \end{align*}
   \item for $\ga_1=0,|\check{\ga}|=1$,
         \begin{align*}
         |F_W^{(\ga)}|&\leq|\check{\nabla}(ST^{\nu}_{\mu}\pa_{\mu}A_{\nu})|+
         e^{s}|\check{\nabla}(JS^2\pa_1K)|\\
         &+e^{-\frac{s}{2}}|\be_1\be_{\tau}
         a\check{\nabla}(U\cdot N)|+|\check{\nabla}G_W\pa_1W|\\
         &\leq M^{\frac{5}{4}}\ta(1,0)+\ep^{\frac{1}{8}}\ta(\frac{5}{8},0)+\ep^{\frac{1}{3}}
         \ta(0,\frac{1}{3})\\
         &\leq\ep^{\frac{1}{3}}\ta(0,\frac{5}{24});
         \end{align*}
   \item for $\ga=(2,0,0)$,
         \begin{align*}
         |F_W^{(\ga)}|&\leq|\pa_1^2(ST^{\nu}_{\mu}\pa_{\mu}A_{\nu})|+e^s|\pa_1^2(JS^2
         \pa_1K)|+e^{-\frac{s}{2}}|\be_1\be_{\tau}a\pa_1^2(U\cdot N)|\\
         &+|\pa_1^2G_W\pa_1W|+|\pa_1G_W\pa_1^2W|\\
         &\leq\ep^{\frac{1}{3}}\ta(1,\frac{1}{3})+\ep\ta(\frac{4}{3},0)
         +M^{\frac{1}{2}}\ta(1-\frac{4}{2m-7},\frac{1}{15})\\
         &+M^{\frac{1}{2}}\ta(1,\frac{1}{3})+M\ta(\frac{1}{2},\frac{1}{3})\fe^{\frac{1}{2}}\\
         &\leq M\ta(\frac{1}{2},\frac{1}{3})\fe^{\frac{1}{2}}+M^{\frac{1}{2}}\ta(1-\frac{4}{2m-7},\frac{1}{15});
         \end{align*}
   \item for $\ga_1=1,|\check{\ga}|=1$,
         \begin{align*}
         |F_W^{(\ga)|}|&\leq|\pa_1\check{\nabla}(ST^{\nu}_{\mu}\pa_{\mu}A_{\nu})|
         +e^s|\pa_1\check{\nabla}(JS^2\pa_1K)|+e^{-\frac{s}{2}}|\be_{1}\be_{\tau}a
         \pa_1\check{\nabla}(U\cdot N)|\\
         &+|\pa_1G_W\pa_1\check{\nabla}W|+|\check{\nabla}G_W\pa_1^2W|+|\pa_1\check{\nabla}G_W
         \pa_1W|+|\check{\nabla}(JW)\pa_1\check{\nabla}W|\\
         &\leq M^{\frac{1}{4}}\ta(\frac{3}{2}-\frac{3}{2m-5},0)+M^{\frac{1}{2}}\ta(1-
         \frac{4}{2m-7},\frac{1}{15})\\
         &+M^{\frac{2}{3}}\ta(\frac{1}{2},\frac{1}{3})+\ep^{\frac{1}{4}}\ta(0,
         \frac{1}{3})+M^{\frac{2}{3}}\ta(0,\frac{1}{3})\\
         &\leq M^{\frac{2}{3}}\ta(0,\frac{1}{3});
         \end{align*}
   \item for $\ga_1=0,|\check{\ga}|=2$,
         \begin{align*}
         |F_W^{(\ga)}&\leq|\check{\nabla}^2(ST^{\nu}_{\mu}\pa_{\mu}A_{\nu})|
         +e^s|\check{\nabla}^2(JS^2\pa_1K)|+e^{-\frac{s}{2}}|\be_{1}\be_{\tau}a
         \check{\nabla}^2(U\cdot N)|\\
         &+|\check{\nabla}G_W\pa_1\check{\nabla}W|+|\check{\nabla}^2G_W
         \pa_1W|+|\check{\nabla}(JW)\pa_1\check{\nabla}W|\\
         &\leq  M^{\frac{1}{4}}\ta(1-\frac{2}{2m-7},0)+M^{\frac{1}{2}}\ta(
         \frac{5}{8}-\frac{9}{4(2m-7)},0) \\
         &+M\ta(\frac{1}{2},\frac{1}{3})+M^{\frac{2}{3}}\ta(0,\frac{1}{3})\leq
         M^{\frac{2}{3}}\ta(0,\frac{1}{3})+M^{\frac{1}{2}}\ta(
         \frac{5}{8}-\frac{9}{4(2m-7)},0).
         \end{align*}
   \end{itemize}
   Similarly, for the forcing terms of $(Z,A,K)$, it holds that
   \begin{align*}
   \begin{split}
   |F_Z^{(\ga)}|&\leq e^{-\frac{s}{2}}|\pa^{\ga}(ST^{\nu}_{\mu}\pa_{\mu}A_{\nu})|+
   e^{\frac{s}{2}}|\pa^{\ga}(S^2J\pa_1K)|+e^{-s}|\be_1\be_{\tau}a\pa^{\ga}(U\cdot N)|\\
   &+\sum_{0\leq\be<\ga}|\pa^{\ga-\be}G_Z\pa^{\be}\pa_1Z|+\sum_{
   0\leq\be<\ga}|\pa^{\ga-\be}(JW)\pa^{\be}\pa_1Z|+\text{l.o.ts},
   \end{split}\\
   \begin{split}
   |F_{A_v}^{(\ga)}|&\leq e^{-\frac{s}{2}}|\pa^{\ga}(ST^{\nu}_{\mu}\pa_{\mu}S)|+e^{-\frac{s}{2}}|\pa^{\ga}
   (S^2T^{\nu}_{\mu}\pa_{\mu}K|+e^{-s}|\be_1\be_{\tau}a\pa^{\ga}A_{\nu}|\\
   &+\sum_{0\leq\be<\ga}|\pa^{\ga-\be}G_U\pa^{\be}\pa_1A_{\nu}|+
   1_{|\ga|\geq2}\sum_{0\leq\be<\ga}|\pa^{\ga-\be}(JW)\pa^{\be}\pa_1A_{\nu}|+\text{l.o.ts},
   \end{split}\\
   |F_{K}^{(\ga)}|&\leq
   \sum_{0\leq\be<\ga}|\pa^{\ga-\be}G_U\pa^{\be}\pa_{\mu}K|+
   1_{|\ga|\geq2}\sum_{0\leq\be<\ga}|\pa^{\ga-\be}(JW)\pa^{\be}\pa_1K|+\text{l.o.ts}.
   \end{align*}
   Therefore, it follows from the bootstrap assumptions, Lemma\ref{highorderwzak},\ref{UNSestimates},\ref{prelimi} and Lemma\ref{dampingtermsestimates} that
   \begin{itemize}
   \item for $|\ga|=0$, $|F_Z|\leq M^{\frac{1}{4}}e^{-s}$;
   \item for $\ga_1=1,|\check{\ga}|=0$,
             \begin{align*}
             |F_Z^{(\ga)}|&\leq e^{-\frac{s}{2}}|\pa_1(ST_{\mu}^{\nu}\pa_{\mu}A_{\nu}|+e^{\frac{s}{2}}|S^2J\pa_1^2K|+e^{-s}|\be_1\be_{\tau}a\pa_1(U\cdot N)|\\
         &+|\pa_1G_Z\pa_1Z|+\text{l.o.ts}\\
         &\leq M\ta(2-\frac{3}{2m-5},0)+\ep^{\frac{1}{8}}\ta(\frac{3}{2},\frac{1}{15})
             +M\ta(2,0)\\
             &\leq\ep^{\frac{1}{8}}\ta(\frac{3}{2},\frac{1}{15});
             \end{align*}
   \item for $\ga_1=0,|\check{\ga}|=1$,
         \begin{align*}
         |F_Z^{(\ga)}|&\leq e^{-\frac{s}{2}}|S||T^{\nu}_{\mu}||\check{\nabla}\pa_{\mu}A_{\nu}|+
         e^{\frac{s}{2}}|JS^2||\pa_1\check{\nabla}K|\\
         &+|\check{\nabla}G_Z\pa_1Z|+e^{-s}|\be_1\be_{\tau}a\check{\nabla}(U\cdot N)|+\text{l.o.ts}\\
         &\leq M\ta(\frac{3}{2},0)+\ep^{\frac{1}{8}}\ta(\frac{9}{8},0)+\ep^{\frac{1}{4}}
         \ta(\frac{3}{2},0)\leq\ep^{\frac{1}{8}}\ta(\frac{9}{8},0);
         \end{align*}
   \item for $\ga_1=2,|\check{\ga}|=0$,
         \begin{align*}
         |F_Z^{(\ga)}|&\leq e^{-\frac{s}{2}}|ST^{\nu}_{\mu}||\pa_1^2\check{\nabla}A_{\nu}|+
         e^{\frac{s}{2}}|JS^2||\pa_1^3K|+e^{-s}|\be_1\be_{\tau}a\pa_1^2(U\cdot N)|\\
         &+|\pa_1^2G_Z\pa_1Z|+|\pa_1Z\pa_1^2Z|+|\pa_1^2(JW)\pa_1Z|+\text{l.o.ts}\\
         &\leq M\ta(2-\frac{3}{2m-5},0)+M\ta(\frac{3}{2}-\frac{3}{2m-5},0)
         \\
         &+M\ta(2,0)+M\ta(\frac{3}{2},\frac{1}{3})\leq M\ta(\frac{3}{2}-\frac{3}{2m-5},0);
         \end{align*}
   \item for $\ga_1=1,|\check{\ga}|=1$,
         \begin{align*}
         |F_Z^{(\ga)}|&\leq e^{-\frac{s}{2}}|ST^{\nu}_{\mu}||\pa_1\pa_{\mu}\check{\nabla}A_{\nu}|
         +e^{\frac{s}{2}}|JS^2||\pa_1^2\check{\nabla}K|\\
         &+|\check{\nabla}G_Z\pa_1^2Z|+|\check{\nabla}(JW)\pa_1\check{\nabla}
         Z|+l.o.t\\
         &\leq M\ta(2-\frac{5}{2m-5},0)+M\ta(\frac{3}{2}-\frac{4}{2m-7},\frac{1}{15})
         +M\ta(\frac{3}{2},0)\\
         &\leq M\ta(\frac{3}{2}-\frac{4}{2m-7},\frac{1}{15})+M\ta(\frac{3}{2},0);
         \end{align*}
   \item for $\ga_1=0,|\check{\ga}|=2$,
         \begin{align*}
         |F_Z^{(\ga)}|&\leq e^{-\frac{s}{2}}|ST^{\nu}_{\mu}||\pa_{\mu}\check{\nabla}A_{\nu}|
         +e^{\frac{s}{2}}|JS^2||\pa_1\check{\nabla}^2K|\\
         &+|\check{\nabla}G_Z\pa_1\check{\nabla}Z|+|\check{\nabla}(JW)||
         \pa_1\check{\nabla}Z|+\text{l.o.ts}\\
         &\leq \ta(\frac{3}{2}-\frac{2}{2m-7},0)+M\ta(\frac{9}{8}-\frac{9}{4(2m-7)},0)\\
         &+\ep^{\frac{1}{4}}\ta(\frac{3}{2},0)+M\ta(\frac{3}{2},\frac{1}{6})\leq
         M\ta(\frac{9}{8}-\frac{9}{4(2m-7)}).
         \end{align*}
   \end{itemize}
   And for $F_{A_{\nu}}^{(\ga)}$ and $F_K^{(\ga)}$, it holds that
   \begin{itemize}
   \item for $|\ga|=0$, $|F_{A_{\nu}}| \leq M^{\frac{1}{2}}\ta(1,0)$;
   \item for $\ga_1=1,|\check{\ga}|=0$, $|F_{K}^{(\ga)}|\leq|\pa_1h_U^{\mu}||\pa_{\mu}K|\leq\ep^{\frac{1}{8}}\ta(\frac{3}{2},\frac{1}{6})$;
   \item for $\ga_1=0,|\check{\ga}|=1$,
         \begin{align*}
         |F_{A_{\nu}}^{(\ga)}|&\leq e^{-\frac{s}{2}}|ST^{\nu}_{\mu}||\check{\nabla}\pa_{\mu}S|+e^{-\frac{s}{2}}|S^2T^{\nu}_{\mu}||\check{\nabla}\pa_{\mu}K|\\
         &+|\check{\nabla}G_U\pa_1 A_{\nu}|+\text{l.o.ts}\leq M\ta(1,\frac{1}{6})\\
         |F_K^{(\ga)}|&\leq|\check{\nabla}h_U^{\mu}||\pa_{\mu}K|+\text{l.o.ts}\leq\ep^{\frac{1}{8}}\ta(\frac{3}{2},0);
         \end{align*}
   \item for $\ga_1=2,|\check{\ga}|=0$,
         \begin{align*}
         |F_K^{(\ga)}|&\leq|\pa_1^2h_U^{\mu}\pa_{\mu}K|+|\pa_1h_U^{\mu}\pa_1\pa_{\mu}K|+|\pa_1^2(JW)\pa_1K|+\text{l.o.ts}\\
         &\leq M^{\frac{1}{3}}\ep^{\frac{1}{4}}\ta(\frac{3}{2},\frac{1}{3})\fe^{\frac{1}{2}};
         \end{align*}
   \item for $\ga_1=1,|\check{\ga}|=1$,
         \begin{align*}
         |F_K^{(\ga)}|&\leq|\pa_1\check{\nabla}h_U^{\mu}\pa_{\mu}K|+|\pa_1h_U^{\mu}\pa_{\mu}\check{\nabla}K|+|\check{\nabla}h_U^{\mu}\pa_1\pa_{\mu}K|\\
         &+|\pa_1\check{\nabla}(JW)\pa_1K|+|\pa_{\mu}(JW)\pa_1^2K|+\text{l.o.ts}\\
         &\leq\ep^{\frac{1}{8}}\ta(\frac{3}{2},\frac{1}{6});
         \end{align*}
   \item for $\ga_1=0,|\check{\ga}|=2$,
         \begin{align*}
         \begin{split}
         |F_{A_{\nu}}^{(\ga)}|&\leq e^{-\frac{s}{2}}|ST^{\nu}_{\mu}||\check{\nabla}^2\pa_{\mu}S|+e^{-\frac{s}{2}}|S^2T^{\nu}_{\mu}|\check{\nabla}^2\pa_{\mu}K|\\
         &+|\check{\nabla}G_U||\pa_1\check{\nabla}A_{\nu}|+|\check{\nabla}(JW)||\pa_1\check{\nabla}A_{\nu}|+\text{l.o.ts}\\
         &\leq M^{\frac{1}{4}}\ta(1-\frac{1}{2m-7},\frac{1}{6});
         \end{split}\\
         |F_K^{(\ga)}|&\leq|\check{\nabla}h_U^{\mu}||\check{\nabla}\pa_{\mu}K|+|\check{\nabla}(JW)||\pa_1\check{\nabla}K|+\text{l.o.ts}\leq\ep^{\frac{1}{8}}\ta(\frac{13}{8},0).
         \end{align*}
   \end{itemize}
   For $|\ga|=1$ and $|y|\leq\mcl$, it follows from \eqref{tildeFWga} and Lemma\ref{dampingtermsestimates} that
   \begin{align*}
   |\tilde{F}_W^{(\ga)}|&\les|\nabla\tilde{F}_W|+|\nabla G_W||\pa_1\tilde{W}|+
   |\nabla(J\pa_1\bar{W})||\tilde{W}|+\text{l.o.ts}\\
   &\les\left\{\begin{aligned}
   &|\pa_1\tilde{F}_W|+\ep^{\frac{1}{10}}\ta(0,\frac{1}{3}),\\
   &|\check{\nabla}\tilde{F}_W|+\ep^{\frac{1}{10}}\ta(0,\frac{1}{3}).
   \end{aligned}\right.
   \end{align*}
   Due to \eqref{estimatesFWga}, \eqref{estpabarW} and Lemma\ref{dampingtermsestimates}, it holds that
   \begin{align*}
   |\pa_1\tilde{F}_W|&\leq|\pa_1\mathcal{F}_W|+(|(\be_{\tau}J-1)\pa_1\bar{W}|+|
   \pa_1G_W|)|\pa_1\bar{W}|+(|\be_{\tau}J-1|+|G_W|)|\pa_1^2\bar{W}|\\
   &+|\pa_1h_W^{\mu}||\pa_{\mu}\bar{W}|+|h_W^{\mu}||\pa_1\pa_{\mu}\bar{W}|\\
   &\les \ta(1,\frac{1}{15})+\ep^{\frac{1}{2}}\ta(0,\frac{1}{2})+(M\ta(\frac{1}{2}-\frac{4}{2m-7},0)+M^{\frac{1}{4}}e^{-s}|y_1|+\ep^{\frac{1}{3}}|\check{y}|)\ta(0,\frac{1}{3})
   +\ta(1,\frac{1}{6})\\
   &\leq\ep^{\frac{1}{3}}\ta(0,\frac{1}{3}),
   \end{align*}
   and
   \begin{align*}
   |\check{\nabla}\tilde{F}_W|&\leq|\check{\nabla}\mathcal{F}_W|+
   (|\check{\nabla}[(\be_{\tau}J-1)\bar{W}]|+|\check{\nabla}G_W|)|\pa_1\bar{W}|\\
   &\leq(|(\be_{\tau}J-1)|+|G_W|)|\pa_1\check{\nabla}\bar{W}|+
   |\check{\nabla}h_W^{\mu}||\pa_{\mu}\bar{W}|+|h_W^{\mu}||\pa_{\mu}\check{\nabla}
   \bar{W}|\\
   &\les\ep^{\frac{1}{6}}\ta(0,\frac{1}{3}),
   \end{align*}
   which completes the proof for the estimates $\tilde{F}_W^{(\ga)}$ with $|\ga|=1$. The estimate for $|\ga|=4$ and $|y|\leq l$ is similar. For $\tilde{F}_W^{(\ga),0}$ with $|\ga|=3$, note that\footnote{This can be derived from \eqref{definitionbarWy} and \eqref{burgers2}. Moreover, one can show that $\nabla^{2n}\bar{W}^0=0$ for all $n\in \mathcal{N}$, see\cite{BSV3Disentropiceuler}.}
   \bee\label{barWtaylor}
   \bar{W}(y)=-y_1+y_1y_2^2+y_1y_3^2-3y_1^5-y_1y_2^4-
   y_1y_3^4-4y_1^3y_2^2-4y_1^3y_3^2-2y_1y_2^2y_3^2+O(|y|^6),
   \ee
   which implies that $\nabla^{2n}\bar{W}^0=0$ for $n=1,2$. Therefore, evaluating \eqref{tildeFWga} at $y=0$ yields
   \begin{align*}
   |\tilde{F}_W^{(\ga),0}|&\les|\pa^{\ga}\mathcal{F}_W^0|+|\nabla^3(\be_{\tau}J-1)|
   +|\nabla^3G_W^0|\\
   &+\left(|\nabla[(\be_{\tau}J-1)\bar{W}]^0|+|\nabla G_W^0|+|\nabla^3h_W^{\mu,0}|+|\nabla h_W^{\mu,0}|\right)\nabla^3\tilde{W}^0\\
   &\leq\ta(\frac{1}{2}-\frac{4}{2m-7},0)+\ep^{\frac{1}{4}}|a\be_1|\ta(\frac{1}{2},0)
   \leq\ta(\frac{1}{2}-\frac{4}{2m-7},0).
   \end{align*}
   \end{pf}
   \section{\textbf{Recover the bootstrap assumptions on modulation variables}}\label{section6}
   \hs The estimates for the modulation variables are one of the key aspects of this work and differ significantly from those in \cite{BSV3Dfulleuler}. Therefore, we will focus on the evolution of these modulation variables to understand how the damping term affects shock formation. To this end, we propose 10 constraints on $\pa^{\ga}W(0,s)$ for $|\ga|\leq2$, which effectively characterize the information about the shock\footnote{See also \eqref{burgersinitial} for the 1D case.}.
   \bee\label{constraints}
   W(0,s)=0,\hs \pa_1W(0,s)=-1,\hs \check{\nabla}W(0,s)=0,\hs \pa_1\nabla W(0,s)=0,\hs \check{\nabla}^2W(0,s)=0.
   \ee
   Initially, these constraints are satisfied, which provide us with 10 ODEs for the modulation variables. By solving these ODEs and appropriately selecting the modulation variables, one can ensure that \eqref{constraints} holds for later values of $s$. Substituting \eqref{constraints} into \eqref{pagaw} and evaluating at $y=0$ yields the following set of equations:
   \begin{align}
   G_W^0+F_W^0-2\be_1\be_{\tau}e^{-\frac{s}{2}}a(U\cdot N)^0&=e^{-\frac{s}{2}}\be_{\tau}\dot{\kappa},\label{W0}\\
   \pa_1F_W^0+\pa_1G_W^0-2\be_1\be_{\tau}e^{-\frac{s}{2}}a\pa_1(U\cdot N)^0&=\be_{\tau}-1,\label{pa1W0}\\
   \pa_{\mu}F_W^0+\pa_{\mu}G_W^0-2\be_1\be_{\tau}e^{-\frac{s}{2}}a\pa_{\mu}(U\cdot N)^0&=0,\label{pamuW0}\\
   \pa_{1i}F_W^0+\pa_{1i}G_W^0-2\be_1\be_{\tau}e^{-\frac{s}{2}}a\pa_{1i}(U\cdot N)^0&=G_W^0\pa_{1i1}W^0+h_W^{\mu,0}\pa_{1i\mu}W^0,\label{pa1iW0}\\
   \pa_{\nu\ga}F_W^0+\pa_{\nu\ga}G_W^0-2\be_1\be_{\tau}e^{-\frac{s}{2}}a\pa_{\nu\ga}(U\cdot N)^0&=G_W^0\pa_{1\nu\ga}W^0+h_W^{\mu,0}\pa_
   {\mu\nu\ga}W^0.\label{panugaW0}
   \end{align}
   Evaluating $\pa^{\ga}G_W$ and $\pa^{\ga}F_W$ at $y=0$ leads to
   \begin{align}
   \dfrac{1}{\be_{\tau}}G_W^0&=e^{\frac{s}{2}}\left(\kappa+\be_2Z^0-2\be_1R_{j1}\dot{\xi}_j\right),\label{GW0}\\
   \dfrac{1}{\be_{\tau}}\pa_1G_W^0&=\be_2e^{\frac{s}{2}}\pa_1Z^0,\label{pa1GW0}\\
   \dfrac{1}{\be_{\tau}}\pa_{\mu}G_W^0&=\be_2e^{\frac{s}{2}}\pa_{\mu}Z^0+
   2\be_1\dot{Q}_{1\mu}+2\be_1R_{j\ga}\dot{\xi}_j\fe_{\ga\mu},\label{pamuGW0}\\
   \dfrac{1}{\be_{\tau}}\pa_{11}G_W^0&=\be_2e^{\frac{s}{2}}\pa_{11}Z^0,\label{pa11GW0}\\
   \dfrac{1}{\be_{\tau}}\pa_{1\mu}G_W^0&=\be_2e^{\frac{s}{2}}\pa_{1\mu}Z^0-
   2\be_1e^{-\frac{3}{2}s}\dot{Q}_{\ga1}\fe_{\ga\mu},\label{pa1muGW0}\\
   \dfrac{1}{\be_{\tau}}\pa_{\mu\nu}G_W^0&=e^{-\frac{s}{2}}\left(
   -\dot{\fe}_{\mu\nu}+\be_2e^s\pa_{\mu\nu}Z^0-2\be_1(\dot{Q}_{\ga\mu}\fe_{\ga\nu}
   +\dot{Q}_{\ga\nu}\fe_{\ga\mu}+R_{j1}\dot{\xi}_jN_{1,\mu\nu}^0)+e^{-\frac{s}{2}}
   \frac{G_W^0}{\be_{\tau}}J^0_{,\mu\nu}\right),\label{pamunuGW0}
   \end{align}
   and
   \begin{align}
   \dfrac{1}{\be_{\tau}}\pa^{\be}F_W^0&=-\be_3(\kappa-Z^0)\pa^{\be}\pa_{\ga}A_{\ga}^0
   +\frac{1}{4}\be_4(\kappa-Z^0)^2e^s\pa^{\be}\pa_{1}K^0+\text{l.o.ts},\label{pamunuFW0}
   \end{align}
   for $|\be|\leq 2$. It follows from \eqref{pa1iW0} that\footnote{The estimates of $G_W^0$ and $h_W^0$ in Lemma\ref{dampingtermsestimates} are not enough to recover the bootstrap assumptions for the modulation variables.}
   \begin{align}
   G_W^0&=(\pa_{1i}F_W^0+\pa_{1i}G_W^0-\be_1\be_{\tau}e^{-\frac{s}{2}}a\pa_1\pa_iZ^0)\cdot (H_{i1}^{0})^{-1},\\
   h_W^{\mu,0}&=(\pa_{1i}F_W^0+\pa_{1i}G_W^0-\be_1\be_{\tau}e^{-\frac{s}{2}}a\pa_1\pa_iZ^0)\cdot (H_{i\mu}^{0})^{-1},
   \end{align}
   where $H^0$ is the following matrix:
   \bee
   H^0(s):=(\pa_1\nabla^2W)(0,s).
   \ee
   Note that \eqref{pa30tildew} and the property $\pa_1\nabla^2\bar{W}(0)\sim \text{Id}_{3\times3}$, which implies
   \bee
   |H^0(s)|\sim 1.
   \ee
   Therefore,
   \bee\label{GW0hW0}
   \bes
   |G_W^0|+|h_W^{\mu,0}|&\leq|\pa_{1i\mu}A_{\mu}^0|+e^s|\pa_{11i}K^0|+e^{\frac{s}{2}}|\pa_{1i}Z^0|
   +e^{-\frac{3}{2}s}|\dot{Q}||\fe|+e^{-\frac{s}{2}}|\be_1a\pa_{1i}Z^0|+\text{l.o.ts}\\
   &\leq\ta(\frac{3}{2}-\frac{5}{2m-5},0)+\ta(1-\frac{4}{2m-7},\frac{1}{15})+\ep^{\frac{1}{2}}\ta(1,0)+\ta(\frac{3}{2},0)\\
   &\leq\ta(1-\frac{4}{2m-7},0).
   \end{split}
   \ee
   \subsection*{\textbf{$\dot{\tau}$ estimates}}
   It follows from \eqref{pa1W0} that
   \begin{equation}\label{dottau1}
   \bes
   \dot{\tau}&=\dfrac{1}{\be_{\tau}}\left(\pa_1F_W^0+\pa_1G_W^0\right)-2\be_1ae^{-\frac{s}{2}}\pa_1(U\cdot N)^0\\
   &=\dfrac{1}{\be_{\tau}}\left(\pa_1F_W^0+\pa_1G_W^0\right)+a\be_1e^{-s}-\be_1ae^{-\frac{s}{2}}\pa_1Z^0\\
   &=\be_1ae^{-s}+\be_2e^{\frac{s}{2}}\pa_1Z^0-\be_3(\kappa-Z^0)\pa_{1\mu}A_{\mu}
   +\frac{1}{4}\be_4(\kappa-Z^0)^2e^s\pa_{11}K^0-\be_1ae^{-\frac{s}{2}}\pa_1Z^0+\text{l.o.ts}\\
   &=\be_1ae^{-s}+O(e^{-s})=\be_1a(\tau(t)-t)+O(e^{-s}).
   \end{split}
   \end{equation}
   If the damping term vanishes, then $\dot{\tau}=O(e^{-s})$, which resulting in $O(\ep^2)$ bounds in the estimates of $|\tau(t)|$. Here we investigate the impact of the damping term and thus consider the following ODE:
   \bee
   \dot{\tilde{\tau}}(t)=\be_1a(\tilde{\tau(t)}-t),\hs \tilde{\tau}(-\ep)=0,
   \ee
   which gives us
   \bee
   \tilde{\tau}(t)=t+\ep e^{\be_1a(t+\ep)}+\dfrac{1}{\be_1a}(1-e^{\be_1a(t+\ep)}),
   \ee
   and then
   \bee
   \tau(t)=t+\ep e^{\be_1a(t+\ep)}+\dfrac{1}{\be_1a}(1-e^{\be_1a(t+\ep)})+O(\ep^2).
   \ee
   Therefore,
   \begin{itemize}
   \item If $\dfrac{1}{\ep}\leq\dfrac{\be_1a}{2}$, then $\tau(t)>t$ and $\dot{\tau}-1>0$ for all $t$. Consequently, it follows from \eqref{jocabian} that the solution is smooth both in the self-similar coordinates and in the rectangular coordinates. Thus, one obtains a global solution to the Euler system\eqref{euler};
   \item if $\dfrac{1}{\ep}\geq2\be_1a$\footnote{If $a=0$, then the damping term vanishes and this case is automatically satisfied, which implies for standard compressible Euler equations with the initial data given by Lemma\ref{initialdata}, shock formation is inevitable for sufficiently small $\ep$.}, then $\tau(T_{\ast})=T_{\ast}$, which implies a shock forms at $t=T_{\ast}=-\frac{1}{\be_1a}\ln(1-\ep\be_1a)-\ep+O(\ep^2)$. Furthermore, the shock time is shifted compared with the work in \cite{BSV3Dfulleuler}.
   \end{itemize}
   \subsection*{\textbf{$\dot{\kappa}$ estimates}}
   It follows from \eqref{W0} that
   \begin{align}
   \dot{\kappa}(t)&=\dfrac{1}{\be_{\tau}}e^{\frac{s}{2}}(F_W^0
   +G_W^0)-2\be_1a(U\cdot N)^0\\
   &=-\be_1a\kappa+\dfrac{1}{\be_{\tau}}e^{\frac{s}{2}}G_W^0-\be_3(\kappa-Z^0)e^{\frac{s}{2}}\pa_{\mu}A_{\mu}^0
   +\frac{1}{4}\be_4(\kappa-Z^0)^2e^{\frac{3}{2}s}\pa_{11}K^0-\be_1aZ^0+\text{l.o.ts}\\
   &=-\be_1a\kappa+O(\ep^{\frac{1}{5}}),
   \end{align}
   which implies that
   \bee
   \kappa=\kappa_0e^{-\be_1a(t+\ep)}+O(\ep^{\frac{1}{5}})\dfrac{1}{\be_1a}
   \left(1-e^{-\be_1a(t+\ep)}\right).
   \ee
   Since as $t\to T_{\ast}$, $e^{\be_1a(t+\ep)}\to \frac{1}{1-\ep\be_1a}\cdot e^{O(\ep^2)}$. Then, it follows that $\dfrac{1}{\be_1a}
   \left(1-e^{-\be_1a(t+\ep)}\right)=O(\ep)$. Therefore,
   \bee\label{kappaestimates}
   \kappa=\kappa_0e^{-\be_1a(t+\ep)}+O(\ep^{\frac{6}{5}}).
   \ee
   Therefore, one can conclude that if $a$ is positive, the wave amplitude will decay as $t$ increases, approaching $\kappa_0(1-\ep\be_1a)+O(\ep^{\frac{6}{5}})$ at the shock point while if $a$ is negative, the wave amplitude will grow as $t$ increases.
   \subsection*{\textbf{$\dot{\xi}$ estimates}}
   It follows from the definition of $G_W^0, h_W^{\mu,0}$ \eqref{defGWGZGU}, the bootstrap assumptions\eqref{bspaZ},\eqref{bspaAmu} and\eqref{GW0hW0} that\footnote{Different form the equations\eqref{W0}-\eqref{pamuW0}, we don't postulate the precise value for $\pa_{ijk}W^0$. So, one can not obtain the estimates for $\xi$ directly from the equation\eqref{pa1iW0} as well as the estimates for $\fe$. The use for \eqref{pa1iW0}-\eqref{panugaW0} is to derive the accurate estimates for $G_W^0$ and $h_W^0$ (see \eqref{GW0hW0}).}
   \begin{equation}
   \bes
   \dot{\xi}_i&=(\dot{\xi}_jR_{jk})R_{ik}\\
   &=\dfrac{R_{i1}}{2\be_1}\left(\kappa+\be_2Z^0-\frac{1}{\be_{\tau}}G_W^0\right)+R_{j\mu}\left(A_{\mu}^0-\frac{1}{2\be_1\be_{\tau}}e^{\frac{s}{2}}h_W^{\mu,0}\right)\\
   &=O(1)\kappa+O(\ep^{\frac{1}{2}})=O(1)\kappa_0e^{-\be_1a(t+\ep)}+O(\ep^{\frac{1}{2}}).
   \end{split}
   \end{equation}
   Note that the location of shock is defined as $\xi_{\ast}=\xi(T_{\ast})$. Then, it follows that
   \begin{align}
   \xi_{\ast}&=\int_{-\ep}^{T_{\ast}}\dot{\xi}=O(1)\int_{-\ep}^{T_{\ast}}\kappa_0e^{-\be_1a(t+\ep)} dt+O(\ep^{\frac{3}{2}})\\
   &=\kappa_0\frac{1}{\be_1a}(1-(1-\ep\be_1a))+O(\ep^{\frac{3}{2}})=O(\ep)<M^{\frac{1}{4}}\ep,
   \end{align}
   which is independent of $a$. Therefore, the damping effect doesn't shift the shock location.
   \subsection*{\textbf{$\dot{n}$,$\dot{\fe}$ estimates}}
   It follows from \eqref{pamuGW0} and \eqref{pamuW0} that
   \begin{equation}
   \bes
   \dot{Q}_{1\mu}&=-\frac{1}{2\be_1\be_{\tau}}\pa_{\mu}F_w^0+\frac{1}{2}e^{-\frac{s}{2}}a\pa_{\mu}Z^0\\
   &-\frac{\be_2}{2\be_1}e^{\frac{s}{2}}\pa_{\mu}Z^0-A_{\ga}^0\fe_{\ga\mu}+\dfrac{1}{2\be_1\be_{\tau}}e^{\frac{s}{2}}h_W^{\ga,0}\fe_{\ga\mu}.
   \end{split}
   \end{equation}
   It follows from the definition of $\dot{Q}$ \eqref{definitionQ} that
   \bee\label{dotnsimQ}
   \left(\begin{array}{cc}
   1+\frac{n_2^2}{n_1(1+n_1)} & \frac{n_2n_3}{n_1(1+n_1)}\\
   \frac{n_2n_3}{n_1(1+n_1)} & 1+\frac{n_3^2}{n_1(1+n_1)}
   \end{array}\right)\cdot
   \left(\begin{array}{c}
   \dot{n}_2\\
   \dot{n}_3
   \end{array}\right)=\left(Id+\frac{\check{n}\otimes\check{n}}{n_1(1+n_1)}\right)\cdot \dot{\check{n}}=
   \left(\begin{array}{c}
   \dot{Q}_{12}\\
   \dot{Q}_{13},
   \end{array}\right).
   \ee
   These together with Lemma\ref{forcingtermsestimates} lead to
   \begin{align*}
   \dot{n}_{\mu}&\sim \dot{Q}_{1\mu}\sim-\pa_{\mu\nu}A_{\nu}^0+e^s\pa_{1\mu}K^0-ae^{-\frac{s}{2}}\pa_{\mu}Z^0-A_{\ga}^0\fe_{\ga\mu}+e^{\frac{s}{2}}h_W^{\ga,0}\fe_{\ga\mu},
   \end{align*}
   which implies\footnote{One also obtains $|\dot{Q}(t)|\leq M^2\ep^{\frac{1}{2}}$.}
   \begin{align}
     |\dot{n}_{\mu}|&\leq Me^{-s}+\ep^{\frac{1}{8}}e^{-\frac{5}{8}s}+M\ep^{\frac{1}{2}}+M^2\ep\leq M\ep^{\frac{1}{2}}<M^2\ep^{\frac{1}{2}}+O(\ep)\footnotemark,\label{estimatedotn}\\
    |n_{\mu}|&\leq\int_{-\ep}^{T_{\ast}}|\dot{n}_{\mu}|\leq M\ep^{\frac{3}{2}}<M^2\ep^{\frac{3}{2}},
   \end{align}
   \footnotetext{Indeed, the term $O(\ep)$ in \eqref{estimatedotn} should be $a\ep^{\frac{3}{2}}$, which may depend on $a$. However, if a shock forms, then $|a\ep^{\frac{3}{2}}|\les\ep^{\frac{1}{2}}$, where the implicit constant is independent of $a$.}
   due to the bootstrap assumptions\eqref{bsmodulation3}, \eqref{bspaZ}-\eqref{bspaK}, where the estimates are independent of $a$. Therefore, the damping term does not affect the blow up direction.\\
   \hs It follows from \eqref{pamunuGW0} that
   \bee\label{equationdotfe}
   \bes
   \dot{\fe}_{\mu\nu}&=-\frac{1}{\be_{\tau}}e^{\frac{s}{2}}\pa_{\mu\nu}G_W^0+\be_2e^s\pa_{\mu\nu}Z^0-2\be_1(\dot{Q}_{\zeta\mu}\fe_{\zeta\nu}+R_{j1}\dot{\xi}_jN_{1,\mu\nu}^0)+e^{-\frac{s}{2}}\frac{G_W^0}{\be_{\tau}}J^0_{,\mu\nu},
   \end{split}
   \ee
   which implies
   \bee
   \bes
    |\dot{\fe}|&\leq M+M+M^2\ep+\ep^{\frac{1}{2}}e^{-s}\leq 2M<M^2,\\
    |\fe|&\leq 2M\ep<M^2\ep,
    \end{split}
   \ee
   due to Lemma\ref{dampingtermsestimates}, \ref{prelimi}, the bootstrap assumptions\eqref{bsmodulation1},\eqref{bsmodulation2} and \eqref{bspaZ}.
   \begin{remark}
   The equations \eqref{W0}-\eqref{pamunuFW0} give system of $10$ ODEs for the modulation variables $(\dot{\tau},\dot{\kappa},\dot{\xi},\dot{\check{n}},\dot{\fe})$, where the coefficients are at least $C^2$ due to the bootstrap assumptions. Therefore, unique $C^1$ existence for the modulation variables is guaranteed in a small time and then one can determine the evolution of these variables for later $t$.
   \end{remark}
   \subsection{\textbf{Damping effect to the modulation variables}}\label{section6.1}
   In conclusion,
   \begin{itemize}
   \item if $\dfrac{1}{\ep}\leq\dfrac{\be_1a}{2}$, then for fixed $\ep$(initial small data)\footnote{Note that $\frac{1}{\ep}$ measures the maximum of $\pa_1\til{w}(-\ep,x)$.}, the damping effect is strong enough to prevent the shock formation while for fixed $a$, initial small data leads to a global solution. Therefore, Therefore, a global solution can be obtained in both the self-similar coordinates and the physical variables. In this case (see also \eqref{blowupquantity}),
       \bee
       N\cdot\nabla_{\tilde{x}}\tilde{w}(\xi(t),t)=\dfrac{1}{\tau(t)-t}>0, \text{for all}\ t\geq-\ep.
       \ee
   \item If $\dfrac{1}{\ep}\geq2\be_1a$, the damping effect is so weak that a shock forms in finite time for fixed $\ep$ while for fixed $a$, the initial large data leads to the shock formation in finite time. In this case, one can explicitly compute the shock time $T_{\ast}$ and location $\xi_{\ast}$. While the shock location and direction remain unchanged, the damping effect shifts the shock time and alters the wave amplitude compared with the undamped case, with the variation in shift or change depending on the sign of $a$. In particular, if $a$ is positive, the damping term will delay the shock formation and reduce the wave amplitude. Conversely, if $a$ is negative, the (anti)-damping term will lead to an immediate shock and amplify the wave amplitude.
   \end{itemize}
   \section{\textbf{Bounds on Lagrangian trajectories and vorticity variation}}\label{section7}
   Define Lagrangian flows as follows:
   \begin{align}
   \dfrac{d}{ds}\Phi_W(y,s)&=V_W(\Phi_W(y,s),s),\hs \Phi_W(y,s_0)=y,\label{definitionPhiW}\\
   \dfrac{d}{ds}\Phi_Z(y,s)&=V_Z(\Phi_Z(y,s),s),\hs \Phi_Z(y,s_0)=y,\label{definitionPhiZ}\\
   \dfrac{d}{ds}\Phi_U(y,s)&=V_U(\Phi_U(y,s),s),\hs \Phi_U(y,s_0)=y.\label{definitionPhiU}
   \end{align}
   Denote $\Phi^{y_0}(s)$ to be the trajectory of either $\Phi_W,\Phi_Z$ or $\phi_U$ emanating from $y_0$ at time $s$, i.e., $\Phi^{y_0}(s)=\Phi(y_0,s)$, $\phi^{y_0}(s_0)=y_0.$ The following lemma recovers the bootstrap assumption \eqref{bssupport}.
   \begin{lem}
   Let $\Phi$ be either $\Phi_W,\Phi_Z$ or $\phi_U$. Then
   \begin{align}
   |\Phi_1(s)|&\leq\dfrac{3}{2}\ep^{\frac{1}{2}}e^{\frac{3}{2}s},\label{phi1}\\
   |\check{\Phi}(s)|&\leq\dfrac{3}{2}\ep^{\frac{1}{6}}e^{\frac{s}{2}},
   \end{align}
   for all $s\geq-\log\ep.$
   \end{lem}
   \begin{pf}
   The proof for \eqref{phi1} with the case of $\Phi=\Phi_W$ will be given and the remaining proofs are the same. It follows from \eqref{definitionPhiW} and the definition of $V_W$ \eqref{defGWGZGU} that
   \bee
   \dfrac{\pa}{\pa s}\Phi_1=G_W\circ\Phi+\dfrac{3}{2}\Phi_1+\be_{\tau}JW\circ\Phi,
   \ee
   which implies
   \bee
   \dfrac{d}{ds}\left(e^{-\frac{3}{2}s}\Phi_1(s)\right)=e^{-\frac{3}{2}s}\left(G_W+\be_{\tau}JW\right)\circ\Phi.
   \ee
   Integrating from $-\log\ep$ to $s$ yields
   \bee
   \bes
   e^{-\frac{3}{2}s}\Phi_1(s)-\ep^{\frac{3}{2}}y_0&=\int_{-\log\ep}^se^{-\frac{3}{2}s'}\left(G_W+\be_{\tau}JW\right)\circ\Phi\ ds'\\
   &\leq\int_{-\log\ep}^s2\ep^{\frac{1}{6}}Me^{-s'}\ ds'\leq 2M\ep^{\frac{7}{6}}
   \end{split}
   \ee
   which implies
   \bee
   |\Phi_1(s)|\leq\dfrac{3}{2}\ep^{\frac{1}{2}}e^{\frac{3}{2}s}.
   \ee
   \end{pf}
   \subsection{\textbf{Lower bounds for the Lagrangian trajectories}}
   \begin{prop}\label{lowerboundslagrangian}
   \begin{itemize}
   \item Let $y_0\in \mbr^3$ and $|y_0|\geq l$. Let also $s_0\geq-\log\ep$. Then, the trajectory $\Phi_W^{y_0}$ moves away from the origin at an exponential rate with
   \bee\label{phiw}
   |\check{\Phi}_W^{y_0}(s)|\geq |y_0|e^{\frac{s-s_0}{5}}.
   \ee
   \item Let $\Phi(s)$ be either $\Phi_Z^{y_0}(s)$ or $\Phi_U^{y_0}(s)$. Then, for any $y_0$, there exists an $s_{\ast}\geq-\log\ep$ such that
   \bee\label{phizu}
   |\Phi_1(s)|\geq\min\left\{|e^{\frac{s}{2}}-e^{\frac{s_{\ast}}{2}}|,e^{\frac{s}{2}}\right\}
   \ee
   by choosing $\kappa_0$ suitably large.
   \end{itemize}
   \end{prop}
   \begin{remark}
   Note that the Lagrangian trajectories of $(W,Z,A_v,K)$ rapidly approach spatial infinity with an exponential growth, which allows one to utilize the spatial decay of different damping and forcing terms to achieve integrable time decay. Furthermore, the estimates of $\Phi_W(s)$ depends on the distance of the initial position, while $(Z,A_v,K)$ do not exhibit such dependence. This distinction requires us to estimate $W$ separately in various regions. Different from the work in\cite{BSV3Dfulleuler}, the first component of the Lagrangian trajectory of $W$ may decrease initially but then exhibit exponential growth after a short period due to the damping effect, while the other two components do not.
   \end{remark}
   \begin{pf}
   \begin{itemize}
   \item[(1)] By Gronwall inequality, in order to prove \eqref{phiw}, it suffices to show
   \bee
   \dfrac{1}{2}\dfrac{d}{ds}|\check{\Phi}_W^{y_0}(s)|^2\geq\dfrac{1}{5}|\check{\Phi}_W^{y_0}|^2.
   \ee
    This is equivalent to show
   \bee
   \check{y}\cdot \check{V}_W(y)\geq\dfrac{2}{5}|\check{y}|^2, \text{for}\ |y|\geq l,
   \ee
   due to the definition of $\Phi_W$. It follows from Lemma\eqref{dampingtermsestimates} that
   \begin{align}
   \check{y}\cdot\check{V}_W(y)&=\dfrac{1}{2}|\check{y}|^2+y_{\mu}h_W^{\mu}\geq\dfrac{1}{2}|\check{y}|^2-e^{-\frac{s}{2}}|\check{y}|\geq\dfrac{2}{5}|\check{y}|^2,
   \end{align}
   where in the last inequality, $|y|\geq l$ is used.
   \item[(2)] Suppose $|\Phi_1(s)|\leq e^{\frac{s}{2}}$. Then, it suffices to prove $|\Phi_1(s)|\geq|e^{\frac{s}{2}}-e^{\frac{s_{\ast}}{2}}|$, i.e.,
   $|\Phi_1(s)|\geq e^{\frac{s}{2}}-e^{\frac{s_{\ast}}{2}}$ or $|\Phi_1(s)|\leq e^{\frac{s_{\ast}}{2}}-e^{\frac{s}{2}}$. It can be shown by Gronwall inequality that the former case is contradict to the assumption $|\Phi_1(s)|\leq e^{\frac{s}{2}}$\footnote{If so, then $\frac{d}{ds}\Phi_1\geq e^{\frac{s}{2}}$ provided $\Phi_1(s_{\ast})\geq0$, which can be integrated to show $\Phi_1(s)\geq e^{\frac{s}{2}}+(y_0)_1-e^{\frac{s_0}{2}}>e^{\frac{s}{2}}$ by taking $|y_0|$ suitably.}. To show $|\Phi_1^{y_0}(s)|\leq e^{\frac{s_{\ast}}{2}}-e^{\frac{s}{2}}$, it suffices to show
   \bee\label{inequalityPhi1}
   \dfrac{d}{ds}\Phi_1^{y_0}(s)\leq-\dfrac{1}{2}e^{\frac{s}{2}}
   \ee
    provided $\Phi_1(s_{\ast})\leq0$. It follows from Lemma\eqref{dampingtermsestimates} that
   \begin{align*}
   \dfrac{d}{ds}\Phi_{Z_1}(s)&=V_{Z_1}\circ\Phi_Z(s)=(\be_2\be_{\tau}JW+G_Z)\circ\Phi_Z(s)+\dfrac{3}{2}\Phi_{Z_1}(s)\\
   &\leq-(1-\be_2)e^{\frac{s}{2}}\kappa_0+\be_2\be_{\tau}JW\circ\Phi_Z(s)+\ep^{\frac{1}{2}}e^{\frac{s}{2}}+\frac{3}{2}e^{\frac{s}{2}}\\
   &\leq-(1-\be_2)e^{\frac{s}{2}}\kappa_0+\be_2\be_{\tau}J\left(|\Phi_{Z_1}(s)|+\ep^{\frac{1}{11}}|\Phi_{Z_{\mu}}(s)|\right)+\left(\dfrac{3}{2}+\ep^{\frac{1}{2}}\right)e^{\frac{s}{2}}\\
   &\leq-(1-\be_2)e^{\frac{s}{2}}\kappa_0+C(3+\ep^{\frac{1}{2}})e^{\frac{s}{2}}<-\dfrac{1}{2}e^{\frac{s}{2}},
   \end{align*}
   provided $\kappa_0\geq\dfrac{C}{1-\be_2}$, and similarly for $\Phi_{U_1}(s)$. Integrating \eqref{inequalityPhi1} from $s_{0}$ to $s$ yields
   \bee
   \Phi_1^{y_0}(s)-\Phi_1^{y_0}(s_0)\leq e^{\frac{s_0}{2}}-e^{\frac{s}{2}}.
   \ee
   To finish the proof, it suffices to show the existence of $s_{\ast}$ and there are following two cases.
   \begin{itemize}
   \item If $\Phi_1(s_0)\leq 0$, then take $s_{\ast}=s_0$.
   \item If $\Phi_1(s_0)>0$, since $\dfrac{d}{ds}\Phi_1^{y_0}(s)\leq-\dfrac{1}{2}e^{\frac{s}{2}}\leq-M$, then there exists an $s_{\ast}\leq-3\log\ep$ such that $\Phi_1^{y_0}(s_{\ast})=0$ where $s_{\ast}$ is independent of $s$. Then, integrating \eqref{inequalityPhi1} from $s_{\ast}$ to $s$ leads to
       \bee
       \Phi_1^{y_0}(s)\leq e^{\frac{s_{\ast}}{2}}-e^{\frac{s}{2}}.
       \ee
   \end{itemize}
   \end{itemize}
   \end{pf}
   The following lemma is the key lemma of the function $\ta(\a,\be)$.
   \begin{lem}\label{keylemma}
   Let $\ta(\a,\be)(y,s)$ be defined in section4.1. 
   Then, the following properties hold.
   \begin{itemize}
   \item Let $|y_0|\geq l$. Suppose one of the following the conditions holds:
         \begin{align}
         \a\geq0 \ \text{or}\ \be\leq 0,\ \text{and}\ 3\be+\a>0;\label{case1w}\\
         \a\leq 0\ \text{and}\ \be\geq 0,\ 2\be+5\a>0\label{case2w}.
         \end{align}
         Then,
         \bee
         \int_{s_0}^{+\infty}\ta(\a,\be)\circ\Phi_W^{y_0}\ ds\les 1,
         \ee
         where $\Phi_W^{y_0}(s_0)=y_0$ and the implicit constant only depends on $\a,\be$ and $|y_0|$ but not on $M$.
   \item Let $\Phi$ be either $\Phi_Z$ or $\Phi_U$. Suppose one of the following the conditions holds:
         \begin{align}
         \a\geq 0 \ \text{or}\ \be\leq 0,\ \text{and}\ 3\be+\a>0;\label{case1z}\\
         -\frac{1}{2}\leq\a\leq 0\ \text{and}\ \be\geq 0,\ \be+\a>0.\label{case2z}
         \end{align}
         Then,
         \bee
         \int_{-\log\ep}^{+\infty}\ta(\a,\be)\circ\Phi^{y_0}\ ds\les 1,
         \ee
         where $\Phi^{y_0}(-\log\ep)=y_0$ and the implicit constant only depends on $\a,\be$.
   \end{itemize}
   \end{lem}
   \begin{pf}
   Suppose \eqref{case1w} hold. If $\a\geq0$, then it follows from Proposition\ref{lowerboundslagrangian} that
   \bee\label{case1w1}
   \ta(\a,\be)\circ\Phi_W^{y_0}\leq\ep^{\frac{\a}{3}}\eta^{-\frac{1}{3}(\a+3\be)}\circ\Phi_W^{y_0}
   \leq\ep^{\frac{\a}{3}}\left(1+|y_0|^2e^{-\frac{2(s-s_0)}{5}}\right)^{-\frac{1}{3}(\a+3\be)},
   \ee
   while if $\be\leq0$, it holds that
   \bee\label{case1w2}
   \ta(\a,\be)\circ\Phi_W^{y_0}\leq \ep^{-\frac{3\be}{2}}e^{-(\a+3\be)s},
   \ee
   due to remark\ref{taprop}. Direct computation yields, in either case, $\int_{s_0}^{+\infty}\ta(\a,\be)\circ\Phi_W^{y_0}\les 1$. Suppose \eqref{case2w} hold. Then, it follows from Proposition\ref{lowerboundslagrangian} that
   \begin{align*}
   \int_{s_0}^{+\infty}\ta(\a,\be)\circ\Phi_W^{y_0}\ ds&\leq\int_{s_0}^{+\infty}e^{-\a s}\left(1+|y_0|^2e^{-\frac{2(s-s_0)}{5}}\right)^{-\be}\ ds\\
   &\leq|y_0|^{-2\be}\int_{s_0}^{+\infty}e^{-\frac{2\be+5\a}{5}}e^{\frac{2\be s_0}{5}}\ ds\\
   &\leq|y_0|^{-2\be}e^{\a s_0}\les 1.
   \end{align*}
   Suppose \eqref{case1z} hold. Then $\int_{-\log\ep}^{+\infty}\ta(\a,\be)\circ\Phi^{y_0}\ ds\les 1$ with $\Phi_Z$ or $\Phi_U$ can be derived in a similar way as in $\Phi_W^{y_0}$. If \eqref{case2z} is satisfied, then it follows from \eqref{phizu} that $|\Phi_1^{y_0}|\geq\min\{e^{\frac{s}{2}},|e^{\frac{s}{2}}-e^{\frac{s_{\ast}}{2}}|\}$. Then, there are the following two cases.
   \begin{itemize}
   \item[Case1]
   If $|\Phi_1^{y_0}|\geq e^{\frac{s}{2}}$, then
   \begin{align*}
   \int_{-\log\ep}^{+\infty}\ta(\a,\be)\circ\Phi^{y_0}\ ds&\leq\int_{-\log\ep}^{+\infty}e^{-\a s}(1+e^{\frac{s}{2}})^{-2\be}\ ds\\
   &\leq \int_{-\log\ep}^{+\infty}e^{-(\be+\a)s}\ ds\leq C(\a,\be).
   \end{align*}
   \item[Case2]
   If $|\Phi_1^{y_0}|\geq|e^{\frac{s}{2}}-e^{\frac{s_{\ast}}{2}}|$, then
   \begin{align*}
   \int_{-\log\ep}^{+\infty}\ta(\a,\be)\circ\Phi^{y_0}\ ds&\leq \int_{-\log\ep}^{+\infty}e^{-\a s}(1+|e^{\frac{s}{2}}-e^{\frac{s_{\ast}}{2}}|)^{-2\be}\ ds\\
   &=\int_{\ep^{-\frac{1}{2}}}^{+\infty}2 p^{-2\a-1}(1+|p-e^{\frac{s_{\ast}}{2}}|)^{-2\be}\ dp\\
   &\les\int_{\ep^{-\frac{1}{2}}}^{+\infty}p^{-2\a-2\be-1}+(1+|p-e^{\frac{s_{\ast}}{2}}|)^{-2\a-2\be-1}\ dp\les 1.
   \end{align*}
   \end{itemize}
   \end{pf}
   \begin{remark}
   Note that for the case\eqref{case1w},\eqref{case2w} and \eqref{case1z}, the above constant $1$ can be refined\footnote{For example, it follows from \eqref{case1w1} that $\int_{s_0}^{+\infty}\ta(\a,\be)\circ\Phi_W^{y_0}\les \ep^{\frac{\a}{3}}$ while \eqref{case1w2} yields $\int_{s_0}^{+\infty}\ta(\a,\be)\circ\Phi_W^{y_0}\les \ep^{\frac{3\be}{2}}$.}, which will be used later.
   \end{remark}
   \subsection{\textbf{Vorticity bounds and variation}}
   We now establish the bounds for the specific vorticity. Recall the equation for the specific vorticity $\mathring{\zeta}$:
   \bee\label{vorticity}
   \bes
   &\pa_t\mathring{\zeta}-2\be_1\mathring{\zeta}\dot{Q}+2\be_1\left[JN\cdot(\mathring{u}+\mathring{v})-\dfrac{\dot{f}}{2\be_1}\right]\pa_1\mathring{\zeta}
   +2\be_3(\mathring{u}_{\mu}+\mathring{v}_{\mu})\pa_{\mu}\mathring{\zeta}
   -2\be_1[JN\cdot\mathring{\zeta}\pa_1\mathring{u}+\mathring{\zeta}_{\mu}\pa_{\mu}\mathring{u}]\\
   &=2\be_4\dfrac{\mathring{\si}}{\mathring{\p}}\nabla_{\tilde{x}}\mathring{\si}\times\nabla_{\tilde{x}}\mathring{k}-2\be_1a\mathring{\zeta}.
   \end{split}
   \ee
   Decompose \eqref{vorticity} along the tangential direction $(T^2,T^3)$ as follows.
   \begin{align}
   \pa_t(\mr{\zeta}\cdot T^2)+(v\cdot\nabla)(\mr{\zeta}\cdot T^2)&=F_{21}(\mr{\zeta}\cdot N)+F_{2\mu}(\mr{\zeta}\cdot T^{\mu})+G_2-2\be_1a(\mr{\zeta}\cdot T^2),\label{vorticityT2}\\
   \pa_t(\mr{\zeta}\cdot T^3)+(v\cdot\nabla)(\mr{\zeta}\cdot T^3)&=F_{31}(\mr{\zeta}\cdot N)+F_{3\mu}(\mr{\zeta}\cdot T^{\mu})+G_3-2\be_1a(\mr{\zeta}\cdot T^3),\label{vorticityT3}
   \end{align}
   where
   \begin{align*}
   F_{21}&=N\pa_tT^2+v_jN\pa_jT^2+2\be_1\dot{Q}_{ij}T_j^2N_i+2\be_1[J\pa_1a_2+N_{\mu}(\pa_{\mu}a_2-\mr{u}\pa_{\mu}T^2)],\\
   F_{22}&=2\be_1T^2_{\mu}(\pa_{\mu}a_2-\mr{u}\pa_{\mu}T^2),\\
   F_{23}&=T^3\pa_tT^2+v_jT^3\pa_jT^2+2\be_1\dot{Q}_{ij}T^2_jT^3_i+2\be_1T^3_{\mu}(\pa_{\mu}a_2-\mr{u}\pa_{\mu}T^2),\\
   G_2&=2\be_4\dfrac{\mr{\sigma}}{\mr{\p}}\left(\nabla_{\tilde{x}}\mr{\sigma}\times\nabla_{\tilde{x}}\mr{k}\right)\cdot T^2=2\be_4\dfrac{\mr{\sigma}}{\mr{\p}}
   \left(\pa_{T^3}\mr{\sigma}\pa_N\mr{k}-\pa_N\mr{\sigma}\pa_{T^3}\mr{k}\right),
   \end{align*}
   with the bounds $|F_{21}|\leq M, |F_{2\mu}|\leq M\ep^{\frac{1}{2}}\leq \ep^{\frac{1}{4}}$. Similarly, one has
    \begin{align*}
    F_{31}&=N\pa_tT^3+v_jN\pa_jT^3+2\be_1\dot{Q}_{ij}T_j^3N_i+2\be_1[J\pa_1a_3+N_{\mu}(\pa_{\mu}a_3-\mr{u}\pa_{\mu}T^3)],\\
   F_{32}&=T^2\pa_tT^3+v_jT^2\pa_jT^3+2\be_1\dot{Q}_{ij}T^3_jT^2_i+2\be_1T^2_{\mu}(\pa_{\mu}a_3-\mr{u}\pa_{\mu}T^3),\\
   F_{33}&=2\be_1T^3_{\mu}(\pa_{\mu}a_3-\mr{u}\pa_{\mu}T^3),\\
   G_3&=2\be_4\dfrac{\mr{\sigma}}{\mr{\p}}\left(\nabla_{\tilde{x}}\mr{\sigma}\times\nabla_{\tilde{x}}\mr{k}\right)\cdot T^3=2\be_4\dfrac{\mr{\sigma}}{\mr{\p}}
   \left(\pa_{N}\mr{\sigma}\pa_{T^2}\mr{k}-\pa_{T^2}\mr{\sigma}\pa_{N}\mr{k}\right),
    \end{align*}
    with the bounds $|F_{31}|\leq M,|F_{3\mu}|\leq \ep^{\frac{1}{4}}$. It follows from the definition of $\mr{\zeta}$ that along normal direction $N$,
   \begin{equation}
   \bes
   \mr{\zeta}\cdot N&=\dfrac{curl\mr{u}\cdot N}{\tilde{\p}}=\dfrac{1}{\tilde{\p}}(\pa_{T^2}\mr{u}\cdot T^3-\pa_{T^3}\mr{u}\cdot T^2)\\
   &=\dfrac{1}{\tilde{\p}}(T^2_i\pa_{\tilde{x_i}}\mr{u}\cdot T^3-T^3_i\pa_{\tilde{x_i}}\mr{u}\cdot T^2)\\
   &=\dfrac{1}{\tilde{\p}}\left[T^2\cdot JN\pa_1\mr{u}\cdot T^3+T^2_{\mu}\pa_{\mu}\mr{u}\cdot T^3-T^2\cdot JN\pa_1\mr{u}\cdot T^2-T^3_{\mu}\pa_{\mu}\mr{u}\cdot T^2\right]\\
   &=\dfrac{1}{\tilde{\p}}\left[T^2_{\mu}\pa_{\mu}a_3-T^2_{\mu}\mr{u}\cdot\pa_{\mu}T^3-T_{\mu}^3\pa_{\mu}a_2+T^3_{\mu}\mr{u}\cdot\pa_{\mu}T^2\right],
   \end{split}
   \end{equation}
   which implies
   \bee\label{vorticitynormal}
   |\mr{\zeta}\cdot N|\les M^{\frac{1}{2}}\ep^{\frac{1}{2}}\leq\ep^{\frac{1}{4}}.
   \ee
   Define the Lagrangian flow associated to $v$ as follows:
   \bee
   \dfrac{d}{dt}\fai(x,t)=v\circ\fai,\hs \fai(x,-\ep)=x.
   \ee
   Along $\fai$, denote
   \bee
   Q_2=(\mr{\zeta}\cdot T^2)\circ\fai,\hs  Q_3=(\mr{\zeta}\cdot T^3)\circ\fai,\hs Q_1=(\mr{\zeta}\cdot N)\circ\fai.
   \ee
   Then, \eqref{vorticityT2} and \eqref{vorticityT3} can be written as
   \begin{align}
   \dfrac{d}{dt}Q_2&=F_{21}\circ\fai Q_1+F_{2\mu}\circ\fai Q_{\mu}+G_2\circ\fai-2\be_1a Q_2,\\
   \dfrac{d}{dt}Q_3&=F_{31}\circ\fai Q_1+F_{3\mu}\circ\fai Q_{\mu}+G_3\circ\fai-2\be_1a Q_3.
   \end{align}
   That is,
   \bee\label{vorticity2}
   \dfrac{1}{2}\dfrac{d}{dt}(Q_2^2+Q_3^2)=F_{\mu1}\circ\fai Q_1Q_{\mu}+F_{\mu\nu}\circ\fai Q_{\mu}Q_{\nu}+G_{\mu}\circ\fai Q_{\mu}-2\be_1a(Q_2^2+Q_3^2).
   \ee
   Let $y=(Q_2^2+Q_3^2)^{\frac{1}{2}}$. Then, by the estimates for $F_{ij}$ above and \eqref{vorticitynormal}, $y$ satisfies:
   \bee
   \dfrac{d}{dt}y+2\be_1ay\leq(1+\ep)^{\frac{1}{4}}y+\ep^{\frac{1}{4}}+|G_2\circ\fai|+|G_3\circ\fai|,
   \ee
   which can be integrated to obtain
   \bee\label{vorticity3}
   \bes
   y(t)&\leq e^{-2\be_1a(t+\ep)}y(-\ep)+\int_{-\ep}^te^{-2\be_1a(t-t')}(\ep^{\frac{1}{4}}+|G_2\circ\fai|+|G_3\circ\fai|) dt'\\
   &:= e^{-2\be_1a(t+\ep)}y(-\ep)+\text{\Rmnum{1}}+\text{\Rmnum{2}}+\text{\Rmnum{3}},
   \end{split}
   \ee
   where $\text{\Rmnum{1}},\text{\Rmnum{2}},\text{\Rmnum{3}}$ are the integrals given in order. For the term \Rmnum{1}, it can be computed as
   \bee
   \text{\Rmnum{1}}=\int_{-\ep}^{t}\ep^{\frac{1}{4}}e^{-2\be_1a(t-t')}\ dt'=\dfrac{1}{2\be_1a}\ep^{\frac{1}{4}}(1-e^{-2\be_1a(t+\ep)}).
   \ee
   For the term \Rmnum{2}, it follows from the bootstrap assumptions \eqref{bspaK} and \eqref{bsallyW} that $|G_2\circ\fai|\leq 1+\eta^{-\frac{1}{3}}(y)e^s$ and then
   \begin{align*}
   \text{\Rmnum{2}}&\leq \int_{-\ep}^te^{-2\be_1a(t-t')}\ dt'+\int_{-\ep}^te^s\eta^{-\frac{1}{3}}\circ\fai \ dt'\\
   &=\dfrac{1}{2\be_1a}(1-e^{-2\be_1a(t+\ep)})+\int_{-\log\ep}^s\ta(0,\frac{1}{3})\circ\Phi_U\ ds'\\
   &\leq \dfrac{1}{2\be_1a}(1-e^{-2\be_1a(t+\ep)})+\ep^{\frac{1}{8}},
   \end{align*}
   where in the second step, the relation between $\fai$ and $\Phi_U$ is used. The estimates for term \Rmnum{3} is the same as \Rmnum{2}. Therefore, it holds that
   \bee
   y(t)\leq e^{-2\be_1a(t+\ep)}y(-\ep)+\dfrac{1}{2\be_1 a}(1-e^{-2\be_1a(t+\ep)})+\ep^{\frac{1}{8}}.
   \ee
   \begin{prop}\label{boundsforvorticity}(Concentration of vorticity on non-blow up direction)
  The following bounds for the specific vorticity $\mathring{\zeta}$ hold
   \begin{align}
   ||\mathring{\zeta}\cdot N||_{L^{\infty}}&\leq\ep^{\frac{1}{4}},\\
   ||\mathring{\zeta}\cdot T^{\mu}||_{L^{\infty}}&\leq e^{-2\be_1a(t+\ep)}\mathring{\zeta}(-\ep)\cdot T^{\mu}_0+\dfrac{1}{2\be_1 a}(1-e^{-2\be_1a(t+\ep)})+\ep^{\frac{1}{8}}.
   \end{align}
   Moreover, as $t\to T_{\ast}$, $e^{-2\be_1a(t+\ep)}\to (1-\ep\be_1 a)^2$. Therefore,
   \begin{itemize}
   \item If $a<0$, the instantaneous shock will lead to the vorticity increasing;
   \item if $0<a\leq\dfrac{1}{2\be_1\ep}$, the delayed shock will lead to the dissipation of the vorticity;
   \item if $a\geq\dfrac{2}{\be_1\ep}$, then one obtains the global solution and the exponential decay for the vorticity.
   \end{itemize}
   \end{prop}
   \subsection{\textbf{Initial region leads to the shock formation}}
   Recall the Lagrangian flows in self-similar coordinates $\Phi_W$, $\Phi_Z$, $\Phi_U$. The corresponding Lagrangian flows in physical space are defined as $\fai_w$, $\fai_z$, $\fai_u$\footnote{Precisely, $\fai_w$ is defined by $\dfrac{d}{dt}\fai_w(x,t)=v_w\circ\fai_w$, $\fai_w(x,-\ep)=x$ and similar for $\fai_z$ and $\fai_u$.}, which are related with $\Phi_W$, $\Phi_Z$, $\Phi_U$ as follows.
   \begin{align}
   \fai_{w_1}&=e^{-\frac{3}{2}s}\Phi_{W_1},\hs \fai_{w_{\mu}}=e^{-\frac{s}{2}}\Phi_{W_{\mu}},\\
   \fai_{z_1}&=e^{-\frac{3}{2}s}\Phi_{Z_1},\hs \fai_{z_{\mu}}=e^{-\frac{s}{2}}\Phi_{Z_{\mu}},\\
   \fai_{u_1}&=e^{-\frac{3}{2}s}\Phi_{U_1},\hs \fai_{u_{\mu}}=e^{-\frac{s}{2}}\Phi_{U_{\mu}},
   \end{align}
   with the initial condition
   \bee
   \fai_{w}(x_0,t_0)=\fai_z(x_0,t_0)=\fai_u(x_0,t_0)=x_0.
   \ee
   \begin{lem}
   Let $\fai^{x_0}(t):=\fai(x_0,t)$ be either $\fai_w^{x_0}$, $\fai_z^{x_0}$ or $\fai_u^{x_0}$ emanating from the initial point $x_0$. If $\lim_{t\to T_{\ast}}\fai^{x_0}(t)\to 0$, then it holds that
   \bee\label{initialshcokregion1}
   |(x_0)_1|,|(x_0)_{\mu}|\leq-\frac{1}{\be_1 a}\ln(1-\ep\be_1 a)+O(\ep^2).
   \ee
   In particular, for $\fai=\fai_z$ or $\fai_u$, the following bounds hold
   \begin{align}
   |(x_0)_1-(1-\be_2)\kappa_0|,|(x_0)_{\mu}|&\leq-\frac{1}{\be_1 a}\ln(1-\ep\be_1 a)+O(\ep^2),\ \text{for}\ \fai=\fai_z,\label{initialshockregionz}\\
   |(x_0)_1-(1-\be_1)\kappa_0|,|(x_0)_{\mu}|&\leq-\frac{1}{\be_1 a}\ln(1-\ep\be_1 a)+O(\ep^2),\ \text{for}\ \fai=\fai_u.\label{initialshockregionu}
   \end{align}
   That is, if one changes the initial data outside the above region, the shock location and shock time won't change.
   \end{lem}
   \begin{pf}
   Letting $t\to T_{\ast}$ in the identity $\fai(x_0,t)-x_0=\int_{-\ep}^t\pa_t\fai(t')\ dt'$ yields
   \bee\label{x0vorticity}
   x_0=-\int_{-\ep}^{T_{\ast}}\pa_t\fai(t')\ dt'.
   \ee
   Note that
   \begin{align}
   \pa_t\fai_1&=(1-\dot{\tau})e^{-\frac{s}{2}}(\pa_s\Phi_{1}-\frac{3}{2}\Phi_{1}),\\
   \pa_t\fai_{\mu}&=(1-\dot{\tau})e^{\frac{s}{2}}(\pa_s\Phi_{\mu}-\frac{1}{2}\Phi_{\mu}).
   \end{align}
   Then,
   \begin{itemize}
   \item for the case $\fai=\fai_w$, it follows from \eqref{defGWGZGU} and Lemma\eqref{dampingtermsestimates} that
         \begin{align*}
         |\pa_t\fai_1|&\leq e^{-\frac{s}{2}}|W|+e^{-\frac{1}{2}}|G_W|\leq 1,\\
         |\pa_t\fai_{\mu}|&\leq e^{\frac{s}{2}}|h_W^{\mu}|\leq 1.
         \end{align*}
         Therefore, it follows from \eqref{x0vorticity} that
         \bee
         |x_0|\leq|\pa_t\fai|(T_{\ast}+\ep)\leq-\frac{1}{\be_1 a}\ln(1-\ep\be_1 a)+O(\ep^2).
         \ee
   \item For the case $\fai=\fai_z$, It follows from \eqref{defGWGZGU} and Lemma\eqref{dampingtermsestimates} that
         \bee
         |\pa_t\fai_1-(1-\be_2)\kappa_0|\leq1,\hs |\pa_t\fai_{\mu}|\leq 1.
         \ee
         Integrating \eqref{x0vorticity} from $-\ep$ to $T_{\ast}$ leads to \eqref{initialshockregionz}. The proof for \eqref{initialshockregionu} is similar as for \eqref{initialshockregionz}.
   \end{itemize}
   \end{pf}
\section{\textbf{Recover the Bootstrap assumptions for the self-similar variables}}\label{section8}

\subsection{\textbf{Recover the Bootstrap assumptions for $\tilde{W}$ in $|y|\leq l$}}
\subsubsection{\textbf{The case for $|\ga|=4$}}
\hs Take $|\ga|=4$ in the equation \eqref{pagatildew}. The damping term becomes
\bee\label{damping4tildew}
D_{\tilde{W}}^{(\ga)}=\dfrac{3}{2}+\ga_1+\be_{\tau}J(\pa_1\bar{W}+\ga_1\pa_1W)\geq\dfrac{3}{2}-(1+C\ep)+\ga_1-\ga_1(1+C\ep)\geq\dfrac{1}{3}.
\ee
Along the Lagrangian flow $\Phi_W$, \eqref{pagatildew} becomes
  \begin{align*}
  \dfrac{d}{ds}\left(\pa^{\ga}\tilde{W}\circ\Phi_W^{y_0}(s)\right)+D_{\tilde{W}}^{(\ga)}\circ\Phi_{\tilde{W}}^{y_0}\pa^{\ga}\tilde{W}\circ\Phi_W^{y_0}&=\tilde{F}_W^{(\ga)}\circ
  \Phi_W^{y_0},
  \end{align*}
  which implies
  \begin{align*}
  \dfrac{d}{ds}\left(\pa^{\ga}\tilde{W}\circ\Phi_W(s)e^{\int_{s_0}^sD_{\tilde{W}}^{(\ga)}\circ\Phi_W^{y_0} ds'}\right)&=
  e^{\int_{s_0}^sD_{\tilde{W}}^{(\ga)}\circ\Phi_W^{y_0} ds'}\tilde{F}_{W}^{(\ga)}\circ\Phi_W^{y_0}.
  \end{align*}
  Integrating from $s_0=-\log\ep$ to $s$ yields
  \begin{align}
  |\pa^{\ga}\tilde{W}\circ\Phi_W^{y_0}(s)|\leq e^{-\int_{s_0}^sD_{\tilde{W}}^{(\ga)}\circ\Phi_W\ ds'}|\pa^{\ga}\tilde{W}(y_0)|+\int_{s_0}^s\tilde{F}_W^{(\ga)}\circ\Phi_We^{-\int_{s'}^sD_{\tilde{W}}^{(\ga)}\circ\Phi_W ds''}ds'.
  \end{align}
  Hence, it follows from Lemma\ref{forcingtermsestimates} and \eqref{damping4tildew} that
  \bee\label{pa4tildew}
  |\pa^{\ga}\tilde{W}\circ\Phi_W(s)|\leq\ep^{\frac{1}{8}}e^{-\frac{3}{2}(s-s_0)}+\ep^{\frac{1}{10}}(\log M)^{|\check{\ga}|-1}<\ep(\log M)^{|\check{\ga}|},
  \ee
   which recovers the Bootstrap assumption\eqref{bapa4tildew} by choosing $M$ large enough.
  \subsubsection{\textbf{The case for $|\ga|\leq 3$}}
  \hs In this case, the sign of the damping terms is not definite. As we are examining estimates in the vicinity of $0$, the behavior of $\pa^{\gamma}\tilde{W}$ is predominantly influenced by the value of $\pa^{\gamma}\tilde{W}(y=0).$
  It follows from \eqref{constraints} and \eqref{barWtaylor} that
  \bee\label{pa0120tildew}
  \tilde{W}(0,s)=\nabla\tilde{W}(0,s)=\nabla^2\tilde{W}(0,s)=0.
  \ee
  Evaluating equation\eqref{pagatildew} at $y=0$ for $|\ga|=3$ yields
  \bee
  \pa_s(\pa^{\ga}\tilde{W})^0=\tilde{F}_W^{(\ga),0}-G_W^0(\pa_1\pa^{\ga}\tilde{W})^0-h_W^{\mu,0}(\pa_{\mu}\pa^{\ga}\tilde{W})^0-(1+\ga_1)(1-\be_{\tau})(\pa^{\ga}\tilde{W})^0.
  \ee
  Then, it follows from Lemma\eqref{forcingtermsestimates}, \eqref{boundsdotQbetau} and \eqref{pa4tildew} which evaluated at $y=0$ that
  \bee
  \pa_s(\pa^{\ga}\tilde{W})^0\leq \ta(\frac{1}{2}-\frac{4}{2m-7},0)+\ep^{\frac{1}{10}}(\log M)^4\ta(\frac{1}{2},0)+CM\ep^{\frac{1}{4}}\ta(1,0)\leq2\ta(\frac{1}{2}-\frac{4}{2m-7},0).
  \ee
  Integrating from $s_0=-\log\ep$ to $s$ leads
  \bee
  (\pa^{\ga}\tilde{W})^0\leq \pa^{\ga}\tilde{W}(0,-\log\ep)+\int_{s_0}^s2\ta(\frac{1}{2}-\frac{4}{2m-7},0) ds'\leq3\ep^{\frac{1}{2}-\frac{4}{2m-7}}<\ep^{\frac{1}{4}}
  \ee
  with $|\ga|= 3$ for sufficiently small $\ep$. Therefore, for $|\ga|=3$, it holds that
  \begin{align*}
  |\pa^{\ga}\tilde{W}|&\leq (\pa^{\ga}\tilde{W})^0+|y_1|\cdot|\pa_1\pa^{\ga}\tilde{W}|+|y_{\mu}|\cdot|\pa_{\mu}\pa^{\ga}\tilde{W}|\\
  &\leq\ep^{\frac{1}{3}}+\ep^{\frac{1}{10}}(\log M)^3|y_1|+\ep^{\frac{1}{10}}(\log M)^4|y_{\mu}|<M\ep^{\frac{1}{4}}+\ep^{\frac{1}{10}}|y|(\log M)^4,
  \end{align*}
  which recovers the Bootstrap assumption \eqref{pa3tildew} for $|\ga|=3$. The estimates for the case $|\ga|\leq2$ can be recovered in the same procedure due to \eqref{pa0120tildew}.
\subsection{\textbf{Weighted estimates for recovering the Bootstrap assumptions for $(W,Z,A,K)$}}
\hs We introduce the following framework to recover the bootstrap assumptions for $(W,Z,A_v,K)$. Let $E$ be any one of $Z,A_v,K$ or $W,\tilde{W}$ with $|y|\geq l$. Note that $E$ satisfies the following type equation:
\bee
\pa_s\pa^{\ga}E+D_E^{(\ga)}\pa^{\ga}E+(V_E\cdot\nabla)\pa^{\ga}E=F_E^{(\ga)}.
\ee
Multiplying $\ta(\a,\be)$ on both sides leads to
\begin{align}
\pa_s(\ta(\a,\be)\pa^{\ga}E)+D_E^{(\ga)'}\ta(\a,\be)\pa^{\ga}E+(V_E\cdot\nabla)\ta(\a,\be)\pa^{\ga}E&=F_E^{(\ga)'}:=\ta(\a,\be)F_E^{(\ga)},\label{framework}
\end{align}
where $D_E^{(\ga)'}$ is the damping term given by
\begin{align*}
D_E^{(\ga)'}&=D_{E}^{(\ga)}+\a+\be\eta^{-1}(y)V_E\cdot(2y_1,6y_{\mu}|y_{\mu}|^4)\\
&=D_E^{(\ga)}+\a+3\be-3\be\eta^{-1}(y)+2\be\eta^{-1}(y)(y_1g_E+3y_{\mu}|y_{\mu}|^4h_E^{\mu}).
\end{align*}
It follows from Lemma\ref{keylemma},\ref{dampingtermsestimates} and the bootstrap assumption\eqref{bsallyW} that
\bee
e^{\int_{s_0}^s|2\be\eta^{-1}(y)(y_1g_E+3y_{\mu}|y_{\mu}|^4h_E^{\mu})\circ\Phi_E|\ ds}\les 1.
\ee
Therefore, applying Gronwall inequality to \eqref{framework} yields
\begin{align*}
\ta(\a,\be)\pa^{\ga}E(y,s)&=\ta(\a,\be)(s_0)\pa^{\ga}E(s_0)e^{-\int_{s_0}^sD_E^{(\ga)'}\circ\Phi_E\ ds'}\\
&+e^{-\int_{s_0}^sD_E^{(\ga)'}\circ\Phi_E\ ds'}\int_{s_0}^sF_{E}^{(\ga)'}\circ\Phi_E e^{\int_{s_0}^{s'}D_E^{(\ga)'}\circ\Phi_E ds''} ds'\\
&\leq\ta(\a,\be)(s_0)|\pa^{\ga}E(s_0)|e^{-\int_{s_0}^s(D_E^{(\ga)}+\a+3\be-3\be\eta^{-1})\circ\Phi_E\ ds'}\\
&+\int_{s_0}^sF_E^{(\ga)'}\circ\Phi_Ee^{-\int_{s'}^s(D_E^{(\ga)}+\a+3\be-3\be\eta^{-1})\circ\Phi_E\ ds''}\ ds',
\end{align*}
which implies
\bee\label{framework2}
\bes
|\pa^{\ga}E(y,s)|&\leq\frac{\ta(\a,\be)(s_0)}{\ta(\a,\be)(s)}|\pa^{\ga}E(s_0)|e^{-\int_{s_0}^s(D_E^{(\ga)}+\a+3\be-3\be\eta^{-1})\circ\Phi_E\ ds'}\\
&+
\ta^{-1}(\a,\be)(s)\int_{s_0}^sF_E^{(\ga)'}\circ\Phi_Ee^{-\int_{s'}^s(D_E^{(\ga)}+\a+3\be-3\be\eta^{-1})\circ\Phi_E\ ds''}\ ds'.
\end{split}
\ee
Note that for $\be\leq 0$, the term $-3\be\eta^{-1}$ in \eqref{framework2} does not contribute to the estimates for $\pa^{\ga}E$ so one can omit this term.
\begin{itemize}
\item For $\pa^{\ga}E=\pa^{\ga}Z,\pa^{\ga}A_v$ or $\pa^{\ga}K$ with $\pa^{\ga}K\neq\pa^2_1K$, $\pa^{\ga}A_v\neq\pa_1A_v$, take $s_0=-\log\ep$, $\be=0$ and $\a\geq-\dfrac{3\ga_1+\ga_2+\ga_3}{2}$\footnote{Actually $\a$ equals to the exactly value of the exponent in the bootstrap assumptions.}. Then, it follows from Lemma\ref{keylemma} that
    \bee
    e^{-\int_{s_0}^s(D_E^{(\ga)}+\a+3\be-3\be\eta^{-1})\circ\Phi_E(s')\ ds'}\leq e^{-\int_{s_0}^s\be_{\tau}J\pa_1W\circ\Phi_E ds'}\leq 1.
    \ee
    Furthermore, due to Lemma\ref{forcingtermsestimates} and \ref{keylemma}, it holds that
    \bee
    \int_{s_0}^sF_E^{(\ga)'}\circ\Phi_Ee^{-\int_{s'}^s(D_E^{(\ga)}+\a+3\be-3\be\eta^{-1})\circ\Phi_E\ ds''}\ ds'\les 1.
    \ee
    \begin{remark}
    The above implicit constant may depend on $M$ (for example, in the case of $\pa^{\ga}E=\pa^2_1Z$, $\int_{s_0}^sF_E^{(\ga)'}\circ\Phi_Ee^{-\int_{s'}^s(D_E^{(\ga)}+\a+3\be-3\be\eta^{-1})\circ\Phi_E\ ds''}\ ds'\leq M^{\frac{1}{2}}<M),$ but it strictly less than the constant in the corresponding bootstrap assumptions.
    \end{remark}
    Therefore,
    \begin{align*}
    |\pa^{\ga}E(y,s)|&\leq\ta^{-1}(\a,\be)\ta(\a,\be)(s_0)|\pa^{\ga}E(-\log\ep)|\\
    &+\ta^{-1}(\a,\be)\int_{s_0}^sF_E^{(\ga)'}\circ\Phi_Ee^{-\int_{s'}^s(D_E^{(\ga)}+\a+3\be-3\be\eta^{-1})\circ\Phi_E\ ds''}\ ds'\\
&\leq e^{-\a s}\ep^{-\a}|\pa^{\ga}E(-\log\ep)|+Ce^{-\a s},
    \end{align*}
    which recovers the bootstrap assumptions by standard computation.
    \item For $\pa^{\ga}E=\pa_1^2K$, take $\a=-2,\be=-\dfrac{1}{15}$. Then, it follows from Lemma\ref{keylemma} and \ref{forcingtermsestimates} that
    \begin{align}
    |\pa_1^{2}K(y,s)|&\leq\ta(2,\dfrac{1}{15})(y,s)\ep^{\frac{1}{4}}\eta^{-\frac{1}{15}}(y)+\ta(2,\dfrac{1}{15})\int_{s_0}^s\ta(-2,-\dfrac{1}{15})F_K^{(\ga)}\circ\Phi_U\ ds'\\
    &\leq \ep^{\frac{1}{4}}\ta(2,\dfrac{1}{15})+\ep^{\frac{1}{4}}\ta(2,\dfrac{1}{15})<\ep^{\frac{1}{8}}\ta(2,\dfrac{1}{15}),
    \end{align}
    which recovers the bootstrap assumptions for $\pa_1^2K$.
    \item For $E=W$ or $\tilde{W}$ with $|y|\geq l$, we divide the region of $y$ into: $l\leq|y|\leq\mcl$ and $|y|\geq\mcl$ by following the remark\ref{wtildew}. 
          For the initial value in\eqref{framework2}, due to the estimate\eqref{phiw}, it holds that
          \begin{itemize}
          \item for $y\in[l,\mcl]$ and $s\geq-\log\ep$, there exist a pair $y_0\in[l,\mcl]$ and $s_0\in[-\log\ep,s]$ such that $\Phi_E^{y_0}(s)=y$. Moreover,
                \begin{align*}
                &\text{if}\ |y_0|=l,\ \text{then}\ s_0>-\log\ep;\\
                &\text{if}\ |y_0|\geq l,\ \text{then}\ s_0=-\log\ep.
                \end{align*}
          \item For $y\geq\mcl$ and $s\geq-\log\ep$, there exist a pair $y_0\geq\mcl$ and $s_0\in[-\log\ep,s]$ such that $\Phi_E^{y_0}(s)=y$. Moreover,
                \begin{align*}
                &\text{if}\ |y_0|=\mcl,\ \text{then}\ s_0>-\log\ep;\\
                &\text{if}\ |y_0|\geq \mcl,\ \text{then}\ s_0=-\log\ep.
                \end{align*}
          \end{itemize}
          In the following, we take $\ga=(2,0,0)$ as an example to recover the assumption for $\pa_1^{2}W$ on $|y|\geq l$ and the proofs for others are similar. In this case, taking $\ta(\a,\be)=\ta(0,-\frac{1}{3})\fe^{-\frac{1}{2}}$ with $\fe=\ta(0,\frac{1}{2})+\ta(\frac{4}{5},-\frac{1}{5})$ in \eqref{framework} leads to
          \bee\label{framepaga12W}
          \pa_s(\ta(0,-\frac{1}{3})\fe^{-\frac{1}{2}}\pa^{\ga}W)+D_W^{(\ga)'}(\ta(0,-\frac{1}{3})\fe^{-\frac{1}{2}})\pa^{\ga}W+(V_W\cdot\nabla)(\ta(0,-\frac{1}{3})\fe^{-\frac{1}{2}})\pa^{\ga}W=F_W^{(\ga)'},
          \ee
          where
          \begin{align*}
          D_W^{(\ga)'}&=D_W^{(\ga)}-\frac{2}{5}\fe^{-1}\ta(\frac{4}{5},-\frac{1}{5})-(1-\eta^{-1}+\frac{2}{3}\eta^{-1}(y_1g_W+3y_{\mu}|y_{\mu}|^4h_W^{\mu})\\
          &+\fe^{-1}[\ta(\frac{4}{5},-\frac{1}{5})(\frac{3}{10}-\frac{3}{10}\eta^{-1}+\frac{1}{5}\eta^{-1}(y_1g_W+3y_{\mu}|y_{\mu}|^4h_W^{\mu})\\
          &-\ta(0,\frac{1}{2})(\frac{3}{4}-\frac{3}{4}\eta^{-1}+\frac{1}{2}\eta^{-1}(y_1g_W+3y_{\mu}|y_{\mu}|^4h_W^{\mu}),\\
          F_W^{(\ga)'}&=(\ta(0,-\frac{1}{3})\fe^{-\frac{1}{2}})F_W^{(\ga)}.
          \end{align*}
          Note that $e^{\int_{s_0}^s|\ta(0,1)(y_1g_W+3y_{\mu}|y_{\mu}|^4h_W^{\mu})\circ\Phi_W|\ ds}\les 1$ and $-D_W^{(\ga)'}\leq-\frac{1}{20}+\frac{41}{30}\eta^{-1}(y_1g_W+3y_{\mu}|y_{\mu}|^4h_W^{\mu})+3\be_{\tau}J\pa_1W$, which implies $e^{\int_{s_0}^s-D_W^{(\ga)'}}\les1.$ Then, \eqref{framepaga12W} leads to
          \bee
          \bes
          |\pa^{\ga}W(y,s)|&\les\frac{\ta(0,-\frac{1}{3})\fe^{-\frac{1}{2}}(s_0)}{\ta(0,-\frac{1}{3})\fe^{-\frac{1}{2}}(s)}|\pa^{\ga}W(s_0)|+\ta(0,\frac{1}{3})\fe^{\frac{1}{2}}(s)
          \int_{s_0}^s\ta(0,-\frac{1}{3})\fe^{-\frac{1}{2}}F_W^{(\ga)}\circ\Phi_We^{-\int_{s'}^sD_W^{(\ga)'}\circ\Phi_W} ds'.
          \end{split}
          \ee
           Then, it follows from  Lemma\ref{forcingtermsestimates} and \ref{keylemma} that
          \bee\label{pa12w}
          |\pa^{\ga}W(y,s)|\leq\eta^{-\frac{1}{3}}(y)\fe^{\frac{1}{2}}(y,s)|\eta^{\frac{1}{3}}(y_0)\fe^{-\frac{1}{2}}(y_0,s_0)\pa^{\ga}W(y_0,s_0)|+M^{\frac{1}{2}} e^{-\frac{s_0}{12}}\fe^{\frac{1}{2}}(y,s)\eta^{-\frac{1}{3}}(y).
          \ee
          due to the estimate
          \bee
          \int_{s_0}^s\ta(0,-\frac{1}{3})\fe^{-\frac{1}{2}}F_W^{(\ga)}\circ\Phi_W\les \int_{s_0}^sM^{\frac{1}{2}}\ta(\frac{3}{5}-\frac{7}{2m-4},-\frac{1}{6})\circ\Phi_W \les M^{\frac{1}{2}}e^{-\frac{s_0}{12}}.
          \ee
          There are the following cases.
          \begin{itemize}
          \item In the region $l\leq|y|\leq\mcl$,
                \been[(1)]
                \item if $|y_0|\geq l$, then $s_0=-\log\ep$ and by the initial data assumptions, it holds that
                      \bee
                      |\pa^{\ga}W(y,s)|\leq\eta^{-\frac{1}{3}}(y)\fe^{\frac{1}{2}}(y,s)+M^{\frac{1}{2}}\ep^{\frac{1}{12}}\eta^{-\frac{1}{3}}(y)\fe^{\frac{1}{2}}(y,s)< M^{\frac{1}{3}}\eta^{-\frac{1}{3}}\fe^{\frac{1}{2}};
                      \ee
                \item while if $|y_0|=l$, then $s_0\geq-\log\ep$, and one are able to obtain by using \eqref{pa3tildew} for $|\ga|=2$
                      \begin{align*}
                      |\pa^{\ga}W(y,s)|&\leq\eta^{-\frac{1}{3}}(y)\fe^{\frac{1}{2}}(y,s)\left(\eta^{\frac{1}{3}}(y_0)|\pa^{\ga}\bar{W}(y_0)|
                      +\eta^{\frac{1}{3}}(y_0)|\pa^{\ga}\tilde{W}(y_0)|\right)+M^{\frac{1}{2}}\ep^{\frac{1}{12}}\eta^{-\frac{1}{3}}\fe^{\frac{1}{2}}\\
                      &\leq(\ep^{\frac{1}{11}}l^{10}+\frac{3}{4})\eta^{-\frac{1}{3}}\fe^{\frac{1}{2}}+M^{\frac{1}{2}}\ep^{\frac{1}{12}}\eta^{-\frac{1}{3}}\fe^{\frac{1}{2}}
                      <M^{\frac{1}{3}}\eta^{-\frac{1}{3}}(\eta^{-\frac{1}{2}}+e^{-\frac{4}{5}s}\eta^{\frac{2}{5}})^{\frac{1}{2}}.
                      \end{align*}
                \een
                In conclusion, for $l\leq|y|\leq\mcl$, it holds that
                \bee\label{pa12w2}
                |\pa^{\ga}W(y,s)|\leq M^{\frac{1}{4}}\eta^{-\frac{1}{3}}(\eta^{-\frac{1}{2}}+e^{-\frac{4}{5}s}\eta^{\frac{2}{5}})^{\frac{1}{2}}
                <M^{\frac{1}{3}}\eta^{-\frac{1}{3}}(\eta^{-\frac{1}{2}}+e^{-\frac{4}{5}s}\eta^{\frac{2}{5}})^{\frac{1}{2}}.
                \ee
          \item In the region $|y|\geq\mcl$,
                \been[(1)]
                \item if $|y_0|\geq \mcl$, then $s_0=-\log\ep$ and by the initial data assumptions, it holds that
                      \bee
                      |\pa^{\ga}W(y,s)|\leq\eta^{-\frac{1}{3}}(y)\fe^{\frac{1}{2}}(y,s)+M^{\frac{1}{2}}\ep^{\frac{1}{12}}\eta^{-\frac{1}{3}}(y)\fe^{\frac{1}{2}}(y,s)< M^{\frac{1}{3}}\eta^{-\frac{1}{3}}(\eta^{-\frac{1}{2}}+e^{-\frac{4}{5}s}\eta^{\frac{2}{5}})^{\frac{1}{2}};
                      \ee
                \item while if $|y_0|=\mcl$, then $s_0\geq-\log\ep$ and it follows from \eqref{pa12w2} that
                      \begin{align*}
                      |\pa^{\ga}W(y,s)|&\leq\eta^{-\frac{1}{3}}(y)\fe^{\frac{1}{2}}(y,s)|\eta^{\frac{1}{3}}(y_0)\pa^{\ga}W(y_0,s_0)|
                      +M^{\frac{1}{2}}\ep^{\frac{1}{12}}\eta^{-\frac{1}{3}}\fe^{\frac{1}{2}}(y,s)\\
                      &\leq2M^{\frac{1}{4}}\eta^{-\frac{1}{3}}(\eta^{-\frac{1}{2}}+e^{-\frac{4}{5}s}\eta^{\frac{2}{5}})^{\frac{1}{2}}<M^{\frac{1}{3}}\eta^{-\frac{1}{3}}(
                      \eta^{-\frac{1}{2}}+e^{-\frac{4}{5}s}\eta^{\frac{2}{5}})^{\frac{1}{2}}.
                      \end{align*}
                \een
          \end{itemize}
        Hence, one recovers the bootstrap assumptions for $\pa_1^2W$ in the region $|y|\geq l$. 
\end{itemize}
To deal with the case $\pa^{\ga}E=\pa_1A_v$, one uses the following lemma. 
\begin{lem}
The following identities for $\pa_1A_v$ hold.
\begin{align}
e^{\frac{3}{2}s}J\pa_1A_2&=(\a e^{-\frac{K}{2}}S)^{\frac{1}{\a}}\Omega\cdot T^3+\dfrac{1}{2}T^2_{\mu}\left(\pa_{\mu}W+e^{\frac{s}{2}}\pa_{\mu}Z\right)-e^{\frac{s}{2}}N_{\mu}\pa_{\mu}A_2\\
&-\left[(U\cdot N)N+A_vT^v\right]\cdot(N_i\pa_{x_i}T^2-T^2_i\pa_{x_i}N),\label{boundsfora2}\\
e^{\frac{3}{2}s}J\pa_1A_3&=-(\a e^{-\frac{K}{2}}S)^{\frac{1}{\a}}\Omega\cdot T^2+\dfrac{1}{2}T^3_{\mu}\left(\pa_{\mu}W+e^{\frac{s}{2}}\pa_{\mu}Z\right)-e^{\frac{s}{2}}N_{\mu}\pa_{\mu}A_3\\
&+\left[(U\cdot N)N+A_vT^v\right]\cdot(N_i\pa_{x_i}T^3-T^3_i\pa_{x_i}N).\label{boundsfora3}
\end{align}
As a consequence, it follows from Proposition\ref{boundsforvorticity} and bootstrap assumptions\eqref{pagaZ}-\eqref{pagaallyW} that
\bee
e^{\frac{3}{2}s}|\pa_1A_v|\les\kappa_0\ep^{\frac{1}{10}}+(1+M^{\frac{1}{2}}\ep^{\frac{1}{2}})+(\kappa_0+\ep^{\frac{1}{4}})\ep< M^{\frac{1}{4}},
\ee
which recovers the bootstrap assumptions for $\pa_1A_v$.
\end{lem}
\begin{pf}
We only establish \eqref{boundsfora3} and the computation for \eqref{boundsfora2} is same. Note that in the orthonomral basis $(N,T^2,T^3)$,
\bee
curl_{\tilde{x}}\mr{u}\cdot T^2=\pa_{T^3}\mr{u}(x,t)\cdot N-\pa_N\mr{u}(x,t)\cdot T^3,
\ee
and
\begin{align}
\pa_{T^3}\mr{u}&=T^3_1\pa_{x_1}\mr{u}-T^3_{v}f_{,v}\pa_{x_1}\mr{u}+T^3_v\pa_{x_v}\mr{u}=JN\cdot T^3\pa_{x_1}\mr{u}+T^3_v\pa_{x_v}\mr{u}=T^3_v\pa_{x_v}\mr{u}(x,t),\\
\pa_{N}\mr{u}&=N_1\pa_{x_1}\mr{u}-N_{v}f_{,v}\pa_{x_1}\mr{u}+N_v\pa_{x_v}\mr{u}=JN\cdot N\pa_{x_1}\mr{u}+N_v\pa_{x_v}\mr{u}=J\pa_{x_1}\mr{u}+N_v\pa_{x_v}\mr{u}(x,t).
\end{align}
Therefore, it follows that
\begin{align}
&(\a e^{-\frac{K}{2}}S)^{\frac{1}{\a}}\Omega\cdot T^2=curl_{\tilde{x}}\mr{u}\cdot T^2=T^3_v\pa_{x_v}\mr{u}(x,t)-J\pa_{x_1}\mr{u}-N_v\pa_{x_v}\mr{u}(x,t)\\
&=\dfrac{1}{2}T^3_v\pa_{x_v}(w+z)-\mr{u}\cdot T^3_v\pa_{x_v}N-J\pa_{x_1}a_3-N_v\pa_{x_v}a_3+\mr{u}\cdot N_v\pa_{x_v}T^3\\
&=\dfrac{1}{2}T^3_v\pa_{y_v}(e^{-\frac{s}{2}}W+Z)-e^{\frac{3}{2}s}J\pa_{y_1}A_3-N_ve^{\frac{s}{2}}\pa_{y_v}A_3+\left[(U\cdot N)N+A_vT^v\right]\cdot(N_i\pa_{x_i}T^3-T^3_i\pa_{x_i}N).
\end{align}
\end{pf}
\section{\textbf{The energy estimates}}\label{section9}
\hs By far, to complete the proof for the main theorem, it suffices to prove the energy estimate prop\eqref{energy}, which only relies on the Bootstrap assumptions. To this end, introduce the following semi-norm:
\bee
E_m^2(s)=\sum_{\ga=m}\lam^{\check{\ga}}\left(||\pa^{\ga}U||^2_{L^2}+||H\pa^{\ga}\mathcal{P}||_{L^2}^2+\kappa_0^2||\pa^{\ga}H||_{L^2}^2\right),
\ee
where $\lam\in(0,1)$ is a constant to absorb the various coefficients. Obviously,
\bee
E_m^2\sim(||U||_{\dot{H}^m}^2+||\mathcal{P}||_{\dot{H}^m}^2+||H||_{\dot{H}^m}^2),
\ee
due to the estimate $|H-1|\leq\ep.$ Taking $\pa^{\ga}$ for $|\ga|=m$ to the $(U,\mathcal{P},H)$ system \eqref{UPHsystem} yields the following equations
\begin{align}
\begin{split}
&\pa_s(\pa^{\ga}U_i)-2\be_1\be_{\tau}e^{-s}\dot{Q}_{ij}(\pa^{\ga}U_j)+(V_U\cdot\nabla)\pa^{\ga}U_i+D_{\ga}\pa^{\ga}U_i+2\be_1\be_{\tau}e^{-s}a\pa^{\ga}U_i\\
&+2\be_3\be_{\tau}\mathcal{P}H^2(JN_ie^{\frac{s}{2}}\pa_1\pa^{\ga}\mathcal{P}+\da^{iv}e^{-\frac{s}{2}}\pa_v\pa^{\ga}\mathcal{P})+2\be_3\be_{\tau}H^2\pa^{\ga}\mathcal{P}
(JN_ie^{\frac{s}{2}}\pa_1\mathcal{P}+\ga_1JN_ie^{\frac{s}{2}}\pa_1\mathcal{P})=F_{U_i}^{(\ga)},\label{pagau}
\end{split}\\
\begin{split}
&\pa_s\pa^{\ga}\mathcal{P}+(V_U\cdot\nabla)\pa^{\ga}\mathcal{P}+D_{\ga}\pa^{\ga}\mathcal{P}+2\be_3\be_{\tau}\mathcal{P}(JN\cdot e^{\frac{s}{2}}\pa_1\pa^{\ga}U+e^{-\frac{s}{2}}
\pa_{\mu}\pa^{\ga}U_{\mu})\\
&+2\be_3\be_{\tau}e^{\frac{s}{2}}\pa^{\ga}\mathcal{P}JN\cdot\pa_1U+2\ga_1\be_3\be_{\tau}e^{\frac{s}{2}}\pa_1UJN\cdot\pa^{\ga}U=F_{\mathcal{P}}^{(\ga)},\label{pagap}
\end{split}\\
&\pa_s\pa^{\ga}H+(V_U\cdot\nabla)\pa^{\ga}H+D_{\ga}\pa^{\ga}H=F_{H}^{(\ga)},\label{pagah}
\end{align}
where
\bee
D_{\ga}=\dfrac{1}{2}|\ga|+\ga_1(1+\pa_1g_U),
\ee
and the forcing terms are given as
\begin{align}
\begin{split}
F_{U_i}^{(\ga)}&=-[\pa^{\ga},V_U\cdot\nabla]U_i+D_{\ga}\pa^{\ga}U_i-2\be_3\be_{\tau}e^{-\frac{s}{2}}[\pa^{\ga},\mathcal{P}H^2]\pa_{v}\mathcal{P}\\
&+2\be_3\be_{\tau}e^{\frac{s}{2}}(JN_iH^2\pa^{\ga}\mathcal{P}\pa_1\mathcal{P}+\ga_1JN_i\pa_1\mathcal{P}\pa^{\ga}\mathcal{P}-[\pa^{\ga},\mathcal{P}H^2JN_i]\pa_1\mathcal{P}),
\end{split}\\
\begin{split}
F_{\mathcal{P}}^{(\ga)}&=-[\pa^{\ga},V_U\cdot\nabla]\mathcal{P}+D_{\ga}\pa^{\ga}\mathcal{P}-2\be_3\be_{\tau}e^{-\frac{s}{2}}[\pa^{\ga},\mathcal{P}H^2]\pa_v\mathcal{P}\\
&+2\be_3\be_{\tau}e^{\frac{s}{2}}(JN\cdot \pa_1U\pa^{\ga}\mathcal{P}+\ga_1\pa_1\mathcal{P}JN_i\pa^{\ga}U_i-[\pa^{\ga},\mathcal{P}H^2JN_i]\pa_1U_i),
\end{split}\\
F_{H}^{(\ga)}&=-[\pa^{\ga},V_U\cdot\nabla]H+D_{\ga}\pa^{\ga}H.
\end{align}
Then, applying the standard Freidrich's energy estimates leads to the following proposition.
\begin{prop}\label{energyinequality}
There exist a universal constant $C$ such that the following energy inequality holds
\begin{align}
&\dfrac{\pa}{\pa s}E^2_m(s)+(|\ga|-C)E_m^2\leq\\
&2\sum_{|\ga|=m}\lam^{\check{\ga}}\int\pa^{\ga}U_i\cdot F_{U_i}^{(\ga)}+H^2\pa^{\ga}\mathcal{P}\cdot F_{\mathcal{P}}^{(\ga)}+\kappa_0^2\pa^{\ga}H\cdot F_{H}^{(\ga)}.
\end{align}
\end{prop}
\begin{pf}
Multiplying $\lam^{|\check{\ga}|}\pa^{\ga}U_i$,  $\lam^{|\check{\ga}|}H^2\pa^{\ga}\mathcal{P}$, $\lam^{|\check{\ga}|}\kappa_0^2\pa^{\ga}H$ to $[\eqref{pagau},\eqref{pagap},\eqref{pagah}$ respectively, adding them up, integrating over $R^3$ and using the skew-symmetric of $\dot{Q}$ lead to
\begin{equation}\label{energy1}
\bes
&\dfrac{\pa}{\pa s}\left[\lam^{|\check{\ga}|}(||\pa^{\ga}U_i||^2_{L^2}+||H^2\pa^{\ga}\mathcal{P}||_{L^2}^2+\kappa_0^2||\pa^{\ga}H||_{L^2}^2)\right]\\
&+\underbrace{\lam^{|\check{\ga}||}\int
(2D_{\ga}-2div V_U)(|\pa^{\ga}U_i|^2+H^2|\pa^{\ga}\mathcal{P}|^2+\kappa_0^2|\pa^{\ga}H|^2)}_{I}\\
&+\underbrace{2\be_1\be_{\tau}e^{-s}a\lam^{|\check{\ga}|}\int|\pa^{\ga}U_i|^2}_{II}\\
&+2\be_3\be_{\tau}\lam^{|\check{\ga}|}\int\mathcal{P}H^2\left[(JN_ie^{\frac{s}{2}}\pa_1\pa^{\ga}
\mathcal{P}\pa^{\ga}U_i+\da^{iv}e^{-\frac{s}{2}}\pa_v\pa^{\ga}\mathcal{P}\pa^{\ga}U_i)+(JN_ie^{\frac{s}{2}}\pa_1\pa^{\ga}U_i\pa^{\ga}\mathcal{P}+e^{-\frac{s}{2}}\pa_{\mu}\pa^{\ga}
U_{\mu}\pa^{\ga}\mathcal{P})\right]\\
&+2\be_3\be_{\tau}\lam^{|\check{\ga}|}\int H^2\left(JN_ie^{\frac{s}{2}}(1+\ga_1)\pa_1\mathcal{P}\pa^{\ga}P\pa^{\ga}U_i+e^{\frac{s}{2}}JN_i(\pa_1U_i|\pa^{\ga}\mathcal{P}|^2
+\ga_1\pa_1\mathcal{P}\pa^{\ga}U_i\pa^{\ga}\mathcal{P})\right)\\
&=2\lam^{\check{\ga}}\int\pa^{\ga}U_i\cdot F_{U_i}^{(\ga)}+H^2\pa^{\ga}\mathcal{P}\cdot F_{\mathcal{P}}^{(\ga)}+\kappa_0^2\pa^{\ga}H\cdot F_{H}^{(\ga)},
\end{split}
\end{equation}
where the second line of \eqref{energy1} denoted to be term I, the third line of \eqref{energy1} denoted to be term II and the fourth line and fifth line denoted to be term III.
\begin{itemize}
\item For the damping term I, it holds that
\begin{align*}
2D_{\ga}-div V_U&=|\ga|+2\ga_1(1+\pa_1g_U)-div V_U\\
&=|\ga|-\dfrac{5}{2}+(2\ga-1)\pa_1g_U+2\ga_1-\pa_{\mu}h_U^{\mu}\\
&\geq|\ga|-\dfrac{5}{2}-\be_1\be_{\tau}(2\ga_1-1)+2\ga_1-e^{-\frac{s}{2}}.
\end{align*}
\item For the term II, it can be bounded directly as
\begin{equation*}
II\geq-2\be_1\be_{\tau}e^{-s}|a|||\pa^{\ga}U_i||_{L^2}^2\geq-\be_{\tau}||\pa^{\ga}U_i||_{L^2}^2.
\end{equation*}
\item For the term III, integrating by parts leads to
\begin{align*}
III&=-2\be_3\be_{\tau}\lam^{|\check{\ga}|}\int\pa^{\ga}\mathcal{P}(2H\pa^{\ga}HJN_ie^{\frac{s}{2}}\pa^{\ga}U_i+\pa_{\mu}e^{-\frac{s}{2}}\pa^{\ga}U_i\\
&+2\be_3\be_{\tau}\lam^{|\check{\ga}|}\int2\ga_1H^2e^{\frac{s}{2}}JN_i\pa^{\ga}U_i\pa^{\ga}\mathcal{P}+e^{\frac{s}{2}}H^2JN_i\pa_1U_i|\pa^{\ga}\mathcal{P}|^2\\
&\geq-(2\be_3\be_{\tau}\ga_1+2\be_3\be_{\tau}+e^{-s})\lam^{|\check{\ga}|}(||\pa^{\ga}U_i||^2_{L^2}+||H^2\pa^{\ga}\mathcal{P}||_{L^2}^2+\kappa_0^2||\pa^{\ga}H||_{L^2}^2).
\end{align*}
\end{itemize}
Taking summation for \eqref{energy1} with $|\ga|=m$ and combing the above results yield
\begin{equation}
\bes
&\dfrac{\pa}{\pa s}E^2_m(s)+D_{\ga}^{'}E_m^2\leq\\
&2\sum_{|\ga|=m}\lam^{|\check{\ga}|}\int\pa^{\ga}U_i\cdot F_{U_i}^{(\ga)}+H^2\pa^{\ga}\mathcal{P}\cdot F_{\mathcal{P}}^{(\ga)}+\kappa_0^2\pa^{\ga}H\cdot F_{H}^{(\ga)},
\end{split}
\end{equation}
where
\begin{equation}
\bes
D_{\ga}'&=|\ga|-\dfrac{5}{2}+2\ga_1-\be_1\be_{\tau}(2\ga_1-1)-2\be_3\be_{\tau}\ga_1-(2\be_3+1)\be_{\tau}-e^{-\frac{s}{2}}\\
&=|\ga|-\dfrac{5}{2}+2\ga_1-2\be_{\tau}\ga_1-3\be_3\be_{\tau}-e^{-\frac{s}{2}}\geq|\ga|-C.
\end{split}
\end{equation}
Here the constant $C$ can be taken as a universal constant by choosing $|\ga|$ large enough.
\end{pf}
For the forcing terms in \eqref{energy1}, one has the following lemma.
\begin{lem}\label{Forcingterms}
Let $m $ be sufficiently large and $\lam=\dfrac{\da^2}{16m^2}$. For $0<\da\leq\dfrac{1}{32}$, there exists a universal constant $C$ such that
\begin{align}
2\sum_{|\ga|=m}\lam^{|\check{\ga}|}\int\left|F_{U_i}^{(\ga)}\pa^{\ga}U_i\right|&\leq (5+C\da)E_m^2+e^{-s}M^{4m-1},\label{forcingu}\\
2\sum_{|\ga|=m}\lam^{|\check{\ga}|}\int\left|F_{\mathcal{P}}^{(\ga)}H^2\pa^{\ga}\mathcal{P}\right|&\leq (2+C\da)E_m^2+e^{-s}M^{4m-1},\label{forcingp}\\
2\sum_{|\ga|=m}\lam^{|\check{\ga}|}\int\left|F_{H}^{(\ga)}\kappa_0^2\pa^{\ga}H\right|&\leq (2+C\da)E_m^2+e^{-s}M^{4m-1},\label{forcingh}
\end{align}
by taking $\ep$ sufficiently small in terms of $m,\da,\lam, M,\kappa_0$.
\end{lem}
\begin{pf}
Decompose the forcing terms as
\bee
F_{U_i}^{(\ga)}=F_{U_i}^{(m)}+F_{U_i}^{(m-1)},\hs F_{\mathcal{P}}^{(\ga)}=F_{\mathcal{P}}^{(m)}+F_{\mathcal{P}}^{(m-1)},\hs F_{H}^{(\ga)}=F_{H}^{(m)}+F_{H}^{(m-1)},
\ee
where the index $m$ and $m-1$ represent the terms of order of derivatives equal to $m$ and $\leq m$, respectively. Precisely\footnote{In the following, $A$ is a lower order term (l.o.t) compared with $B$ means $||A(y,s)||_{L^2}=\ep^{\a}e^{-\be s}||O(B(y,s))||_{L^2}$ where $\a,\be\geq0$ and $\a^2+\be^2\neq 0$.},
\begin{align*}
F_{U_i}^{(m)}&=F_{U_i,(1)}^{(m)}+F_{U_i,(2)}^{(m)}+F_{U_i,(3)}^{(m)},
\end{align*}
where
\begin{align*}
\begin{split}
F_{U_i,(1)}^{(m)}&=-(\ga_{\mu}\pa_{\mu}g_U\pa_1\pa^{\ga-e_{\mu}}U_i+\pa^{\ga}g_U\pa_1U_i+\ga_j\pa_jh_U^{\mu}\pa_{\mu}^{\ga-e_j}U_i+\pa^{\ga}h_U^{\mu}\pa_{\mu}U_i)\\
&=-(\ga_{\mu}\pa_{\mu}g_U\pa_1\pa^{\ga-e_{\mu}}U_i+\pa^{\ga}g_U\pa_1U_i+\text{l.o.ts});
\end{split}\\
\begin{split}
F_{U_i,(2)}^{(m)}&=-2\be_3\be_{\tau}e^{-\frac{s}{2}}(\pa^{\ga}(\mathcal{P}H^2)\pa_v\mathcal{P}+\ga_j\pa_j(\mathcal{P}H^2)\pa_v\pa^{\ga-e_j}\mathcal{P});
\end{split}\\
\begin{split}
F_{U_i,(3)}^{(m)}&=-2\be_3\be_{\tau}e^{\frac{s}{2}}(\ga_1\mathcal{P}\pa_1(H^2)JN_i\pa^{\ga}\mathcal{P}+\ga_{\mu}\pa_{\mu}(\mathcal{P}H^2JN_i)\pa^{\ga-e_{\mu}}\pa_1\mathcal{P}+\mathcal{P}\pa^{\ga}(H^2JN_i)\pa_1\mathcal{P})\\
&=-2\be_3\be_{\tau}e^{\frac{s}{2}}(\ga_1\mathcal{P}\pa_1(H^2)JN_i\pa^{\ga}\mathcal{P}+\mathcal{P}\pa^{\ga}(H^2JN_i)\pa_1\mathcal{P})+\text{l.o.ts},
\end{split}
\end{align*}
where $F_{U_i,(2)}^{(m)}$ is a l.o.t compared with $F_{U_i,(3)}^{(m)}$.
\begin{align*}
F_{U_i}^{(m-1)}&=-\sum_{|\be|=2,\be\leq\ga}^{|\ga|-1}C_{\ga}^{\be}(\pa^{\be}g_{U}\pa_1\pa^{\ga-\be}U_i+\pa^{\be}h_{U}^{\mu}\pa_{\mu}\pa^{\ga-\be}U_i)-2\be_3\be_{\tau}e^{-\frac{s}{2}}\sum_{|\be|=2,\be\leq\ga}^{|\ga|-1}C_{\ga}^{\be}\pa^{\be}(\mathcal{P}H^2)\pa_i\pa^{\ga-\be}\mathcal{P}\\
&-2\be_3\be_{\tau}e^{\frac{s}{2}}\sum_{|\be|=2,\be\leq\ga}^{|\ga|-1}C_{\ga}^{\be}\pa^{\be}(\mathcal{P}H^2JN_i)\pa_1\pa^{\ga-\be}\mathcal{P}:=F_{U_i,(1)}^{(m-1)}+F_{U_i,(2)}^{(m-1)}+F_{U_i,(3)}^{(m-1)},
\end{align*}
where $F_{U_i,(1)}^{(m-1)},F_{U_i,(2)}^{(m-1)},F_{U_i,(3)}^{(m-1)}$ are the terms given in order and $F_{U_i,(2)}^{(m-1)}$ is a l.o.t compared with $F_{U_i,(3)}^{(m)}$.
\begin{align*}
F_{\mathcal{P}}^{(m)}&=F_{\mathcal{P},(1)}^{(m)}+F_{\mathcal{P},(2)}^{(m)}+F_{\mathcal{P},(3)}^{(m)},\hs
F_{\mathcal{P},(1)}^{(m)}=-(\ga_{\mu}\pa_{\mu}g_U\pa_1\pa^{\ga-e_{\mu}}\mathcal{P}+\pa^{\ga}g_U\pa_1\mathcal{P}+\text{l.o.t}),\\ F_{\mathcal{P},(3)}^{(m)}&=-2\be_3\be_{\tau}e^{\frac{s}{2}}\ga_{\mu}\pa_{\mu}(\mathcal{P}JN_i)\pa_1\pa^{\ga-e_{\mu}}U_i.
\end{align*}
And $F_{\mathcal{P},(2)}^{(m)}$ is a l.o.t compared with$F_{\mathcal{P},(3)}^{(m)}$.
\begin{align*}
F_{\mathcal{P}}^{(m-1)}&=-\sum_{|\be|=2,\be\leq\ga}^{|\ga|-1}C_{\ga}^{\be}(\pa^{\be}g_{U}\pa_1\pa^{\ga-\be}\mathcal{P}+\pa^{\be}h_{U}^{\mu}\pa_{\mu}\pa^{\ga-\be}\mathcal{P})
-2\be_3\be_{\tau}e^{-\frac{s}{2}}\sum_{|\be|=2,\be\leq\ga}^{|\ga|-1}C_{\ga}^{\be}\pa^{\be}\mathcal{P}\pa_{\mu}\pa^{\ga-\be}U_{\mu}\\
&-2\be_3\be_{\tau}e^{\frac{s}{2}}\sum_{|\be|=2,\be\leq\ga}^{|\ga|-1}C_{\ga}^{\be}\pa^{\be}(\mathcal{P}JN_i)\pa_1\pa^{\ga-\be}U_i
:=F_{\mathcal{P},(1)}^{(m-1)}+F_{\mathcal{P},(2)}^{(m-1)}+F_{\mathcal{P},(3)}^{(m-1)},
\end{align*}
where $F_{\mathcal{P},(1)}^{(m-1)},F_{\mathcal{P},(2)}^{(m-1)},F_{\mathcal{P},(3)}^{(m-1)}$ are the terms given in order and $F_{\mathcal{P},(2)}^{(m-1)}$ is a l.o.t compared with $F_{\mathcal{P},(3)}^{(m-1)}$.
\begin{align*}
F_H^{(m)}&=-(\ga_{\mu}\pa_{\mu}g_U\pa_1\pa^{\ga-e_{\mu}}H+\pa^{\ga}g_U\pa_1H+\text{l.o.t}),\\
F_H^{(m-1)}&=-\sum_{|\be|=2,\be\leq\ga}^{|\ga|-1}C_{\ga}^{\be}(\pa^{\be}g_U\pa_1\pa^{\ga-\be}H+\pa^{\be}h_U^{\mu}\pa_{\mu}\pa^{\ga-\be}H).
\end{align*}
In the following, the proof for \eqref{forcingu} will be given and the proofs for \eqref{forcingp} and \eqref{forcingh} are the same. For the first term in $F_{U_i,(1)}^{(m)}$, note that $\pa_1\pa^{\ga-e_{\mu}}\neq\pa^{\ga}$ and
\bee
\lam^{|\check{\ga}|}\pa_1\pa^{\ga-e_{\mu}}U_i=\lam\lam^{|\check{\ga}|-1}\pa_1\pa^{\ga-e_{\mu}}U_i,
\ee
which implies
\begin{align*}
||\ga_{\mu}\pa_{\mu}g_U\pa_1\pa^{\ga-e_{\mu}}U_i||_{L^2}&\les2\lam^{\frac{1}{2}}\sum_{|\ga|=m}||\pa_{\mu}g_U||_{L^{\infty}}\lam^{\frac{|\check{\ga}|-1}{2}}||\pa_1\pa^{\ga-e_{\mu}}U_i||_{L^2}\lam^{\frac{|\check{\ga}|}{2}}||\pa^{\ga}U_i||_{L^2}\\
&\les(1+\ep^{\frac{1}{4}})\lam^{\frac{1}{2}}E_m^2\leq\da E_m^2.
\end{align*}
For the second term, note that
\bee
\pa^{\ga}g_U=\be_1\be_{\tau}e^{\frac{s}{2}}\pa^{\ga}U\cdot N+\be_1\be_{\tau}e^{\frac{s}{2}}\pa_{\mu}(JN)\pa^{\ga-e_{\mu}}U_i.
\ee
Then, it follows from Lemma\ref{GNS}
\begin{equation}
\sum_{|\ga|=m}\lam^{|\check{\ga}|}e^{\frac{s}{2}}||\pa^{\ga}U\cdot N||_{L^2}||\pa_1U_i||_{L^{\infty}}\les (1+\ep^{\frac{1}{4}})E_m,
\end{equation}
\begin{equation}
\bes
&\sum_{|\ga|=m}\lam^{|\check{\ga}|}e^{\frac{s}{2}}||\pa_1U_i||_{L^{\infty}}||\pa_{\mu}JN||_{L^{\infty}}||\pa^{\ga-e_{\mu}}U_i||_{L^2}\les
\ep e^{-\frac{s}{2}}||\pa^{\ga}U||_{L^2}^{\frac{2m-5}{2m-3}}||U||_{L^{\infty}}^{\frac{2}{2m-3}}\\
&\leq\da E_m+\dfrac{1}{\da}\left(\ep e^{-\frac{s}{2}}||U||_{L^{\infty}}^{\frac{2}{2m-3}}\right)^{\frac{2m-3}{2}}\leq\da E_m+\ep e^{-s}.
\end{split}
\end{equation}
Therefore, this term can be bounded as
\bee
||\pa^{\ga}g_U\pa_1U_i||_{L^2}\leq (1+2\da)E_m^2+\ep e^{-s}
\ee
To sum up, it holds that
\bee
2\sum_{|\ga|=m}\lam^{|\check{\ga}|}\int|F_{U_i,(1)}^{(m)}\cdot\pa^{\ga}U_i|\leq (1+C\da)E_m^2+\ep e^{-s}.
\ee
For the first term in $F_{U_i,(3)}^{(m)}$, it can be bounded by using same technique as estimating $\pa_1\pa^{\ga-e_{\mu}}U_i\pa_{\mu}g_U$ and then
\begin{equation*}
||2\be_3\be_{\tau}e^{\frac{s}{2}}(\ga_1\mathcal{P}\pa_1(H^2)JN_i\pa^{\ga}\mathcal{P}||_{L^2}\les 2\sum_{|\ga|=m}\lam^{\frac{1}{2}}||\pa_{\mu}(\mathcal{P}H^2JN_i)||_{L^{\infty}}\lam^{\frac{|\check{\ga}|}{2}}||\pa^{\ga}U_i||_{L^2}\lam^{\frac{|\check{\ga}|-1}{2}}
||\pa_1\pa^{\ga-e_{\mu}}\mathcal{P}||_{L^2}\leq \da E_m^2.
\end{equation*}
For the second term, note that $|\pa_1\mathcal{P}|\leq (1+\ep)e^{-\frac{s}{2}}$, $|\mathcal{P}|\leq\kappa_0$ and
\bee
\pa^{\ga}(H^2JN_i)=\pa^{\ga}(H^2)JN_i+\pa_{\mu}(JN_i)\pa^{\ga-e_{\mu}}(H^2).
\ee
Then,
\bee
||\pa^{\ga}(H^2)JN_i2\be_3\be_{\tau}e^{\frac{s}{2}}\mathcal{P}\pa_1\mathcal{P}||_{L^2}\leq (1+\ep)\kappa_0||H||_{L^{\infty}}||H||_{\dot{H}^m}\leq \kappa_0 E_m\leq E_m,
\ee
\begin{equation}
\bes
&||\pa^{\ga-e_{\mu}}(H^2)\pa_{\mu}(JN_i)2\be_3\be_{\tau}e^{\frac{s}{2}}\mathcal{P}\pa_1\mathcal{P}||_{L^2}\leq\ep e^{-\frac{s}{2}}(1+\ep)\kappa_0||H||_{L^{\infty}}||H||_{
\dot{H}^{m-1}}\\
\leq &\ep e^{-\frac{s}{2}}\kappa_0||H||_{\dot{H}^{m}}^{\frac{2m-5}{2m-3}}||H||_{L^{\infty}}^{\frac{2}{2m-3}}\leq \kappa_0^{-1}\ep e^{-\frac{s}{2}} E_m^{\frac{2m-5}{2m-3}}\leq\da E_m+\ep e^{-s}.
\end{split}
\end{equation}
Therefore, it can be bounded as
\bee
2\sum_{|\ga|=m}\lam^{|\check{\ga}|}\int|F_{U_i,(3)}^{(m)}\cdot\pa^{\ga}U_i|\leq (1+C\da)E_m^2+\ep e^{-s}.
\ee
To sum up, for the terms of order $m$, they can be bounded as
\bee\label{L2FUim}
2\sum_{|\ga|=m}\lam^{|\check{\ga}|}\int|F_{U_i}^{(m)}\cdot\pa^{\ga}U_i|\leq (2+C\da)E_m^2+\ep e^{-s}.
\ee
For the first term in $F_{U_i,(1)}^{(m-1)}$, the following term will be estimated
\bee
-\sum_{|\be|=2}^{|\ga|-1}C_{\ga}^{\be}\pa^{\be}g_U\pa_1\pa^{\ga-\be}U_i,
\ee
while the estimates for the term $-\sum_{|\be|=2}^{|\ga|-1}C_{\ga}^{\be}\pa^{\be}h_U^{\mu}\pa_{\mu}\pa^{\ga-\be}U_i$ are similar. Note that
\bee
\pa^{\be}g_U=\be_1\be_{\tau}e^{\frac{s}{2}}J\pa^{\be}U\cdot N+\text{l.o.t},
\ee
Then, it follows from Lemma\ref{sobolev2} that
\bee\label{uim-11}
||\sum_{|\be|=2}^{|\ga|-1}C_{\ga}^{\be}\pa^{\be}g_U\pa_1\pa^{\ga-\be}U_i||_{L^2}\les e^{\frac{as}{2}}||\pa^{\ga}U_i||_{L^2}^a
||\pa^{\ga}U_i||_{L^2}^b||D^2g_U||_{L^q}^{1-a}||D^2U||_{L^q}^{1-b}.
\ee
For the second derivatives terms in\eqref{uim-11}, it follows from Lemma\ref{UNSestimates} and \ref{dampingtermsestimates} that
\bee
|D^2g_U|\leq M\eta^{-\frac{1}{6}}(y),\hs |D^2U|\leq Me^{-\frac{s}{2}}\eta^{-\frac{1}{6}}(y),
\ee
which implies
\begin{align*}
||D^2g_U||_{L^q}^{1-a}\leq M^{1-a}||\eta^{-\frac{1}{6}}(y)||_{L^q}^{1-a}\leq M^{1-a},\hs ||D^2U||_{L^q}^{1-b}\leq M^{1-b}e^{-\frac{1-b}{2}s}.
\end{align*}
Therefore,
\bee
||\sum_{|\be|=2}^{|\ga|-1}C_{\ga}^{\be}\pa^{\be}g_U\pa_1\pa^{\ga-\be}U_i||_{L^2}\leq M^{2-(a+b)}e^{\frac{a+b-1}{2}s}E_m^{a+b}.
\ee
Collecting the above results leads to
\begin{equation}\label{L2FUim-11}
\bes
&2\sum_{|\ga|=m}\lam^{|\check{\ga}|}\int|\sum_{|\be|=2}^{|\ga|-1}C_{\ga}^{\be}\pa^{\be}g_U\pa_1\pa^{\ga-\be}U_i\cdot\pa^{\ga}U_i|\leq M^{2-(a+b)}e^{\frac{a+b-1}{2}s}E_m^{1+a+b}\\
\leq&\da E_m^2+\dfrac{1}{\da}\left(M^{2-(a+b)}e^{\frac{a+b-1}{2}s}\right)^{\frac{2}{1-(a+b)}}\leq\da E_m^2+e^{-s}M^{2m-1}.
\end{split}
\end{equation}
For the term $F_{U_i,(3)}^{(m-1)}$, note that
\bee
\pa^{\be}(\mathcal{P}H^2JN_i)=H^2JN_i\pa^{\be}\mathcal{P}+2JN_i\mathcal{P}H\pa^{\be}H+\text{l.o.t}
\ee
and
\bee
|D^2(\mathcal{P}H^2JN_i)|,|D^2\mathcal{P}|\leq Me^{-\frac{s}{2}}\eta^{-\frac{1}{6}}(y).
\ee
Then, it holds that
\begin{equation}
\bes
||F_{U_i,(3)}^{(m-1)}||_{L^2}&\les e^{\frac{s}{2}}(H^2||\pa^{\ga}\mathcal{P}||_{L^2}+\kappa_0^2||\pa^{\ga}H||_{L^2})^a(H^2||\pa^{\ga}\mathcal{P}||_{L^2})^b\cdot||D^2(\mathcal{P}H^2JN_i)||_{L^q}^{1-a}||D^2\mathcal{P}||_{L^q}^{1-b}\\
&\leq e^{\frac{a+b-1}{2}s} E_m^{a+b}M^{1-(a+b)},
\end{split}
\end{equation}
due to Lemma\ref{sobolev2}. Therefore,
\bee\label{L2FUim-13}
2\sum_{|\ga|=m}\lam^{|\check{\ga}|}\int|F_{U_i,(3)}^{(m-1)}\pa^{\ga}U_i|\leq e^{\frac{a+b-1}{2}s}E_m^{1+a+b}M^{1-(a+b)}\leq \da E_m^2+e^{-s}M^{2m-1}.
\ee
In conclusion, collecting \eqref{L2FUim}, \eqref{L2FUim-11} and \eqref{L2FUim-13} yields
\bee
2\sum_{|\ga|=m}\lam^{|\check{\ga}|}\int\left|F_{U_i}^{(\ga)}\pa^{\ga}U_i\right|\leq (5+C\da)E_m^2+e^{-s}M^{4m-1}.
\ee
\end{pf}
Then, there exists a universal constant $C$ such that
\bee
\dfrac{\pa}{\pa s}E_m^2(s)+(|\ga|-C)E_m^2\leq e^{-s}M^{4m-1}.
\ee
due to Proposition\ref{energyinequality} and Lemma\ref{forcingtermsestimates}. Then, by taking $|\ga|$ large enough and Gronwall inequality, one obtains
\begin{equation}\label{energyinequality2}
\bes
E_m^2(s)&\leq e^{-2(s-s_0)}E_m^2(s_0)+3e^{-s}M^{4m-1}(1-e^{-(s-s_0)})\\
&\leq 2\kappa_0^2\ep^{-1}e^{-2s}+3e^{-s}M^{4m-1}(1-\ep^{-1}e^{-s}).
\end{split}
\end{equation}
Therefore, the $\dot{H}^m$ bounds for $(W,Z,A_v,K)$ is the direct consequence of \eqref{energyinequality2} and the following lemma.
\begin{lem}
For $\ep$ sufficiently small in terms of $\kappa_0,M,m$, the following bounds hold
\bee
\lam^m\left(||U||_{\dot{H}^m}^2+||S||_{\dot{H}^m}^2+||K||_{\dot{H}^m}^2-e^{-2s}\right)\leq E_m^2\leq\kappa_0^2\left(||U||_{\dot{H}^m}^2+||S||_{\dot{H}^m}^2+||K||_{\dot{H}^m}^2+e^{-2s}\right)
\ee
for all $a\geq-\log\ep$. As a consequence,
\bee
\kappa_0^{-2}E_m^2-e^{-2s}\leq e^{-s}||W||_{\dot{H}^m}^2+||Z,A_v,K||_{\dot{H}^m}^2\leq 4\lam^{-m}E_m^2+4e^{-2s}.
\ee
\end{lem}
\begin{pf}
To prove this lemma, due to the estimates $|\mathcal{P}-\frac{\kappa_0}{2}|\leq\ep^{\frac{1}{8}},$ $|H-1|\leq\ep$ and the equivalence of $E_m^2$ and $\dot{H}^m$ norm of $(U,\mathcal{P},H)$, it suffices to prove that there exists a universal constant $C$ and a small $\da$ such that
\begin{align}
||\pa^{\ga}H-\pa^{\ga}K||_{L^2}&\leq C\ep||\pa^{\ga}K||_{L^2}+e^{-s},\label{knorm}\\
||\pa^{\ga}S-H\pa^{\ga}\mathcal{P}-\mathcal{P}\pa^{\ga}H||_{L^2}&\leq\da(||H||_{\dot{H}^m}+||\mathcal{P}||_{\dot{H}^m})+\da^{-\frac{2m-7}{2}}e^{-\frac{2m-3}{4}s}.\label{snorm}
\end{align}
For \eqref{knorm}, it can be proved by induction. When $|r|=1$, $||\pa H-\pa K||_{L^2}=|H-1|||\pa K||_{L^2}\leq C\ep||\pa K||_{L^2}$. Assume \eqref{knorm} hold for $|\ga|=m-1$. Then, for $|\ga|=m$, due to $\pa^{\ga}H-\pa^{\ga}K=\pa^{\ga-1}(H\pa K)-\pa^{\ga}K$, it holds that
\begin{align*}
||\pa^{\ga}H-\pa^{\ga}K||_{L^2}&\leq\left(||\pa^{\ga-1}K||_{L^2}||\pa^{\ga}K||_{L^{\infty}}+||H-1||_{L^{\infty}}||\pa^{\ga}K||_{L^2}\right)\\
&\leq\left(||\pa^{\ga-1}K||_{L^2}+e^{-s}\right)||\pa K||_{L^{\infty}}+C\ep||\pa^{\ga}K||_{L^2}\\
&\leq||\pa^{\ga}K||_{L^2}^{\frac{2m-5}{2m-3}}\ep^{\frac{1}{2}}e^{-\frac{s}{2}}+C\ep||\pa^{\ga}K||_{L^2}+e^{-s}\\
&\leq C\ep||\pa^{\ga}K||_{L^2}+e^{-s}.
\end{align*}
For \eqref{snorm}, it follows from Moser inequality that
\begin{align*}
||\pa^{\ga}S-H\pa^{\ga}\mathcal{P}-\mathcal{P}\pa^{\ga}H||_{L^2}&\leq C\left(||\nabla H||_{L^{\infty}}||\mathcal{P}||_{\dot{H}^m}+||\nabla\mathcal{P}||_{L^{\infty}}||H||_{\dot{H}^{m-1}}
\right)\\
&\leq C\left(||\nabla H||_{L^{\infty}}||\mathcal{P}||_{\dot{H}^m}^{\frac{2m-5}{2m-3}}||\mathcal{P}||_{L^{\infty}}^{\frac{2}{2m-3}}+
||\nabla\mathcal{P}||_{L^{\infty}}||H||_{\dot{H}^m}^{\frac{2m-5}{2m-3}}||H||_{L^{\infty}}^{\frac{2}{2m-3}}\right)\\
&\leq\da(||\mathcal{P}||_{\dot{H}^m}+||H||_{\dot{H}^m})+\da^{-\frac{2m-7}{2}}e^{-\frac{2m-3}{4}s}.
\end{align*}
\end{pf}
\\
\bibliographystyle{plain}
\bibliography{ref}
\end{document}